\documentclass[11pt,twoside,leqno,a4paper]{article}

\usepackage{amsmath, amssymb, amsthm, booktabs, caption, cite, geometry, graphicx, indentfirst, mathtools, paralist, tikz-cd, xcolor}

\usepackage[nottoc,notlot,notlof]{tocbibind}

\usepackage{journal-layouts} 

\JournalLayout{Preprint} 

\definecolor{citeblue}{RGB}{0,0,128}
\definecolor{linkred}{RGB}{128,0,0}
\definecolor{urlblue}{RGB}{0,0,204}

\usepackage{titlesec}
\usepackage[colorlinks=true, citecolor=citeblue, linkcolor=linkred, urlcolor=urlblue, pagebackref=true]{hyperref}

\titleformat{\section}{\filcenter\Large\bfseries}{\normalfont\thesection.}{0.5em}{}
\titleformat{\subsection}{\filcenter\bfseries}{\normalfont\thesubsection.}{.5em}{}

\allowdisplaybreaks 

\newtheorem{mythm}{\,\,\,\,\,\, \normalfont\scshape Theorem}[section]
\newtheorem{mycor}[mythm]{\,\,\,\,\,\, \normalfont\scshape Corollary}
\newtheorem{mylem}[mythm]{\,\,\,\,\,\, \normalfont\scshape Lemma}
\newtheorem{myprop}[mythm]{\,\,\,\,\,\, \normalfont\scshape Proposition}
\newtheorem{mydefn}[mythm]{\,\,\,\,\,\, \normalfont\it Definition}
\newtheorem{myrem}[mythm]{\,\,\,\,\,\, \normalfont\it Remark}
\newtheorem{myassumption}{\,\,\,\,\,\, \normalfont\it Assumption}

\newtheorem{theoremA}{\,\,\,\,\,\, \normalfont\scshape Theorem}

\newtheorem{theoremB}{\,\,\,\,\,\, \normalfont\scshape Theorem}

\newtheorem{theoremBB}{\,\,\,\,\,\, \normalfont\scshape Theorem}

\newtheorem{theoremC}{\,\,\,\,\,\, \normalfont\scshape Theorem}

\newcommand{\pd}[2]{\frac{\partial #1}{\partial #2}} 

\numberwithin{equation}{section}

\DeclareMathOperator{\sgn}{sgn}

\usepackage{etoolbox}
\makeatletter
\patchcmd{\endmyassumption}{\@endpefalse}{}{}{}
\makeatother

\renewcommand*{\backref}[1]{}
\renewcommand*{\backrefalt}[4]{%
  \ifcase #1
  \or
    {\scriptsize (Cited on page~#2.)}%
  \else
    {\scriptsize (Cited on pages~#2.)}%
  \fi}

\usepackage{tocloft} 

\makeatletter

\renewcommand{\l@section}{\@dottedtocline{1}{0em}{2.3em}}
\renewcommand{\l@subsection}{\@dottedtocline{2}{1.5em}{3.2em}}
\renewcommand{\@dotsep}{10000}

\makeatother

\newcommand{\indown}{\rotatebox[origin=c]{-90}{$\in$}} 

\usepackage{authblk}

\AtEndDocument{%
\par
\bigskip
\begin{tabular}{@{}l@{}}%
\textsc{\small Sajjad Bakrani}\\
{\small Instituto de Ci\^encias Matem\'aticas e Computa\c{c}\~ao (ICMC),}\\ {\small Universidade de S\~ao Paulo, S\~ao Carlos, Brazil}\\
\textit{\small E-mail}: \href{mailto:sajjad.bakrani@icmc.usp.br}{\texttt{\small sajjad.bakrani@icmc.usp.br}}\\
\href{https://orcid.org/0000-0001-7814-0992}{\small https://orcid.org/0000-0001-7814-0992}
\end{tabular}}

\newcommand\blfootnote[1]{%
\begingroup
\renewcommand\thefootnote{}\footnote{#1}%
\addtocounter{footnote}{-1}%
\endgroup
}

\begin{document}

\title{\bf Spatial switching near a four-loop homoclinic bouquet: a mechanism for multi-bump bound states in coupled nonlinear Schr\"odinger systems}


\author{\textsc{Sajjad Bakrani}}

\maketitle

\begin{abstract}
We study how a localized standing wave of a coupled nonlinear Schr\"odinger (CNLS) system can switch between components from one pulse to the next. The standing wave reduction of the CNLS system yields a four-dimensional Hamiltonian system with two reflection symmetries, in which four elementary pulses correspond to a bouquet of homoclinic loops to a saddle. For saddle eigenvalues $\pm\lambda_{1},\pm\lambda_{2}$ satisfying $0<\lambda_{1} < \lambda_{2} < 2\lambda_1$, and under suitable nondegeneracy assumptions, we construct stable and unstable invariant manifolds of the four-loop bouquet. Apart from trajectories asymptotic to the saddle in the corresponding time direction, these manifolds characterize one-sided trapping near the bouquet within the critical Hamiltonian level, namely the level set containing the saddle equilibrium. Both manifolds have codimension one in this level. We prove that if they intersect transversely away from the bouquet, producing super-homoclinic orbits, then there exist infinitely many multi-bump homoclinic orbits to the saddle. These orbits correspond to multi-bump localized standing waves of the CNLS system.
\end{abstract}



\tableofcontents



\section{Introduction}\label{sec:introduction}

Consider first the cubic coupled nonlinear Schr\"odinger (CNLS) system
\begin{equation}\label{eq:cubic_CNLS}
\begin{aligned}
i\Psi_t+\Psi_{xx}+(\mu_1|\Psi|^2+\beta|\Phi|^2)\Psi&=0,\\
i\Phi_t+\Phi_{xx}+(\beta|\Psi|^2+\mu_2|\Phi|^2)\Phi&=0,
\end{aligned}
\end{equation}
where $\mu_1,\mu_2>0$ and $\beta\in\mathbb R$. Its focusing integrable case $\mu_1=\mu_2=\beta$ is the Manakov system \cite{Manakov1974}.
Coupled Schr\"odinger equations arise, for example, in models of optical polarizations \cite{ChristodoulidesJoseph1988} and of multicomponent condensates \cite{EsryGreeneBurkeBohn1997}. We seek bound states of the form
\begin{equation}\label{eq:cubic_standing_wave_ansatz}
(\Psi(x,t),\Phi(x,t))
=\bigl(e^{i\omega_1t}\psi(x),e^{i\omega_2t}\phi(x)\bigr),
\qquad (\psi(x),\phi(x))\longrightarrow(0,0)
\quad\text{as }|x|\longrightarrow\infty,
\end{equation}
with real profiles and positive frequencies. For \eqref{eq:cubic_CNLS}, the four elementary profiles are $(\psi_\pm,0)$ and $(0,\phi_\pm)$, where
\begin{equation}\label{eq:profiles_cubic_CNLS}
\psi_\pm(x)=\pm\sqrt{2\omega_1/\mu_1}\,
\operatorname{sech}(\sqrt{\omega_1}x),\qquad
\phi_\pm(x)=\pm\sqrt{2\omega_2/\mu_2}\,
\operatorname{sech}(\sqrt{\omega_2}x).
\end{equation}
Changing the sign of an entire component amounts to a constant phase shift in the PDE. In the real profile equation the two signs give distinct orbits; recording them keeps track of the relative signs of pulses in the same component.

In this paper, we are interested in localized profiles of \eqref{eq:cubic_CNLS} consisting of several well-separated pulses, each resembling one of the four elementary profiles in \eqref{eq:profiles_cubic_CNLS}. Under the assumptions stated below, we construct stable and unstable manifolds of the union of the four profiles in the spatial dynamical system. We prove that if they intersect transversely within the critical Hamiltonian level away from the union itself, then infinitely many such multi-bump localized profiles exist. We also describe the possible sequences of components and signs of their successive pulses.

\paragraph{The spatial equation.}
Our results apply more generally to
\begin{equation*}
\begin{aligned}
i\Psi_t+\Psi_{xx}+f_1(|\Psi|^2,|\Phi|^2)\Psi&=0,\\
i\Phi_t+\Phi_{xx}+f_2(|\Psi|^2,|\Phi|^2)\Phi&=0,
\end{aligned}
\end{equation*}
where the smooth real nonlinearities vanish at the origin and have a smooth potential $F(s_{1}, s_{2})$ with $\partial_{s_j}F=f_j$ and $F(0,0)=0$. Substitution of \eqref{eq:cubic_standing_wave_ansatz} gives
\begin{equation*}
\begin{aligned}
\psi^{\prime \prime} - \omega_1\psi + f_1(\psi^2,\phi^2)\psi &=0,\\
\phi^{\prime \prime} - \omega_2\phi + f_2(\psi^2,\phi^2)\phi &=0.
\end{aligned}
\end{equation*}
Defining $(q_1,q_2) := (\psi,\phi)$ and $(p_1,p_2) := (\psi',\phi')$, this is the Hamiltonian system
\begin{equation}\label{eq:CNLS_spatial_Hamiltonian_system}
q_j^{\prime} = \frac{\partial H}{\partial p_j},\qquad
p_j^{\prime} = -\frac{\partial H}{\partial q_j},\qquad j=1,2,
\end{equation}
with
\begin{equation*}
H = \tfrac{1}{2} (p_1^2 + p_2^2) - \tfrac{1}{2}(\omega_{1} q_1^2 +\omega_{2} q_2^2) + \tfrac{1}{2} F(q_1^2, q_2^2).
\end{equation*}
A localized profile is a homoclinic trajectory to the equilibrium $O$ of this system and lies in $\{H=0\}$. Suppose that the two invariant $(q_{1}, p_{1})$- and $(q_{2}, p_{2})$- planes contain two primary homoclinic orbits each. We write
\begin{equation}\label{eq:homoclinics_of_CNLS}
B=\Gamma_{1+}\cup\Gamma_{1-}\cup\Gamma_{2+}\cup\Gamma_{2-}
\end{equation}
for their union, the \emph{homoclinic bouquet}. In the cubic example, following $\Gamma_{1+}$ produces a positive pulse in $\psi$, while following $\Gamma_{2-}$ produces a negative pulse in $\phi$. An orbit passing near these loops produces the corresponding excursions with the other component small. We use \emph{multi-bump} to refer to profiles formed by several such separated excursions.

The decay rates at $O$ are $\lambda_j=\sqrt{\omega_j}$. We consider
\begin{equation}\label{eq:spectral_window}
0<\lambda_1<\lambda_2<2\lambda_1,
\qquad\text{equivalently}\qquad
0<\omega_1<\omega_2<4\omega_1.
\end{equation}
Thus the first pair of loops approaches the saddle along its leading (slower-decaying) directions, and the second pair along its non-leading directions. The following statements summarize the results for the standing-wave problem. We assume that the primary loops are transverse in $\{H=0\}$ and that the global maps satisfy Assumption~\ref{assump:nonzero_coeffs} stated later. All stable and unstable manifolds below refer to evolution in the spatial variable $x$.

\begin{theoremA}[Informal: geometry of the bouquet]
The bouquet has two-dimensional local stable and unstable $C^1$ manifolds in $\{H=0\}$. Apart from trajectories asymptotic to $O$ in the corresponding direction, these manifolds describe exactly which trajectories stay near the bouquet for all increasing or all decreasing $x$. They are tangent to the corresponding manifolds of $O$ along $\Gamma_{1\pm}$ and transverse to them along $\Gamma_{2\pm}$.
\end{theoremA}

We can assign a code to each orbit on the bouquet manifolds that records the successive excursion along the four primary profiles.
\begin{theoremBB}[Informal: selection of pulse sequences]
On either local bouquet manifold, successive primary excursions alternate between the $\phi$ and $\psi$ components. On removing $B$, each local bouquet manifold has either two or four connected components. Each carries one cyclic code of primitive period two or four; the global maps, introduced later, determine the admissible codes, including the signs of the primary loops.
\end{theoremBB}

A periodic code specifies the order of visits; it implies neither a periodic profile nor equal spacing between successive pulses. An orbit that follows a code indefinitely in one direction makes infinitely many excursions and is not a localized standing wave. The next result explains how its itinerary can occur as a finite block in a localized profile.

\begin{theoremC}[Informal: localized multi-bump profiles]
Assume in addition that the global stable and unstable manifolds of $B$ intersect transversely along an orbit away from $B$. Such an orbit is called \emph{super-homoclinic}. For every sufficiently large pair of integers $(n,m)$, there is a bound-state profile that follows a finite block of the past itinerary of the super-homoclinic orbit, shadows its global excursion and then follows a finite block of its future itinerary. The numbers of returns on the two sides can be chosen independently.
\end{theoremC}

\begin{figure}[!htb]
\centering
\includegraphics[width=0.75\linewidth]{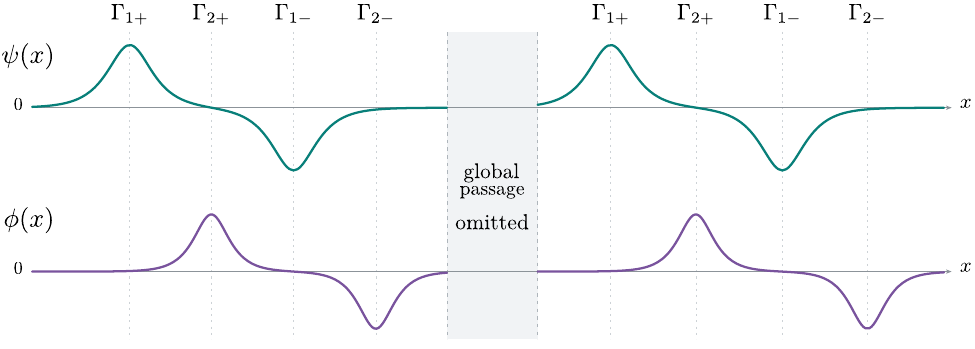}
\caption{Schematic localized CNLS profile shadowing a super-homoclinic passage. Excursions near $\Gamma_{1\pm}$ and $\Gamma_{2\pm}$ produce signed bumps predominantly in $\psi$ and $\phi$, respectively. On either side of the shaded interval, successive dominant bumps alternate between the two components. The shaded interval represents the global excursion, whose profile is omitted. Both components decay to zero as $x\to\pm\infty$.}
\label{fig:CNLS_schematic_profiles}
\end{figure}

Figure~\ref{fig:CNLS_schematic_profiles} shows how to read this result in the two profile components. The past and future itineraries of the chosen super-homoclinic orbit determine the sequences of components and signs in the corresponding blocks of pulses. The integers $n$ and $m$ count returns to specified cross-sections to the bouquet; Theorem~\ref{thm:superhomoclinic} gives the corresponding numbers of bumps.

\paragraph{Relation to earlier work.}
The four-loop problem grew out of work on self-localized electromagnetic fields~\cite{EleonskiiKulagin1989, EleonskyKorolevKulagin92} and on the geometry and coding of homoclinic orbits~\cite{EleonskyKorolev1993, EleonskyCodinghomoclinic1996}. Shilnikov and Turaev~\cite{ShilnikovDimaSuperhomoclinic1997} studied the four-loop bouquet when the two decay rates are equal, i.e. $\lambda_{1} = \lambda_{2}$. They constructed its stable and unstable manifolds, found alternating codes, and obtained multi-pulse homoclinic orbits from a transverse super-homoclinic orbit. In their setting, the bouquet manifolds are tangent to the corresponding manifolds of the equilibrium along all four loops. For the unequal rates in~\eqref{eq:spectral_window}, this tangency persists along the leading loops $\Gamma_{1\pm}$, whereas along the non-leading loops $\Gamma_{2\pm}$ the corresponding manifolds meet transversely in the critical level.

The case of a single and a symmetric pair of non-leading loops was also studied in \cite{Bakrani2022JDE}. With $\gamma=\lambda_1/\lambda_2$, two-dimensional local stable or unstable manifolds can occur in those settings for $\gamma<1/2$, subject to sign conditions on the global maps. However, for $1/2<\gamma<1$, the local stable and unstable sets are trivial.

The study of a four-loop bouquet differs substantially from the two-loop case of \cite{Bakrani2022JDE}. An orbit can leave a neighborhood of the non-leading pair, follow a leading loop, and return without ever leaving a neighborhood of the full four-loop bouquet. To account for these passages, we derive a new normal form and use it to estimate the local maps from each loop family to the other, together with their derivatives. Under our hypotheses, the full bouquet thus has both two-dimensional local manifolds in the range $1/2<\gamma<1$, in which the non-leading pair has neither. These results concern the critical level of the first integral; for the dynamics on nearby levels, see~\cite{DimaShilnikov1989, Turaev2014} for leading loops and~\cite{Bakrani2025JDE} for non-leading loops.

Super-homoclinic orbits and similar notions have also been studied in nonconservative systems~\cite{Homburg1996memoirs}, in systems with a first integral~\cite{DimaMultipulse2001}, near a homoclinic loop to a bifocus in a reversible system~\cite{Rodrigues2019}, and in heteroclinic networks \cite{ChawanyaAshwin2010}. Bandara et al.~\cite{Krauskopf2023generalized} proposed a mechanism for generating super-homoclinic orbits as limits of generalized solitons in a scalar nonlinear Schr\"odinger equation with quartic dispersion.

Multi-pulse solitary waves in coupled Schr\"odinger equations have been studied through tail matching and homoclinic bifurcation methods. Yang~\cite{Yang1998} used tail matching to construct solitary-wave trains. Other constructions exploit specific homoclinic bifurcations: Yew~\cite{Yew2001Multipulses} obtained multipulses at a resonant homoclinic bifurcation, Jackson~\cite{Jackson2006HomoclinicFlip} obtained alternating pulses at an orbit flip, and Champneys and Yang~\cite{ChampneysYang2002Nonlocal} identified a bifurcation at which pulses in one component separate from a central pulse in the other. These mechanisms belong to the broader theory of homoclinic bifurcations surveyed in~\cite{HomburgSandstede2010}. Our approach uses the invariant manifolds of the full four-loop bouquet to determine the admissible sequences of pulse components and signs. A transverse super-homoclinic orbit then gives rise to localized profiles that shadow its past and future itineraries for independently prescribed, sufficiently large numbers of section returns.

Beyond these construction mechanisms, a broader CNLS literature addresses the existence, classification, and stability of vector solitary waves. Variational methods establish the existence of least-energy states~\cite{Lin_Wei_ground_state_CNLS_2005, LinWei2008Erratum, Ambrosetti2007standing, Sirakov2007least} and multi-bump states~\cite{Ambrosetti_multi_bump_2007, Terracini_multiple_phases_2009}; see also~\cite{Byeon_standing_wave_NLS_2026} for attractive couplings. Bifurcation and numerical studies describe families of vector solitary waves~\cite{HaeltermanSheppard1994, ClassificationSolitaryWaves1997}, while experiments demonstrate the occurrence of multihump optical solitons~\cite{MitchellSegevChristodoulides1998}. Once such waves are found, their stability under the time evolution of the PDE is a separate question from the spatial classification developed here. This stability question has been studied for multihump waves in~\cite{OstrovskayaKivsharSkryabinFirth1999, PelinovskyYang2005Multihump}. Related work on nondegenerate vector solitons and breathers in the integrable Manakov system can also be found in~\cite{Pelinovsky_CNLS_breathers_2026}.

\paragraph{Acknowledgments and use of artificial intelligence.}
This work was supported by S\~ao Paulo Research Foundation (FAPESP), grant 2013/07375-0.

I developed the mathematical arguments in this paper and carried out the calculations by hand or with Maple (computer algebra). In particular, I used Maple for the longer symbolic calculations arising from Shilnikov's boundary value method in Appendix~\ref{sec:proofs_technical_lemmas}.

I used OpenAI's ChatGPT (6 Astra and 6.2 Astra) to help revise the exposition, prepare figures and \LaTeX{} and TikZ code, and locate relevant literature. I also used it to cross-check the calculations throughout the manuscript. I reviewed the suggestions and decided which to incorporate. I did not use AI output as a mathematical source, and I take full responsibility for the final manuscript, including its proofs, calculations, and references.

\section{Abstract homoclinic setting and main results}\label{sec:abstract_setting}

The standing-wave reduction in Section~\ref{sec:introduction} leads naturally to a four-dimensional spatial dynamical system. We now isolate the geometric structure used in the proofs and formulate the problem independently of its CNLS origin. The abstraction is useful in two ways.  It shows that the mechanism does not depend on the particular form of the nonlinearities in the CNLS system, and it exposes a homoclinic problem that is of independent interest: how can several homoclinic loops to a real saddle cooperate to organize trajectories that remain near the entire configuration?

To start with, we first put the spatial system \eqref{eq:CNLS_spatial_Hamiltonian_system} in coordinates adapted to the stable and unstable eigendirections at $O$. Consider the linear change of variables
\begin{equation*}
\begin{alignedat}{2}
u_{1} &= \frac{1}{2}\left(p_{1} - \sqrt{\omega_{1}} q_{1}\right),\qquad
u_{2} &&= \frac{1}{2}\left(\sqrt{\omega_{2}} q_{2} - p_{2}\right),\\
v_{1} &= \frac{1}{2}\left(\frac{p_{1}}{\sqrt{\omega_1}} + q_{1}\right),\qquad
v_{2} &&= \frac{1}{2}\left(q_{2} + \frac{p_{2}}{\sqrt{\omega_2}}\right).
\end{alignedat}
\end{equation*}
In these coordinates the linear part is diagonal, and the system becomes
\begin{equation*}
\begin{aligned}
\dot{u}_{1} &= -\sqrt{\omega_1} u_{1} - \frac{1}{2} f_{1}(q_{1}^{2}, q_{2}^{2}) q_{1},\\
\dot{u}_{2} &= -\sqrt{\omega_2} u_{2} + \frac{1}{2} f_{2}(q_{1}^{2}, q_{2}^{2}) q_{2},\\
\dot{v}_{1} &= +\sqrt{\omega_1} v_{1} - \frac{1}{2\sqrt{\omega_1}} f_{1}(q_{1}^{2}, q_{2}^{2}) q_{1},\\
\dot{v}_{2} &= +\sqrt{\omega_2} v_{2} - \frac{1}{2\sqrt{\omega_2}} f_{2}(q_{1}^{2}, q_{2}^{2}) q_{2},
\end{aligned}
\end{equation*}
where $q_{1} = v_{1} - \frac{u_1}{\sqrt{\omega_1}}$ and $q_{2} = v_{2} + \frac{u_2}{\sqrt{\omega_2}}$. Moreover, the Hamiltonian becomes the first integral
\begin{equation}\label{eq:first_integral_CNLS}
H(u_{1}, u_{2}, v_{1}, v_{2}) = 2\sqrt{\omega_{1}} u_{1} v_{1} - 2\sqrt{\omega_{2}} u_{2} v_{2} + \frac{1}{2} F\left((v_{1} - \frac{u_1}{\sqrt{\omega_1}})^{2}, (v_{2} + \frac{u_2}{\sqrt{\omega_2}})^{2}\right).
\end{equation}

Motivated by this system, we now pass to the abstract problem. It is important to note that the time variable $t$ is used for the abstract flow; in the standing-wave problem, the spatial variable $x$ is used. Consider a $C^{\infty}$ system of ODEs on $\mathbb{R}^{4}$ of the form
\begin{equation}\label{eq_system_1}
\begin{aligned}
\dot{u}_1 &= -\lambda_{1} u_{1} + \cdots,\\
\dot{u}_2 &= -\lambda_{2} u_{2} + \cdots,\\
\dot{v}_1 &= \lambda_{1} v_{1} + \cdots,\\
\dot{v}_2 &= \lambda_{2} v_{2} + \cdots,
\end{aligned}
\end{equation}
where $\lambda_{2} \geq \lambda_{1} > 0$ and the dots denote $C^{\infty}$ functions of $(u_{1}, u_{2}, v_{1}, v_{2})$ that, together with their first derivatives, vanish at the origin $O$. Thus $O$ is a saddle equilibrium with two stable and two unstable real eigendirections. We assume that the system admits a ${C}^{\infty}$ first integral of the form
\begin{equation}\label{eq87ywhouyeg3r23}
H = \gamma u_{1} v_{1} - u_{2} v_{2} + \cdots,
\end{equation}
where $\gamma = \lambda_{1} \lambda_{2}^{-1}$ and the dots stand for cubic and higher order terms. We further assume that system \eqref{eq_system_1} is invariant under the two $\mathbb{Z}_{2}$-symmetries
\begin{equation}\label{eq:symmetry_1}
(u_{1}, v_{1}) \leftrightarrow (-u_{1}, -v_{1})
\end{equation}
and
\begin{equation}\label{eq:symmetry_2}
(u_{2}, v_{2}) \leftrightarrow (-u_{2}, -v_{2}).
\end{equation}
We also assume that the first integral is invariant under both symmetries, that is,
\begin{equation*}
\begin{aligned}
H(u_{1}, u_{2}, v_{1}, v_{2}) &= H(-u_{1}, u_{2}, -v_{1}, v_{2}),\\
H(u_{1}, u_{2}, v_{1}, v_{2}) &= H(u_{1}, -u_{2}, v_{1}, -v_{2}).
\end{aligned}
\end{equation*}
Thus the remainder in \eqref{eq87ywhouyeg3r23} actually starts at order four. For the CNLS reduction, take $\lambda_j=\sqrt{\omega_j}$ and divide the first integral \eqref{eq:first_integral_CNLS} by $2\lambda_2$ to obtain this normalization\footnote{The forms~\eqref{eq_system_1} and \eqref{eq87ywhouyeg3r23} are not specific to the CNLS reduction. More generally, suppose that a system on $\mathbb{R}^{4}$ has a hyperbolic equilibrium $O$ and a smooth first integral whose quadratic part at $O$ is non-degenerate. The linear part of the integral vanishes at $O$. In the presence of symmetry~\eqref{eq:symmetry_1}, the equilibrium $O$ has to be a saddle, and suitable local coordinates put the system and the integral into the forms~\eqref{eq_system_1} and \eqref{eq87ywhouyeg3r23}, respectively. See, for example, \cite{BakraniPhDthesis,Turaev2014}.}.

Write $W^s_{\mathrm{glo}}(O)$ and $W^u_{\mathrm{glo}}(O)$ for the global stable and unstable manifolds of $O$. The fixed-point sets of the two symmetries are the coordinate planes $\{u_{2}=v_{2}=0\}$ and $\{u_{1}=v_{1}=0\}$, and both planes are invariant under the flow. Motivated by the four signed primary profiles in \eqref{eq:homoclinics_of_CNLS}, we impose the following hypothesis.

\begin{figure}[!htb]
\centering
\includegraphics[width=0.75\linewidth]{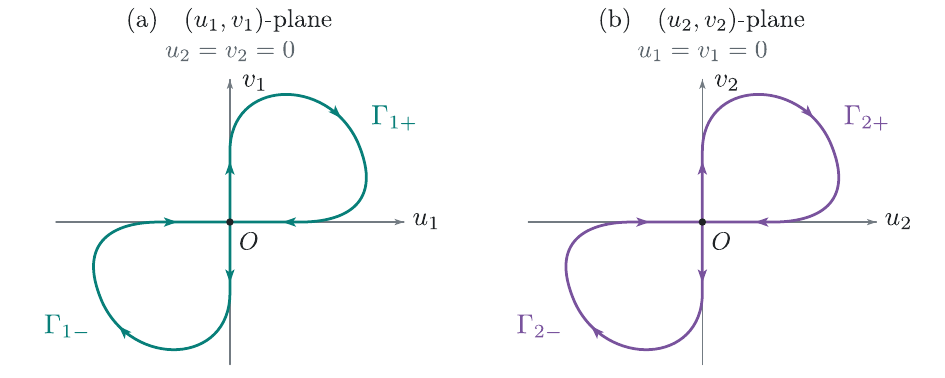}
\caption{The primary homoclinic orbits in the invariant planes. Each symmetry exchanges the two loops in its own plane while fixing the other plane.}
\label{fig:four_loops_in_invariant_planes}
\end{figure}

\begin{myassumption}\label{assump:exist_transverse_homoclinics}
The equilibrium $O$ has four primary homoclinic orbits, two in the $(u_1,v_1)$-plane and two in the $(u_2,v_2)$-plane. Along each orbit, the manifolds $W^s_{\mathrm{glo}}(O)$ and $W^u_{\mathrm{glo}}(O)$ intersect transversely within $\{H=0\}$; see Figure~\ref{fig:four_loops_in_invariant_planes}.
\end{myassumption}

This transversality condition is imposed in the three-dimensional energy level. The next assumption is the spectral condition \eqref{eq:spectral_window}.

\begin{myassumption}\label{assump:gamma_larger_half}
We assume that $2\lambda_{1} > \lambda_{2} > \lambda_{1} > 0$.
\end{myassumption}

\begin{figure}[!tbh]
\centering
\includegraphics[width=\linewidth]{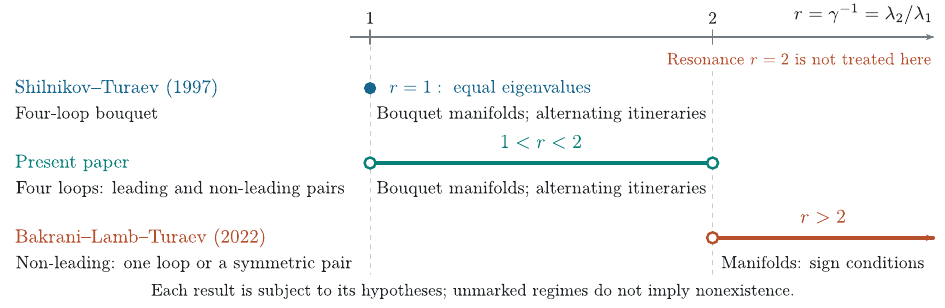}
\caption{Spectral regimes, with $r=\lambda_2/\lambda_1=\gamma^{-1}$. The rows refer to the equal-rate four-loop bouquet \cite{ShilnikovDimaSuperhomoclinic1997}, the four-loop bouquet studied here, and non-leading one- and two-loop configurations \cite{Bakrani2022JDE}. Each row is subject to its own hypotheses; unmarked intervals make no nonexistence claim.}
\label{fig:spectral_regimes}
\end{figure}

In particular, $\frac{1}{2} < \gamma < 1$. This assumption is not merely technical. Even for a single non-leading loop or a non-leading pair, the local
stable and unstable sets behave differently on the two sides of
$\gamma=1/2$; see~\cite{Bakrani2022JDE}. Figure~\ref{fig:spectral_regimes} places this assumption alongside the earlier results for an equal-rate four-loop bouquet and for a single non-leading loop or a non-leading pair. We do not treat the resonance $\lambda_2=2\lambda_1$.

To state the remaining hypotheses and the main results, we fix the local coordinates constructed in Section~\ref{sec:local_analysis_near_equilibrium}. Lemma~\ref{Lem_normal_form_1} puts system~\eqref{eq_system_1} into a normal form that preserves the symmetries and the first integral while straightening the local stable, unstable, strong stable, and strong unstable manifolds of $O$.

Denote the two loops in the $(u_i,v_i)$-plane by $\Gamma_{i+}$ and $\Gamma_{i-}$. We choose signs so that $\Gamma_{i+}$ leaves $O$ along the positive $v_i$-axis and returns along the positive $u_i$-axis; $\Gamma_{i-}$ uses the two negative axes. A sign change in the local coordinates preserves the normal form \eqref{eq:normal_form} but may change the relative sign of the quadratic terms of $H$; we have
\begin{equation}\label{first_integral_kappa}
H = \gamma u_{1} v_{1} - \kappa u_{2} v_{2} + \cdots
\end{equation}
with $\kappa=\pm1$; see~\eqref{eq:first_integral}. For each $i$, the two
loops $\Gamma_{i+}$ and $\Gamma_{i-}$ are exchanged by the corresponding
$\mathbb{Z}_{2}$-symmetry.

\begin{figure}[!htb]
\centering
\includegraphics[width=0.85\linewidth]{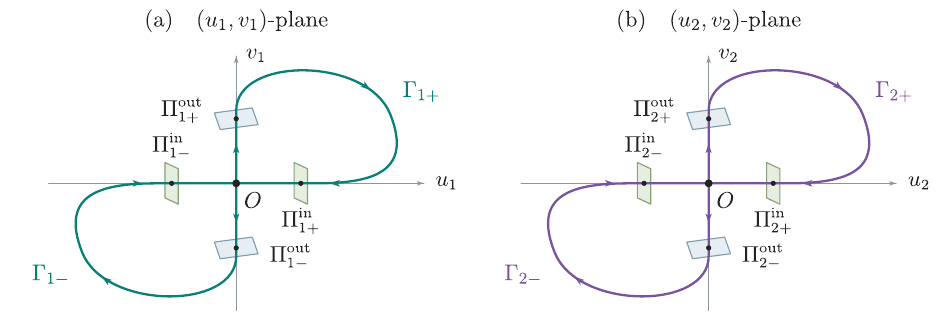}
\caption{The eight cross-sections on the homoclinic orbits. Each $\Gamma_{i\sigma}$ leaves $O$ through $\Pi_{i\sigma}^{\mathrm{out}}$ and returns through $\Pi_{i\sigma}^{\mathrm{in}}$. The shaded patches schematically represent two-dimensional transverse sections in $\{H=0\}$. We choose local coordinates on each section, and the dots mark the coordinate origins.}
\label{fig:eight_cross_sections}
\end{figure}

Fix $\delta>0$ sufficiently small. For each loop $\Gamma_{i\sigma}$, where $i=1,2$ and $\sigma=\pm$, introduce the incoming and outgoing cross-sections
\begin{equation}\label{eq:cross-sections}
\Pi_{i\sigma}^{\mathrm{in}} = \{u_{i} = \sigma\delta\} \cap \{H=0\} \quad \mathrm{and}\quad \Pi_{i\sigma}^{\mathrm{out}} = \{v_{i} = \sigma\delta\} \cap \{H=0\}
\end{equation}
which are transverse to the corresponding homoclinic orbit; see Figure~\ref{fig:eight_cross_sections}. By the implicit function theorem, $(u_{1},v_{1})$ are local coordinates on $\Pi^{\mathrm{in}}_{2\sigma}$ and $\Pi^{\mathrm{out}}_{2\sigma}$, while $(u_{2},v_{2})$ are local coordinates on $\Pi^{\mathrm{in}}_{1\sigma}$ and $\Pi^{\mathrm{out}}_{1\sigma}$; see Lemma~\ref{Lem_coordinates_on_cross_sections}. In these coordinates the intersection point of $\Gamma_{i\sigma}$ with either of its cross-sections is the section origin $(0,0)$. Following the flow from $\Pi^{\mathrm{out}}_{i\sigma}$ along the global part of $\Gamma_{i\sigma}$ to $\Pi^{\mathrm{in}}_{i\sigma}$ defines a global Poincar\'e map $T^{\mathrm{glo}}_{i\sigma}$. Symmetry makes the maps for the two signs identical, so we denote both by $T^{\mathrm{glo}}_{i}$. Taking the symmetries into account, the Taylor expansions of these maps at the section origin $(0,0)$ are given by
\begin{equation}\label{eq_Taylor_expansion_global_map}
\begin{aligned}
T^{\mathrm{glo}}_{1}\begin{pmatrix}
u_{2}\\
v_{2}
\end{pmatrix} &= \begin{pmatrix}
a_{1} u_{2} + b_{1} v_{2} + O(\|(u_{2}, v_{2})\|^{3})\\[3pt]
c_{1} u_{2} + d_{1} v_{2} + O(\|(u_{2}, v_{2})\|^{3})
\end{pmatrix},\\[5pt]
T^{\mathrm{glo}}_{2}\begin{pmatrix}
u_{1}\\
v_{1}
\end{pmatrix} &= \begin{pmatrix}
a_{2} u_{1} + b_{2} v_{1} + O(\|(u_{1}, v_{1})\|^{3})\\[3pt]
c_{2} u_{1} + d_{2} v_{1} + O(\|(u_{1}, v_{1})\|^{3})
\end{pmatrix}
\end{aligned}
\end{equation}
for some constants $a_{i}$, $b_{i}$, $c_{i}$ and $d_{i}$ ($i=1,2$).

\begin{myassumption}\label{assump:nonzero_coeffs}
We assume that all the coefficients $a_{i}$, $b_{i}$, $c_{i}$, and $d_{i}$ for $i = 1,2$ are nonzero.
\end{myassumption}

The transversality required in Assumption~\ref{assump:exist_transverse_homoclinics} follows from Assumption~\ref{assump:nonzero_coeffs}. Once the primary loops exist, their transversality within the critical level is equivalent to $d_i\ne0$ for $i=1,2$. Indeed, the local unstable manifold of $O$ intersects the outgoing section along the vertical axis. Its image under the global map is a curve whose tangent slope at the incoming-section origin is $d_i/b_i$. Since the local stable manifold of $O$ intersects the incoming section along the horizontal axis, the two curves meet transversely precisely when $d_i\ne0$.

\begin{theoremA}[Invariant manifolds and one-sided trapped dynamics]
\label{thm:dynamic_near_loops}
Let $B=\Gamma_{1-}\cup\Gamma_{1+}\cup\Gamma_{2-}\cup\Gamma_{2+}$. Suppose that Assumptions~\ref{assump:exist_transverse_homoclinics},
\ref{assump:gamma_larger_half} and~\ref{assump:nonzero_coeffs} hold. Then $B$ has two-dimensional local stable and unstable invariant manifolds $W^s_{\mathrm{loc}}(B)$ and $W^u_{\mathrm{loc}}(B)$ in $\{H=0\}$, of class $C^1$. Points on the stable manifold approach $B$ in forward time, and points on the unstable manifold approach $B$ in backward time. Moreover:
\begin{enumerate}[(i)]
\item Both manifolds contain $B$. Along $\Gamma_{1+}\cup\Gamma_{1-}$, the manifold $W^s_{\mathrm{loc}}(B)$ is tangent to $W^s_{\mathrm{glo}}(O)$, and $W^u_{\mathrm{loc}}(B)$ is tangent to $W^u_{\mathrm{glo}}(O)$. Along $\Gamma_{2+}\cup\Gamma_{2-}$, both corresponding intersections
are transverse within $\{H=0\}$.
\item There is an open neighborhood $U$ of $B\cup\{O\}$ in $\{H=0\}$ with the following property. For $p\in U\setminus W^s_{\mathrm{glo}}(O)$, the forward semiorbit of $p$ stays in $U$ if and only if $p\in W^s_{\mathrm{loc}}(B)$. For $p\in U\setminus W^u_{\mathrm{glo}}(O)$, the backward semiorbit of $p$ stays in $U$ if and only if $p\in W^u_{\mathrm{loc}}(B)$.
\end{enumerate}
\end{theoremA}

\begin{figure}[!tbh]
\centering
\includegraphics[width=0.75\linewidth]{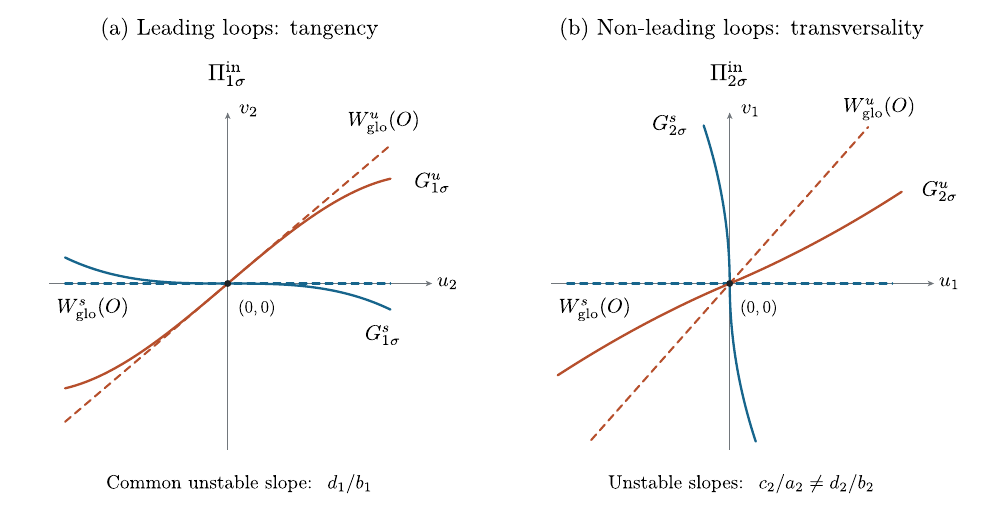}
\caption{Local section geometry of Theorem~\ref{thm:dynamic_near_loops}, for $\sigma=\pm$. Solid curves are $G^{s,u}_{i\sigma} =W^{s,u}_{\mathrm{loc}}(B)\cap\Pi^{\mathrm{in}}_{i\sigma}$; dashed curves are the local traces of $W^{s,u}_{\mathrm{glo}}(O)$ on the incoming sections. Both traces pass through the section origin $(0,0)$.
\textup{(a)} Along the leading loops, the corresponding stable and unstable curves are tangent: the stable tangent is horizontal and the common unstable slope is $d_1/b_1$.
\textup{(b)} Along the non-leading loops, both corresponding pairs are transverse. The stable bouquet curve has a vertical tangent while the stable equilibrium curve is the horizontal axis. The unstable tangent slopes are $c_2/a_2$ and $d_2/b_2$, respectively ($\tfrac{c_{2}}{a_{2}} \neq \tfrac{d_2}{b_2}$ as $T_{2}^{\mathrm{glo}}$ is a diffeomorphism).}
\label{fig:theorem_A_sections}
\end{figure}

Figure~\ref{fig:theorem_A_sections} shows these intersections on the cross-sections. We next describe the connected components of $W^u_{\mathrm{loc}}(B)\setminus B$ and $W^s_{\mathrm{loc}}(B)\setminus B$ by the successive section branches that an orbit visits.

Consider first the unstable manifold of the bouquet $B$. For each $i=1,2$ and $\sigma=\pm$, the intersection $W^{u}_{\mathrm{loc}}(B)\cap\Pi^{\mathrm{in}}_{i\sigma}$ is a curve through the point $(0,0)$ (the intersection of $\Gamma_{i\sigma}$ and $\Pi^{\mathrm{in}}_{i\sigma}$). Removing that point leaves two branches, lying in the upper and lower halves of the cross-section; we label them ${i}_{\sigma}^{+}$ and ${i}_{\sigma}^{-}$, respectively.

Choose an orbit in a connected component of $W^{u}_{\mathrm{loc}}(B)\setminus B$ and record its successive labeled incoming-section intersections in forward time, for as long as these local intersections are defined. Section~\ref{sec:Dynamics_near_four_loops} shows that these labels follow the periodic symbolic pattern generated by the forward transition rule $S_u$ of Definition~\ref{dfn_evolution_function_Su}. The actual forward itinerary may be finite. We call the smallest repeating block of this symbolic pattern the \emph{code} of the component, considered up to cyclic permutation. This code is independent of the chosen orbit within the component. The following theorem lists all possibilities.

\begin{theoremB}[Unstable switching codes]\label{thm:coding_unstable}
Let $k_u$ be the number of connected components of $W^u_{\mathrm{loc}}(B)\setminus B$, and let $\kappa$ be as in \eqref{first_integral_kappa}. Under the hypotheses of Theorem~\ref{thm:dynamic_near_loops}, the component counts and codes are given in Table \ref{tab:unstable_component_codes}.
\end{theoremB}
\begin{table}[!htb]
\centering
\begingroup
\normalfont\normalsize
\setlength{\tabcolsep}{8pt}
\renewcommand{\arraystretch}{1.08}
\newcommand{\UnstableCodePair}[2]{%
  $#1,\quad #2$%
}
\noindent\makebox[\linewidth][c]{%
\begin{tabular}{@{}cccl@{}}
\toprule
$c_{2}$ & $\kappa d_{1}a_{2}$ & $k_u$ & Component codes \\
\midrule
$>0$ & $>0$ & $4$ &
  \UnstableCodePair
    {{1}_{+}^{+}\, {2}_{+}^{+},\quad {1}_{+}^{-}\, {2}_{-}^{+}}
    {{1}_{-}^{+}\, {2}_{+}^{-},\quad {1}_{-}^{-}\, {2}_{-}^{-}} \\
\addlinespace[4pt]
$>0$ & $<0$ & $2$ &
  \UnstableCodePair
    {{1}_{+}^{+}\, {2}_{+}^{+}\, {1}_{+}^{-}\, {2}_{-}^{+}}
    {{1}_{-}^{+}\, {2}_{+}^{-}\, {1}_{-}^{-}\, {2}_{-}^{-}} \\
\addlinespace[4pt]
$<0$ & $>0$ & $2$ &
  \UnstableCodePair
    {{1}_{+}^{+}\, {2}_{+}^{-}\, {1}_{-}^{-}\, {2}_{-}^{+}}
    {{1}_{+}^{-}\, {2}_{-}^{-}\, {1}_{-}^{+}\, {2}_{+}^{+}} \\
\addlinespace[4pt]
$<0$ & $<0$ & $2$ &
  \UnstableCodePair
    {{1}_{+}^{+}\, {2}_{+}^{-}\, {1}_{-}^{+}\, {2}_{+}^{+}}
    {{1}_{+}^{-}\, {2}_{-}^{-}\, {1}_{-}^{-}\, {2}_{-}^{+}} \\
\bottomrule
\end{tabular}%
}
\par
\endgroup
\caption{Component counts and codes for the unstable bouquet
manifold in Theorem~\ref{thm:coding_unstable}.}
\label{tab:unstable_component_codes}
\end{table}

For example, the code $1_{+}^{+}\,2_{+}^{+}\,1_{+}^{-}\,2_{-}^{+}$ gives, in order, a positive $\psi$-bump, a positive $\phi$-bump, a positive $\psi$-bump, and a negative $\phi$-bump in the CNLS system. The two visits to $\Pi^{\mathrm{in}}_{1+}$ are on opposite branches, but both correspond to a positive $\psi$-bump.

\begin{figure}[!htb]
\centering
\includegraphics[width=0.85\textwidth]{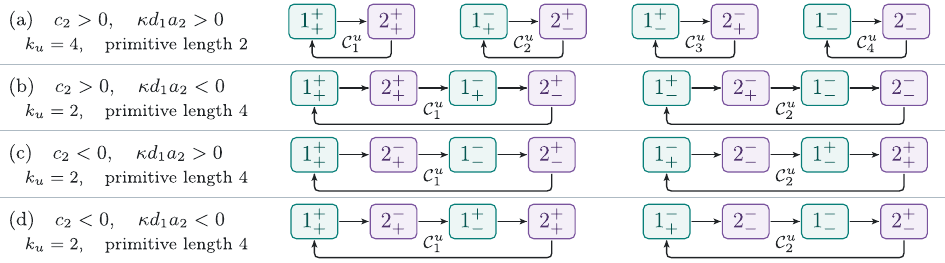}
\caption{Symbolic cycles of the unstable bouquet manifold in Theorem \ref{thm:coding_unstable} (Table \ref{tab:unstable_component_codes}). Arrows follow the forward rule $S_u$ of Table~\ref{table:bijection_Su}. Each cycle represents one connected component $\mathcal C_j^u$ of $W^u_{\mathrm{loc}}(B)\setminus B$; components are numbered within each panel.}
\label{fig:symbolic_cycle_components}
\end{figure}

We also read the stable code in forward time. On $\Pi^{\mathrm{in}}_{1\sigma}$ the labels $1_{\sigma}^{+}$ and $1_{\sigma}^{-}$ denote the right and left branches, respectively; on $\Pi^{\mathrm{in}}_{2\sigma}$ the labels $2_{\sigma}^{+}$ and $2_{\sigma}^{-}$ denote the upper and lower branches. As before, the code is the smallest repeating block, up to cyclic permutation.

\begin{theoremB}[Stable switching codes]\label{thm:coding_stable}
Let $k_s$ be the number of connected components of $W^s_{\mathrm{loc}}(B)\setminus B$. Under the hypotheses of Theorem~\ref{thm:dynamic_near_loops}, the component counts and codes are given in Table~\ref{tab:stable_component_codes}.
\end{theoremB}
\begin{table}[!htbp]
\centering
\begingroup
\normalfont\normalsize
\setlength{\tabcolsep}{6pt}
\renewcommand{\arraystretch}{1.08}
\noindent\makebox[\linewidth][c]{%
\begin{tabular}{@{}ccccl@{}}
\toprule
$d_{1}$ & $b_{2}$ & $\kappa a_{2}$ & $k_s$ & Component codes \\
\midrule
$>0$ & $<0$ & $>0$ & $4$ &
${1}_{+}^{+}\, {2}_{+}^{+},\quad
{1}_{+}^{-}\, {2}_{-}^{+},\quad
{1}_{-}^{+}\, {2}_{+}^{-},\quad
{1}_{-}^{-}\, {2}_{-}^{-}$ \\
\addlinespace[4pt]
$>0$ & $<0$ & $<0$ & $2$ &
${1}_{+}^{+}\, {2}_{-}^{+}\, {1}_{+}^{-}\, {2}_{+}^{+},\quad
{1}_{-}^{+}\, {2}_{-}^{-}\, {1}_{-}^{-}\, {2}_{+}^{-}$ \\
\addlinespace[4pt]
$>0$ & $>0$ & $>0$ & $2$ &
${1}_{+}^{+}\, {2}_{-}^{-}\, {1}_{-}^{-}\, {2}_{+}^{+},\quad
{1}_{+}^{-}\, {2}_{+}^{-}\, {1}_{-}^{+}\, {2}_{-}^{+}$ \\
\addlinespace[4pt]
$>0$ & $>0$ & $<0$ & $2$ &
${1}_{+}^{+}\, {2}_{+}^{-}\, {1}_{-}^{+}\, {2}_{+}^{+},\quad
{1}_{+}^{-}\, {2}_{-}^{-}\, {1}_{-}^{-}\, {2}_{-}^{+}$ \\
\midrule
$<0$ & $<0$ & $<0$ & $4$ &
${1}_{+}^{+}\, {2}_{-}^{+},\quad
{1}_{+}^{-}\, {2}_{+}^{+},\quad
{1}_{-}^{+}\, {2}_{-}^{-},\quad
{1}_{-}^{-}\, {2}_{+}^{-}$ \\
\addlinespace[4pt]
$<0$ & $<0$ & $>0$ & $2$ &
${1}_{+}^{+}\, {2}_{+}^{+}\, {1}_{+}^{-}\, {2}_{-}^{+},\quad
{1}_{-}^{+}\, {2}_{+}^{-}\, {1}_{-}^{-}\, {2}_{-}^{-}$ \\
\addlinespace[4pt]
$<0$ & $>0$ & $<0$ & $2$ &
${1}_{+}^{+}\, {2}_{+}^{-}\, {1}_{-}^{-}\, {2}_{-}^{+},\quad
{1}_{+}^{-}\, {2}_{-}^{-}\, {1}_{-}^{+}\, {2}_{+}^{+}$ \\
\addlinespace[4pt]
$<0$ & $>0$ & $>0$ & $2$ &
${1}_{+}^{+}\, {2}_{-}^{-}\, {1}_{-}^{+}\, {2}_{-}^{+},\quad
{1}_{+}^{-}\, {2}_{+}^{-}\, {1}_{-}^{-}\, {2}_{+}^{+}$ \\
\bottomrule
\end{tabular}%
}
\par
\endgroup
\caption{Component counts and codes for the stable bouquet manifold in Theorem~\ref{thm:coding_stable}.}
\label{tab:stable_component_codes}
\end{table}

\begin{figure}[!htb]
\centering
\includegraphics[width=0.85\textwidth]{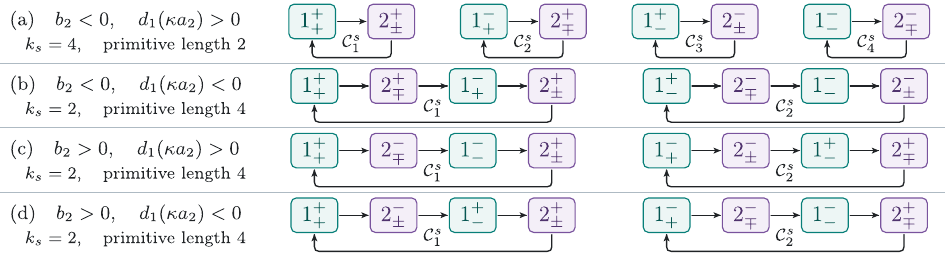}
\caption{Symbolic cycles of the stable bouquet manifold in Theorem~\ref{thm:coding_stable} (Table~\ref{tab:stable_component_codes}). Read the upper sign in every paired subscript for $d_1>0$ and the lower sign for $d_1<0$. Arrows follow the forward rule $S_s$ of Table~\ref{table:bijection_Ss}. Each cycle represents one connected component $\mathcal C_j^s$ of $W^s_{\mathrm{loc}}(B)\setminus B$; components are numbered within each panel.}
\label{fig:stable_symbolic_cycle_components}
\end{figure}

Following Theorems~\ref{thm:coding_unstable} and~\ref{thm:coding_stable}, every admissible itinerary alternates between the two invariant planes. An orbit on either bouquet manifold that shadows $\Gamma_{1+}$ or $\Gamma_{1-}$ must next shadow one of $\Gamma_{2+}$ or $\Gamma_{2-}$, and conversely. Hence no such orbit can remain indefinitely near only the leading pair or only the non-leading pair. Figures~\ref{fig:symbolic_cycle_components} and~\ref{fig:stable_symbolic_cycle_components} illustrate the cycles listed in Theorems~\ref{thm:coding_unstable} and~\ref{thm:coding_stable}, respectively.

To obtain multi-bump homoclinic orbits to $O$ that follow the codes prescribed by Theorems~\ref{thm:coding_unstable} and~\ref{thm:coding_stable}, we extend $W^s_{\mathrm{loc}}(B)$ and $W^u_{\mathrm{loc}}(B)$ under the flow. We denote the resulting global manifolds by $W^s_{\mathrm{glo}}(B)$ and $W^u_{\mathrm{glo}}(B)$. Only the next result requires the following additional hypothesis.

\begin{myassumption}\label{assum:exist_superhomoclinic}
The global stable and unstable manifolds of $B$ intersect transversely in $\{H=0\}$ at a point outside $B$.
\end{myassumption}

We do not verify this assumption for a concrete CNLS system here. However, the fact that $W^{s}_{\mathrm{glo}}(B)$ and $W^{u}_{\mathrm{glo}}(B)$ are codimension-1 inside $\{H=0\}$ makes a transverse intersection expected, and robust once it happens. In addition, the numerical organization of self-localized states reported in \cite{EleonskiiKulagin1989} provided the original indication of such intersection.

\begin{figure}[!htb]
\centering
\includegraphics[width=0.65\linewidth]{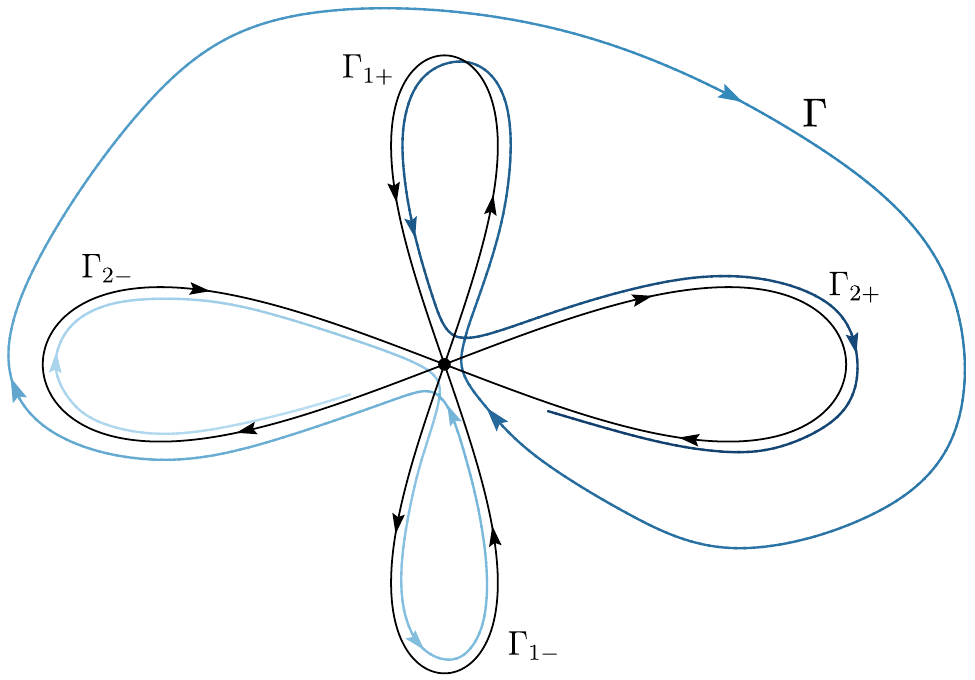}
\caption{Schematic of a super-homoclinic orbit $\Gamma$ (blue) approaching the four-loop bouquet (black) in both time directions. The shading along $\Gamma$ changes from light to dark blue in forward time.}
\label{fig:super_homoclinic}
\end{figure}

An orbit $\Gamma$ at this intersection is called \emph{super-homoclinic}: it is asymptotic to $B$ in both time directions; see Figure \ref{fig:super_homoclinic}. Its past tail follows an unstable component code, read in reverse order when the tail is traversed backward. Its future tail follows a stable component code in forward time. These two codes need not be the same.

Take a super-homoclinic orbit $\Gamma$, and recall $U$ from Theorem \ref{thm:dynamic_near_loops}. Choose a type-1 section $\Pi^{\mathrm{in}}_{1\sigma_1}$ that the past tail visits infinitely often and a type-1 section $\Pi^{\mathrm{in}}_{1\sigma_2}$ that the future tail visits infinitely often. We distinguish a point $p_u$ on the former tail and a later point $p_s$ on the latter so that the negative semiorbit of $p_u$ and the positive semiorbit of $p_s$ stay in $U$, while the intervening segment of $\Gamma$ leaves $U$.

\begin{theoremC}[Homoclinic truncations of a super-homoclinic orbit]\label{thm:superhomoclinic}
Suppose that Assumptions~\ref{assump:exist_transverse_homoclinics}, \ref{assump:gamma_larger_half},\ref{assump:nonzero_coeffs} and~\ref{assum:exist_superhomoclinic} hold. There are tail points $p_u,p_s$ as above, neighborhoods $U_u,U_s$ of these points on their respective sections, and integers $N_u,N_s>0$ such that for every $n\ge N_u$ and $m\ge N_s$ there is a homoclinic orbit $\Gamma_{n,m}$ to $O$ with the following properties:

\begin{enumerate}[(i)]
\item The orbit intersects $U_u$ at a point $q_{n,m}^{u}$ and subsequently intersects $U_s$ at a point $q_{n,m}^{s}$.
\item The segment from $q_{n,m}^{u}$ to $q_{n,m}^{s}$ follows the global passage of $\Gamma$ from $p_u$ to $p_s$, with
\begin{equation*}
q_{n,m}^{u}\longrightarrow p_{u},\qquad
q_{n,m}^{s}\longrightarrow p_{s}\qquad
\text{as }\min\{n,m\}\longrightarrow\infty.
\end{equation*}

\item Traced backward from $q_{n,m}^{u}$, the orbit follows the unstable code of Theorem~\ref{thm:coding_unstable} in reverse order and visits $\Pi^{\mathrm{in}}_{1\sigma_1}$ exactly $n$ further times. After the last visit it makes one more bump along $\Gamma_{1\sigma_1}$ and converges to $O$ in backward time. Traced forward from $q_{n,m}^{s}$, the orbit follows the stable code of Theorem~\ref{thm:coding_stable} and visits $\Pi^{\mathrm{in}}_{1\sigma_2}$ exactly $m$ further times. After the last visit it lies on the local stable manifold of $O$ and converges directly to $O$.

\item Recall $b_2$ and $c_{2}$ from \eqref{eq_Taylor_expansion_global_map}. Define
\begin{equation}\label{eq:bumps_per_visit}
r_u := \begin{cases}
2, & c_2>0,\\
4, & c_2<0,
\end{cases}
\qquad\qquad
r_s := \begin{cases}
2, & b_{2}<0,\\
4, & b_{2}>0.
\end{cases}
\end{equation}
The negative (resp. positive) semiorbit of $q_{n,m}^u$ (resp. $q_{n,m}^s$) contributes $nr_{u} +1$ (resp. $mr_{s}$) bumps to the homoclinic orbit $\Gamma_{n,m}$.
\end{enumerate}

In particular, the system possesses infinitely many multi-bump homoclinic orbits to $O$, with arbitrarily many transitions near the homoclinic bouquet $B$ in both forward and backward time.
\end{theoremC}

\begin{figure}[!htb]
\centering
\includegraphics[width=0.90\linewidth]{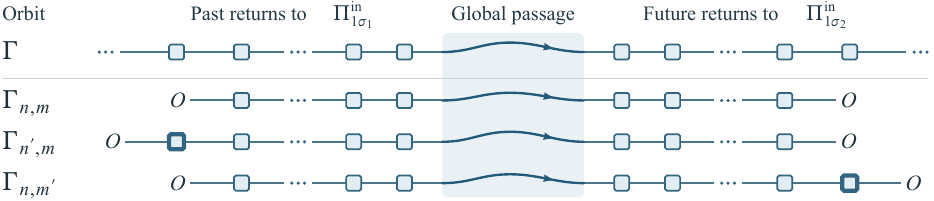}
\caption{Homoclinic orbits shadowing $\Gamma$ (Theorem~\ref{thm:superhomoclinic}). Squares mark section returns in forward order. The bolder squares mark one extra return in the past ($n'=n+1$) or future ($m'=m+1$). $\Gamma$ returns infinitely often on each tail; the shadowing orbits tend to $O$ as $t\to\pm\infty$.}
\label{fig:superhomoclinic_truncations}
\end{figure}

The counts $n$ and $m$ are returns to a section, not necessarily repetitions of a whole component code: an orbit can revisit a section on its other branch. This is why the factors $r_u$ and $r_s$ in \eqref{eq:bumps_per_visit} are determined by the section signs. The central passage is not included in either tail count. Figure~\ref{fig:superhomoclinic_truncations} illustrates how the homoclinic orbits $\Gamma_{n,m}$ can shadow longer portions of either tail of $\Gamma$ while following the same global excursion.

\paragraph{Method and organization.}
We study the dynamics near the homoclinic bouquet by composing Poincar\'e maps between the cross-sections defined in~\eqref{eq:cross-sections}. Each such map is itself a composition: a local map $T_{i\sigma_{0} j\sigma_{1}}^{\mathrm{loc}}$ describes the passage near the equilibrium, and a global map $T_{i}^{\mathrm{glo}}$, introduced above, describes the excursion along the corresponding homoclinic loop outside a small neighborhood of the equilibrium.

The global maps are diffeomorphisms, which we approximate by their Taylor expansions in~\eqref{eq_Taylor_expansion_global_map}. The local maps require a different analysis. The time an orbit spends near the equilibrium can be arbitrarily long. Moreover, an orbit that starts close to a homoclinic orbit on an incoming cross-section need not reach the prescribed outgoing cross-section or return along the chosen loop. It may instead leave the neighborhood of the bouquet. The domains of the local maps therefore have a nontrivial geometry. In particular, the point where a homoclinic orbit meets an incoming cross-section is a singular point of the corresponding local map.

In Section~\ref{sec:local_analysis_near_equilibrium}, we first put the system near the origin into normal form. We then use Shilnikov's boundary value method (see, e.g.,~\cite{PoincareBirkhoff,Dimabook}) to estimate orbit segments that spend a long time near the equilibrium. These estimates allow us to analyze the local Poincar\'e maps. In Section~\ref{sec:Dynamics_near_four_loops}, we use this analysis to study the dynamics near the entire bouquet, including the itineraries on its stable and unstable manifolds and their coding. The proofs of Theorems \ref{thm:dynamic_near_loops}, \ref{thm:coding_unstable}, and ~\ref{thm:coding_stable} are presented in Section~\ref{sec:Dynamics_near_four_loops} and rely on invariant manifolds whose existence is established in Section~\ref{sec:existence_invariant_manifolds}.

In Section~\ref{sec:existence_invariant_manifolds}, we construct the stable and unstable manifolds of the bouquet. The singularities of the Poincar\'e maps prevent us from applying standard invariant-manifold theorems for smooth hyperbolic maps directly. We use the cross-map method of Afraimovich and Shilnikov~\cite{AnnulusPrinciple}. In suitable coordinates, a composition of the Poincar\'e maps admits a cross-map whose differential satisfies the required smallness conditions. The Afraimovich--Shilnikov theory then gives invariant curves on the cross-sections. Following these curves under the flow generates the stable and unstable manifolds of the bouquet. We then prove Theorem~\ref{thm:superhomoclinic} in Section~\ref{sec:proof_superhomoclinic}.

Appendices~\ref{Sec_Proof_Lem_derivatives_local_map_12} and~\ref{Sec_Proof_Lem_derivatives_local_map_21} contain the proofs of two technical lemmas used in Section~\ref{sec:existence_invariant_manifolds}. Appendix~\ref{sec:proofs_technical_lemmas} contains the proofs of further results, most of which involve estimates for Shilnikov boundary value problems. We use Maple for some of the symbolic algebra in the calculations of Appendix \ref{sec:proofs_technical_lemmas}.

\section{Local analysis near the equilibrium}\label{sec:local_analysis_near_equilibrium}

We study the dynamics near the homoclinic orbits by analyzing the associated Poincar\'e maps. The Poincar\'e map is a composition of a local and a global map. In this section, we focus on the local Poincar\'e maps. We first reduce the system near the equilibrium $O$ to a suitable normal form, then estimate the trajectories near $O$ using Shilnikov coordinates, and finally analyze the resulting local maps.

\subsection{Normal form and the estimation of trajectories in Shilnikov coordinates}\label{sec:NF_and_Shilnikov_boundary_value}

We bring the system near $O$ to a normal form and estimate the trajectories near $O$ in this subsection. The technical proofs are postponed. We begin with the statement of the normal form.

\begin{mylem}\label{Lem_normal_form_1}
There exists a near-identity ${C}^{\infty}$ change of coordinates which reduces system \eqref{eq_system_1} near the equilibrium $O$ to the form
\begin{equation}\label{eq:normal_form}
\begin{aligned}
\dot{u}_{1} &= -\lambda_{1} u_{1} + f_{11} u_{1}^{2} v_{1} + f_{12} u_{1} u_{2} v_{2} + f_{13} u_{2}^{3}v_{1}v_{2} + f_{14} u_{2}^{2}v_{1}v_{2}^{2} + f_{15} u_{1}^{3} v_{2}^{2}\\
&\quad + f_{16} u_{2}^{2} v_{1}^{3} + f_{17} u_{1} u_{2}^{2} v_{1}^{2},\\
\dot{u}_{2} &= -\lambda_{2} u_{2} + f_{21} u_{1}^{2}u_{2} + f_{22} u_{1}u_{2}v_{1} + f_{23} u_{2}^{3} + f_{24} u_{2}^{2}v_{2} + f_{25} u_{1}^{2} v_{1}^{2} v_{2}\\
&\quad+ f_{26} u_{1}^{4}v_{2} + f_{27} u_{1}^{3} v_{1} v_{2},\\
\dot{v}_{1} &= +\lambda_{1} v_{1} + g_{11} v_{1}^{2}u_{1} + g_{12} v_{1} v_{2} u_{2} + g_{13} v_{2}^{3}u_{1}u_{2} + g_{14} v_{2}^{2}u_{1}u_{2}^{2} + g_{15} v_{1}^{3}u_{2}^{2}\\
&\quad + g_{16} v_{2}^{2}u_{1}^{3} + g_{17} v_{1} v_{2}^{2} u_{1}^{2},\\
\dot{v}_{2} &= +\lambda_{2} v_{2} + g_{21} v_{1}^{2}v_{2} + g_{22} v_{1}v_{2}u_{1} + g_{23} v_{2}^{3} + g_{24} v_{2}^{2}u_{2} + g_{25} v_{1}^{2} u_{1}^{2} u_{2}\\
&\quad + g_{26} v_{1}^{4} u_{2} + g_{27} v_{1}^{3} u_{1} u_{2},
\end{aligned}
\end{equation}
where $f_{ij}$s and $g_{ij}$s are ${C}^{\infty}$ functions of $(u,v)$. Moreover, first integral \eqref{eq87ywhouyeg3r23} takes the form
\begin{equation}\label{eq34rghuijuy}
H = \gamma u_{1} v_{1} [1 + O(\delta^{2})] -u_{2} v_{2}[1 + O(\delta^{2})] + O(u_{2}^{2} v_{1}^{2}) + O(u_{1}^{2} v_{2}^{2}),
\end{equation}
where $\delta \geq \max\{\lvert u_{1}\rvert, \lvert u_{2}\rvert, \lvert v_{1}\rvert, \lvert v_{2}\rvert\}$ is sufficiently small, and the terms $O(u_{2}^{2} v_{1}^{2})$ and $O(u_{1}^{2} v_{2}^{2})$ are $C^{\infty}$ functions of $(u,v)$. Moreover, system \eqref{eq:normal_form} and first integral \eqref{eq34rghuijuy} remain invariant with respect to symmetries \eqref{eq:symmetry_1} and \eqref{eq:symmetry_2}.
\end{mylem}

This lemma is proved at the end of this subsection. We now start to deal with estimating the trajectories near the origin $O$ of system \eqref{eq:normal_form} using a technique of Shilnikov \cite{PoincareBirkhoff,Dimabook}. Let us first describe his technique. Consider a system
\begin{equation}\label{eq56guyhw}
\begin{aligned}
\dot{u}_{1} &= -\lambda_{1} u_{1} + F_{1}(u_{1}, u_{2}, v_{1}, v_{2}),\\
\dot{v}_{1} &= \lambda_{1} v_{1} + G_{1}(u_{1}, u_{2}, v_{1}, v_{2}),\\
\dot{u}_{2} &= -\lambda_{2} u_{2} + F_{2}(u_{1}, u_{2}, v_{1}, v_{2}),\\
\dot{v}_{2} &= \lambda_{2} v_{2} + G_{2}(u_{1}, u_{2}, v_{1}, v_{2}),
\end{aligned}
\end{equation}
where $\lambda_{2} \geq \lambda_{1} > 0$, and the functions $F_{i}$s and $G_{i}$s, together with their first derivatives vanish at $O$. Take a $\delta > 0$, and define $\Omega := [-\delta, \delta] \times [-\delta, \delta] \times [-\delta, \delta] \times [-\delta, \delta]$. If $\delta$ is small enough, for every $\tau > 0$ and quadruple $(u_{10}, u_{20}, v_{1\tau}, v_{2\tau}) \in \Omega$, system \eqref{eq56guyhw} has a unique solution $(u_{1}^{*}(t), u_{2}^{*}(t), v_{1}^{*}(t), v_{2}^{*}(t))$ that satisfies the boundary conditions
\begin{equation}\label{eq_boundary_conditions_1}
u_{1}^{*}(0) = u_{10}, \quad
u_{2}^{*}(0) = u_{20}, \quad
v_{1}^{*}(\tau) = v_{1\tau}, \quad
v_{2}^{*}(\tau) = v_{2\tau}.
\end{equation}
Fix $\tau > 0$ and $\delta > 0$, and let $\mathcal{I}$ be the set of all continuous functions
\begin{equation*}
t\mapsto (u_{1}(t), u_{2}(t), v_{1}(t), v_{2}(t))
\end{equation*}
from $[0,\tau]$ to $\Omega$. The solution $(u_{1}^{*}(t), u_{2}^{*}(t), v_{1}^{*}(t), v_{2}^{*}(t))$ is a fixed point of the integral operator $\mathfrak{T}: \mathcal{I}\rightarrow \mathcal{I}$ defined by
\begin{equation*}
(u_{1}(t), u_{2}(t), v_{1}(t), v_{2}(t))\mapsto (\overline{u}_{1}(t), \overline{u}_{2}(t), \overline{v}_{1}(t), \overline{v}_{2}(t)),
\end{equation*}
where
\begin{equation}\label{eq5676fuyg8}
\begin{aligned}
\overline{u}_{i}(t) &= e^{-\lambda_{i}t} u_{i0} + \int_{0}^{t} e^{\lambda_{i}(s-t)} F_{i}\left(x(s)\right) ds, \\
\overline{v}_{i}(t) &= e^{-\lambda_{i}(\tau - t)} v_{i\tau} - \int_{t}^{\tau} e^{\lambda_{i}(t-s)} G_{i}\left(x(s)\right) ds
\end{aligned}
\end{equation}
for $i=1,2$ and $x(s) = \left(u_{1}(s), u_{2}(s), v_{1}(s), v_{2}(s)\right)$. However, it is shown that $\mathfrak{T}$ is a contraction on $\mathcal{I}$, and so by Banach fixed point theorem ($\mathcal{I}$ is complete in $C^{0}$-norm), its fixed point is unique and can be estimated by iterating $\mathfrak{T}$. We have used this technique in Section \ref{sec:proofs_technical_lemmas} to prove the following lemma.
\begin{mylem}\label{lem:formula_trajectories}
Assume that $\lambda_{1} < \lambda_{2} <2\lambda_{1}$ (Assumption \ref{assump:gamma_larger_half}). There exists an $M \geq 0$, independent of $\delta$, such that, for all sufficiently small $\delta > 0$ and every quadruple $(u_{10}, u_{20}, v_{1\tau}, v_{2\tau}) \in \Omega$, the solution $(u_{1}(t), u_{2}(t), v_{1}(t), v_{2}(t))$ of system \eqref{eq:normal_form} satisfying boundary conditions \eqref{eq_boundary_conditions_1} is given by
\begin{equation*}
u_{i}(t) = e^{-\lambda_{i} t} u_{i0} + \xi_{i}(x), \qquad
v_{i}(t) = e^{-\lambda_{i} (\tau - t)} v_{i\tau} + \zeta_{i}(x),\qquad (i=1,2),
\end{equation*}
for $x = (t, u_{10}, u_{20}, v_{1\tau}, v_{2\tau})$, where $\xi_{1}$ and $\xi_2$ satisfy
\begin{equation}\label{eq:estimate_trajectory_1}
\begin{aligned}
\lvert \xi_{1}(x)\rvert
&\leq M \Big(e^{-\lambda_{1}t}\lvert u_{10}\rvert \delta^{2} + e^{-\lambda_{1} t +\left(-\lambda_{1}-\lambda_{2}\right) \tau} \lvert u_{20}^{2} v_{1\tau} v_{2\tau}\rvert + e^{-3 \lambda_{1} \tau +\left(3 \lambda_{1}-2 \lambda_{2}\right) t} \lvert u_{20}^{2} v_{1\tau}^{3}\rvert\\
&\quad + e^{-\lambda_{1} t -3 \lambda_{1} \tau} \lvert u_{20}^{3} v_{1\tau}^{2}\rvert + e^{\left(-\lambda_{1}-2 \lambda_{2}\right) \tau +\lambda_{1} t} \lvert u_{20}^{2} v_{1\tau} v_{2\tau}^{2}\rvert \Big),\\
\lvert\xi_{2}(x)\rvert
&\leq M\Big(e^{-\lambda_{2}t}\lvert u_{20}\rvert \delta^{2} + e^{\left(-2 \lambda_{1}-\lambda_{2}\right) \tau +\lambda_{2} t} \lvert u_{10}^{2} v_{1\tau}^{2} v_{2\tau} \rvert + e^{-\lambda_{2} t -\lambda_{2} \tau} \lvert u_{10}^{3} v_{2\tau} \rvert \\
&\quad + e^{-3 \lambda_{2} \tau +\left(-4 \lambda_{1}+3 \lambda_{2}\right) t} \lvert u_{10}^{2} v_{2\tau}^{3} \rvert + e^{\left(-\lambda_{1}-\lambda_{2}\right) \tau +\left(-2 \lambda_{1}+\lambda_{2}\right) t} \lvert u_{10}^{2} v_{1\tau} v_{2\tau} \rvert \Big),
\end{aligned}
\end{equation}
and $\zeta_{1}$ and $\zeta_2$ satisfy
\begin{equation}\label{eq:estimate_trajectory_2}
\begin{aligned}
\lvert \zeta_{1}(x)\rvert
&\leq M \Big(e^{-\lambda_{1}(\tau - t)} \lvert v_{1\tau}\rvert \delta^{2} + e^{-2 \lambda_{2} \tau +\left(-3 \lambda_{1}+2 \lambda_{2}\right) t} \lvert u_{10}^{3} v_{2\tau}^{2} \rvert + e^{-\lambda_{1} t -2 \lambda_{2} \tau} \lvert u_{10} u_{20}^{2} v_{2\tau}^{2} \rvert\\
&\quad + e^{\lambda_{1} t -4 \lambda_{1} \tau} \lvert u_{10}^{2} v_{2\tau}^{3} \rvert + e^{\left(-2 \lambda_{1}-\lambda_{2}\right) \tau +\lambda_{1} t} \lvert u_{10} u_{20} v_{2\tau}^{2} \rvert\Big),\\
\lvert \zeta_{2}(x)\rvert
&\leq M \Big(e^{-\lambda_{2}(\tau - t)} \lvert v_{2\tau}\rvert \delta^{2} + e^{-3 \lambda_{1} \tau +\left(2 \lambda_{1}-\lambda_{2}\right) t} \lvert u_{10} u_{20} v_{1\tau}^{2} \rvert + e^{\lambda_{2} t -2 \lambda_{2} \tau} \lvert u_{20} v_{1\tau}^{3} \rvert\\
&\quad + e^{-\lambda_{2} t -2 \lambda_{1} \tau} \lvert u_{10}^{2} u_{20} v_{1\tau}^{2} \rvert + e^{-4 \lambda_{1} \tau +\left(4 \lambda_{1}-3 \lambda_{2}\right) t} \lvert u_{20}^{3} v_{1\tau}^{2} \rvert\Big).
\end{aligned}
\end{equation}
\end{mylem}

Note that by evaluating \eqref{eq:estimate_trajectory_1} and \eqref{eq:estimate_trajectory_2} at $t=\tau$ and $t=0$, respectively, and simplifying the result further, we obtain
\begin{equation}\label{eq_Poincare_map}
\begin{aligned}
u_{1\tau} &= e^{-\lambda_{1}\tau} u_{10} [1 + O(\delta^{2})] + e^{-\lambda_{2} \tau} O(\lvert u_{20} v_{1\tau} v_{2\tau}\rvert) + e^{-2 \lambda_{2} \tau} O(u_{20}^{2} v_{1\tau}^{2}),\\[2pt]
u_{2\tau} &= e^{-\lambda_{2}\tau} u_{20} [1 + O(\delta^{2})] + e^{-3 \lambda_{1} \tau} O(\lvert u_{10}^{2} v_{1\tau} v_{2\tau}\rvert) + e^{-2\lambda_{1} \tau} O(\lvert u_{10}^{2} v_{1\tau}^{2} v_{2\tau}\rvert)\\
&\quad + e^{-4 \lambda_{1} \tau} O(\lvert u_{10}^{2} v_{2\tau}^{3}\rvert) + e^{-2\lambda_{2} \tau} O(\lvert u_{10}^{3} v_{2\tau}\rvert),\\[2pt]
v_{10} &= e^{-\lambda_{1}\tau} v_{1\tau} [1+O(\delta^{2})] + e^{-\lambda_{2} \tau} O(\lvert u_{10} u_{20} v_{2\tau}\rvert) + e^{-2 \lambda_{2} \tau} O(u_{10}^{2} v_{2\tau}^{2}),\\[2pt]
v_{20} &= e^{-\lambda_{2}\tau} v_{2\tau} [1+O(\delta^{2})] + e^{-2\lambda_{1}\tau} O(\lvert u_{10}^{2} u_{20} v_{1\tau}^{2}\rvert) + e^{-4 \lambda_{1} \tau} O(\lvert u_{20}^{3} v_{1\tau}^{2}\rvert)\\
&\quad  + e^{-3 \lambda_{1} \tau} O(\lvert u_{10} u_{20} v_{1\tau}^{2}\rvert) + e^{-2 \lambda_{2} \tau} O(\lvert u_{20} v_{1\tau}^{3}\rvert),
\end{aligned}
\end{equation}
where $u_{1\tau} := u_{1}(\tau)$, $u_{2\tau} := u_{2}(\tau)$, $v_{10} := v_{1}(0)$, and $v_{20} := v_{2}(0)$.

\begin{proof}[Proof of Lemma \ref{Lem_normal_form_1} (the normal form lemma)]
Following Lemma 3.6 in \cite{BakraniPhDthesis} and taking symmetries \eqref{eq:symmetry_1} and \eqref{eq:symmetry_2} into account, we can reduce system \eqref{eq_system_1} near the origin to the form
\begin{equation}\label{PhD_normal_form}
\begin{aligned}
\dot{u}_{1} &= -\lambda_{1} u_{1} + \alpha_{11} u_{1}^{2} v_{1} + \alpha_{12} u_{1} u_{2} v_{2} + \alpha_{13} u_{2}v_{1}v_{2} + \alpha_{14} u_{1}^{3} v_{2}^{2} + \alpha_{15} u_{2}^{2} v_{1},\\
\dot{u}_{2} &= -\lambda_{2} u_{2} + \alpha_{21} u_{1}^{2}u_{2} + \alpha_{22} u_{1}u_{2}v_{1} + \alpha_{23} u_{2}^{3} + \alpha_{24} u_{2}^{2}v_{2} + \alpha_{25} u_{1}^{2}v_{2},\\
\dot{v}_{1} &= +\lambda_{1} v_{1} + \beta_{11} v_{1}^{2}u_{1} + \beta_{12} v_{1} v_{2} u_{2} + \beta_{13} v_{2}u_{1}u_{2} + \beta_{14} v_{1}^{3}u_{2}^{2} + \beta_{15} v_{2}^{2}u_{1},\\
\dot{v}_{2} &= +\lambda_{2} v_{2} + \beta_{21} v_{1}^{2}v_{2} + \beta_{22} v_{1}v_{2}u_{1} + \beta_{23} v_{2}^{3} + \beta_{24} v_{2}^{2}u_{2} + \beta_{25} v_{1}^{2}u_{2},
\end{aligned}
\end{equation}
where $\alpha_{ij}$s and $\beta_{ij}$s are some $C^{\infty}$ functions of $(u,v)$. Moreover, first integral \eqref{eq87ywhouyeg3r23} is reduced to the form
\begin{equation}\label{eq64guy}
H = \gamma u_{1} v_{1}\left[1 + H_{1}(u,v)\right] - u_{2} v_{2}\left[1 + H_{2}(u,v)\right],
\end{equation}
where $H_{1}$ and $H_{2}$ are functions of $(u,v)$ that vanish at the origin.

Taking symmetries \eqref{eq:symmetry_1} and \eqref{eq:symmetry_2} into account, we can write \eqref{PhD_normal_form} in the form
\begin{equation}\label{eq3789iyg2ref}
\begin{aligned}
\dot{u}_{1} &= P_{1} + \underline{\phi_{1}(u_{2}) u_{2}^{2}v_{1}} + \underline{\phi_{3}(v_{1},v_{2})u_{2}v_{1}v_{2}},\\
\dot{u}_{2} &= P_{2} + \underline{\phi_{2}(v_{2}) u_{1}^{2}v_{2}},\\
\dot{v}_{1} &= Q_{1} + \underline{\psi_{1}(v_{2}) v_{2}^{2}u_{1}} + \underline{\psi_{3}(u_{1},u_{2})v_{2}u_{1}u_{2}},\\
\dot{v}_{2} &= Q_{2} + \underline{\psi_{2}(u_{2}) v_{1}^{2} u_{2}}.
\end{aligned}
\end{equation}
Here, $P_1$, $P_{2}$, $Q_{1}$ and $Q_{2}$ are functions of $(u_{1}, u_{2}, v_{1}, v_{2})$ given by
\begin{equation*}
\begin{aligned}
P_{1} &:= -\lambda_{1} u_{1} + f_{11} u_{1}^{2} v_{1} + f_{12} u_{1} u_{2} v_{2} + f_{13} u_{2}^{3}v_{1}v_{2} + f_{14} u_{2}^{2}v_{1}v_{2}^{2} + f_{15} u_{1}^{3} v_{2}^{2}\\
&\quad + f_{16} u_{2}^{2} v_{1}^{3} + f_{17} u_{1} u_{2}^{2} v_{1}^{2},\\
P_{2} &:= -\lambda_{2} u_{2} + f_{21} u_{1}^{2}u_{2} + f_{22} u_{1}u_{2}v_{1} + f_{23} u_{2}^{3} + f_{24} u_{2}^{2}v_{2} + f_{25} u_{1}^{2} v_{1}^{2} v_{2}\\
&\quad + f_{26} u_{1}^{4}v_{2} + f_{27} u_{1}^{3} v_{1} v_{2},\\
Q_{1} &:= +\lambda_{1} v_{1} + g_{11} v_{1}^{2}u_{1} + g_{12} v_{1} v_{2} u_{2} + g_{13} v_{2}^{3}u_{1}u_{2} + g_{14} v_{2}^{2}u_{1}u_{2}^{2} + g_{15} v_{1}^{3}u_{2}^{2}\\
&\quad + g_{16} v_{2}^{2}u_{1}^{3} + g_{17} v_{1} v_{2}^{2} u_{1}^{2},\\
Q_{2} &:= +\lambda_{2} v_{2} + g_{21} v_{1}^{2}v_{2} + g_{22} v_{1}v_{2}u_{1} + g_{23} v_{2}^{3} + g_{24} v_{2}^{2}u_{2} + g_{25} v_{1}^{2} u_{1}^{2} u_{2}\\
&\quad + g_{26} v_{1}^{4} u_{2} + g_{27} v_{1}^{3} u_{1} u_{2},
\end{aligned}
\end{equation*}
where $f_{ij}$s and $g_{ij}$s are some ${C}^{\infty}$ functions of $(u,v)$. Moreover,
\begin{equation*}
\begin{alignedat}{2}
\phi_{1}(u_{2}) &= \alpha_{15}(0,u_{2},0,0), \qquad & \psi_{1}(v_{2}) &= \beta_{15}(0,0,0,v_{2}),\\[2pt]
\phi_{2}(v_{2}) &= \alpha_{25}(0,0,0,v_{2}), \qquad & \psi_{2}(u_{2}) &= \beta_{25}(0,u_{2},0,0),\\[2pt]
\phi_{3}(v_{1},v_{2}) &= \alpha_{13}(0,0,v_{1},v_{2}), \qquad & \psi_{3}(u_{1},u_{2}) &= \beta_{13}(u_{1},u_{2},0,0).
\end{alignedat}
\end{equation*}
To reduce this system to normal form \eqref{eq:normal_form}, it is enough to eliminate the underlined terms from \eqref{eq3789iyg2ref}. It is done by applying six consecutive near-identity changes of coordinates as follows.

We first eliminate $\phi_{1}(u_{2}) u_{2}^{2} v_{1}$ by applying a change of coordinates of the form
\begin{equation*}
\varphi_{1}(u_{1}, u_{2}, v_{1}, v_{2}) = \left(u_{1} + u_{2}^{2} v_{1} p_{1}(u_{2}), u_{2}, v_{1}, v_{2}\right),
\end{equation*}
where $p_{1}$ is specified below. It is straightforward to see that applying $\varphi_{1}$ keeps the forms of $\dot{u}_{2}$, $\dot{v}_{1}$ and $\dot{v}_{2}$ in system \eqref{eq3789iyg2ref}. The equation of $\dot{u}_{1}$ in the updated system is given by $\dot{u}_{1} = Z_{1} + u_{2}^{2} v_{1} R_{1}(u_{2})$, where $Z_{1}$ can be written in the form of $P_{1} + \phi_{3}(v_{1},v_{2})u_{2}v_{1}v_{2}$, and
\begin{equation*}
\begin{aligned}
R_{1}(u_{2}) &= \left[2\left(\lambda_{1}-\lambda_{2}\right) + 2 u_{2}^{2} f_{23}\left(0, u_{2}, 0, 0\right)\right] p_{1}(u_{2})\\
&\quad + \left[-\lambda_{2} u_{2} + u_{2}^{3} f_{23}\left(0, u_{2}, 0, 0\right)\right] p_{1}^{\prime}\left(u_{2}\right) + \phi_{1}\left(u_{2}\right).
\end{aligned}
\end{equation*}
The goal is to find $p_{1}$ such that $R_{1}$ vanishes. Consider the planar system
\begin{equation*}
\begin{aligned}
\dot{x} &= 2\left(\lambda_{1}-\lambda_{2}\right)x + 2 x u_{2}^{2} f_{23}\left(0, u_{2}, 0, 0\right) + \phi_{1}\left(u_{2}\right),\\
\dot{u}_{2} &= \lambda_{2} u_{2} - u_{2}^{3} f_{23}\left(0, u_{2}, 0, 0\right).
\end{aligned}
\end{equation*}
This system has an equilibrium at $(\frac{\phi_{1}(0)}{2(\lambda_{2} -\lambda_{1})}, 0)$ with the linearization $\bigl(\begin{smallmatrix}
2\left(\lambda_{1} -\lambda_{2}\right) & 0\\
0 & \lambda_{2}
\end{smallmatrix}\bigr)$ (note that $\phi_{1}^{\prime}(0) = 0$ due to the symmetries). This equilibrium possesses a local unstable manifold given by $x = p_{1}(u_{2})$ for some function $p_{1}$. This function $p_{1}$ makes $R_1$ vanish because the constraint $R_{1} = 0$ is the invariance condition of this manifold (see \cite{BakraniPhDthesis}, Section 2.3.1).

We now eliminate $\phi_{2}(v_{2}) u_{1}^{2} v_{2}$ by applying a change of coordinates of the form
\begin{equation*}
\varphi_{2}(u_{1}, u_{2}, v_{1}, v_{2}) = \left(u_{1}, u_{2} + u_{1}^{2} v_{2} p_{2}(v_{2}), v_{1}, v_{2}\right),
\end{equation*}
where $p_{2}$ is given below. Call the system obtained from applying $\varphi_{1}$ System I (this system is of the form of \eqref{eq3789iyg2ref} with the term $\phi_{1}(u_{2}) u_{2}^{2} v_{1}$ removed). It is straightforward to see that applying $\varphi_{2}$ keeps the form of $\dot{u}_{1}$, $\dot{v}_{1}$ and $\dot{v}_{2}$ of System I. The equation of $\dot{u}_{2}$ in the updated system is given by $\dot{u}_{2} = Z_{2} + u_{1}^{2} v_{2} R_{2}(v_{2})$, where $Z_{2}$ is of the form of $P_{2}$, and
\begin{equation*}
\begin{aligned}
R_{2}(v_{2}) &= \left[2\left(\lambda_{2} - \lambda_{1}\right) + v_{2}^{2} g_{23}\left(0,0,0,v_{2}\right)\right] p_{2}(v_{2})\\
&\quad + \left[\lambda_{2} v_{2} + v_{2}^{3} g_{23}\left(0,0,0,v_{2}\right)\right] p_{2}^{\prime}(v_{2}) + \phi_{2}(v_{2}).
\end{aligned}
\end{equation*}
The planar system
\begin{equation*}
\begin{aligned}
\dot{x} &= 2\left(\lambda_{2} - \lambda_{1}\right) x + x v_{2}^{2} g_{23}\left(0,0,0,v_{2}\right) + \phi_{2}(v_{2}),\\
\dot{v}_{2} &= -\lambda_{2} v_{2} - v_{2}^{3} g_{23}\left(0,0,0,v_{2}\right)
\end{aligned}
\end{equation*}
has an equilibrium at $(\frac{\phi_{2}(0)}{2(\lambda_{1} -\lambda_{2})}, 0)$ with the linearization $\bigl(\begin{smallmatrix}
2\left(\lambda_{2} -\lambda_{1}\right) & 0\\
0 & -\lambda_{2}
\end{smallmatrix}\bigr)$ (note that $\phi_{2}^{\prime}(0) = 0$ due to the symmetries). This equilibrium has a stable manifold $\{x = p_{2}(v_{2})\}$ for some function $p_{2}$. Then, this function $p_2$ makes $R_{2}$ vanish, as desired.

To eliminate $\phi_{3}(v_{1},v_{2}) u_{2} v_{1} v_{2}$, we apply a change of coordinates of the form
\begin{equation*}
\varphi_{3}(u_{1}, u_{2}, v_{1}, v_{2}) = \left(u_{1} + u_{2} v_{1} v_{2} p_{3}(v_{1},v_{2}), u_{2}, v_{1}, v_{2}\right),
\end{equation*}
where $p_{3}$ is given below. Call the system obtained from applying $\varphi_{2}\circ\varphi_{1}$ System II (it is of the form of \eqref{eq3789iyg2ref} with the terms $\phi_{1}(u_{2}) u_{2}^{2} v_{1}$ and $\phi_{2}(v_{2}) u_{1}^{2} v_{2}$ removed). It is straightforward to see that applying $\varphi_{3}$ keeps the form of $\dot{u}_{2}$, $\dot{v}_{1}$ and $\dot{v}_{2}$ of System II. The equation of $\dot{u}_{1}$ in the updated system is given by $\dot{u}_{1} = Z_{3} + u_{2} v_{1} v_{2} R_{3}(v_{1},v_{2})$, where $Z_{3}$ is of the form of $P_{1}$, and
\begin{equation*}
\begin{aligned}
R_{3}(v_{1}, v_{2}) &= \left[2\lambda_{1} + v_{1}^{2} g_{21}\left(0,0,v_{1},v_{2}\right) + v_{2}^{2} g_{23}\left(0,0,v_{1},v_{2}\right)\right] p_{3}(v_{1}, v_{2})\\
&\quad + \phi_{3}(v_{1},v_{2}) + \lambda_{1} v_{1}\cdot\frac{\partial p_{3}}{\partial v_{1}}(v_{1}, v_{2})\\
&\quad+ \left[\lambda_{2} v_{2} + v_{1}^{2} v_{2} g_{21}\left(0,0,v_{1},v_{2}\right) + v_{2}^{3} g_{23}\left(0,0,v_{1},v_{2}\right)\right] \frac{\partial p_{3}}{\partial v_{2}}(v_{1}, v_{2}).
\end{aligned}
\end{equation*}
The 3-dimensional system
\begin{equation*}
\begin{aligned}
\dot{x} &= 2\lambda_{1} x + v_{1}^{2} x g_{21}(0,0,v_{1},v_{2}) + v_{2}^{2} x g_{23}(0,0,v_{1},v_{2}) + \phi_{3}(v_{1},v_{2}),\\
\dot{v}_{1} &= -\lambda_{1} v_{1},\\
\dot{v}_{2} &= -\lambda_{2} v_{2} - v_{1}^{2} v_{2} g_{21}\left(0,0,v_{1},v_{2}\right) - v_{2}^{3} g_{23}\left(0,0,v_{1},v_{2}\right)
\end{aligned}
\end{equation*}
has an equilibrium at $(\frac{-\phi_{3}(0,0)}{2\lambda_{1}},0,0)$. The linearization of this system at this equilibrium is
\begin{equation*}
\begin{pmatrix}
2\lambda_{1} & \frac{\partial \phi_{3}}{\partial v_{1}}(0,0) & \frac{\partial \phi_{3}}{\partial v_{2}}(0,0)\\
0 & -\lambda_{1} & 0\\
0 & 0 & -\lambda_{2}
\end{pmatrix}.
\end{equation*}
Note that $\frac{\partial \phi_{3}}{\partial v_{1}}(0,0) = \frac{\partial \phi_{3}}{\partial v_{2}}(0,0) = 0$ due to the symmetries. Similar to the previous two cases, $R_3$ vanishes if we choose $p_{3}$ such that $x = p_{3}(v)$ gives the stable manifold of the equilibrium.

Analogously, we can eliminate the terms $\psi_{1}(v_{2}) u_{1} v_{2}^{2}$, $\psi_{2}(u_{2}) v_{1}^{2} u_{2}$ and $\psi_{3}(u_{1},u_{2})v_{2}u_{1}u_{2}$ by applying three changes of coordinates $\varphi_{4}$, $\varphi_{5}$ and $\varphi_{6}$, respectively, where
\begin{equation*}
\begin{aligned}
\varphi_{4}(u_{1}, u_{2}, v_{1}, v_{2}) &= \left(u_{1}, u_{2}, v_{1} + v_{2}^{2} u_{1} p_{4}(v_{2}), v_{2}\right),\\[1pt]
\varphi_{5}(u_{1}, u_{2}, v_{1}, v_{2}) &= \left(u_{1}, u_{2}, v_{1}, v_{2} + v_{1}^{2} u_{2} p_{5}(u_{2})\right),\\[1pt]
\varphi_{6}(u_{1}, u_{2}, v_{1}, v_{2}) &= \left(u_{1}, u_{2}, v_{1} + v_{2} u_{1} u_{2} p_{6}(u_{1}, u_{2}), v_{2}\right),
\end{aligned}
\end{equation*}
for some appropriate $p_{4}$, $p_{5}$ and $p_{6}$. This proves that system \eqref{eq_system_1} near $O$ can be reduced to the form \eqref{eq:normal_form}.

Let us now see how these changes of coordinates give first integral \eqref{eq34rghuijuy}. Note that
\begin{equation*}
\begin{aligned}
\varphi_{3}\circ \varphi_{2}\circ \varphi_{1}&(u_{1}, u_{2}, v_{1}, v_{2})\\
&= \varphi_{3} \left(u_{1} + O(\lvert u_{2}^{2} v_{1}\rvert),\, u_{2}[1 + O(\delta^{2})] + O(\lvert u_{1}^{2} v_{2}\rvert),\, v_{1},\, v_{2}\right)\\
&= \left(u_{1}[1+O(\delta^{2})] + O(\lvert u_{2}^{2} v_{1}\rvert) + O(\lvert u_{2} v_{2}\delta\rvert),\, u_{2}[1 + O(\delta^{2})] + O(\lvert u_{1}^{2} v_{2}\rvert),\, v_{1},\, v_{2}\right)
\end{aligned}
\end{equation*}
and
\begin{equation*}
\begin{aligned}
\varphi_{6}\circ\varphi_{5}\circ \varphi_{4}&(u_{1}, u_{2}, v_{1}, v_{2})\\
&= \varphi_{6}\left(u_{1},\, u_{2},\, v_{1} + O(\lvert v_{2}^{2} u_{1}\rvert),\, v_{2}[1 + O(\delta^{2})] + O(\lvert v_{1}^{2} u_{2}\rvert)\right)\\
&= \left(u_{1},\, u_{2},\, v_{1}[1+O(\delta^{2})] + O(\lvert v_{2}^{2} u_{1}\rvert) + O(\lvert v_{2} u_{2} \delta\rvert),\, v_{2}[1 + O(\delta^{2})] + O(\lvert v_{1}^{2} u_{2}\rvert)\right).
\end{aligned}
\end{equation*}
Therefore,
\begin{equation*}
\begin{aligned}
\varphi_{6}\circ \cdots \circ \varphi_{1}(u_{1}, u_{2}, v_{1}, v_{2}) &= \Big(u_{1}[1+O(\delta^{2})] + O(\lvert u_{2}^{2} v_{1}\rvert) + O(\lvert u_{2} v_{2}\delta\rvert),\\
&\qquad u_{2}[1 + O(\delta^{2})] + O(\lvert u_{1}^{2} v_{2}\rvert),\\
&\qquad v_{1} [1+O(\delta^{2})] + O(\lvert v_{2}^{2} u_{1}\rvert) + O(\lvert v_{2} u_{2} \delta\rvert),\\
&\qquad v_{2}[1 + O(\delta^{2})] + O(\lvert v_{1}^{2} u_{2}\rvert)\Big).
\end{aligned}
\end{equation*}
Using this formula, it is straightforward to see that updating \eqref{eq64guy} with $\varphi_{6}\circ \cdots \circ \varphi_{1}$ gives \eqref{eq34rghuijuy}.

Following Proposition 2.24 of \cite{BakraniPhDthesis}, all $\varphi_{i}$s as well as all the changes of coordinates made in the proof of Lemma~3.6 of \cite{BakraniPhDthesis} commute with symmetries \eqref{eq:symmetry_1} and \eqref{eq:symmetry_2}. This ends the proof of Lemma~\ref{Lem_normal_form_1}.
\end{proof}

\subsection{Analysis of the local Poincar\'e maps}

We suppose that the orbit $\Gamma_{i+}$ (resp. $\Gamma_{i-}$), where $i=1,2$, leaves $O$ along the positive (resp. negative) side of $v_{i}$-axis and reaches $O$ along the positive (resp. negative) side of $u_{i}$-axis. This convention may require changing the signs of some of the variables $u_{i}$ and $v_{i}$ which keeps the form of system \eqref{eq:normal_form}, but may change first integral \eqref{eq34rghuijuy} to
\begin{equation}\label{eq:first_integral}
H = \gamma u_{1} v_{1} [1 + O(\delta^{2})] - \kappa u_{2} v_{2}[1 + O(\delta^{2})] + O(u_{2}^{2} v_{1}^{2}) + O(u_{1}^{2} v_{2}^{2}),
\end{equation}
where $\kappa = \pm 1$. Recall the cross-sections $\Pi_{i\sigma}^{\mathrm{in}}$ and $\Pi_{i\sigma}^{\mathrm{out}}$ from \eqref{eq:cross-sections}. The following lemma is an immediate consequence of the implicit function theorem and taking the condition $H=0$ into account.
\begin{mylem}\label{Lem_coordinates_on_cross_sections}
Take arbitrary $i\in\{1,2\}$ and $\sigma = \pm$. Let $j=2$ if $i=1$, and $j=1$ if $i=2$. Then, we can choose $(u_{j}, v_{j})$-coordinates on $\Pi^{\mathrm{in}}_{i\sigma}$ and $\Pi^{\mathrm{out}}_{i\sigma}$, i.e. $(u_{i},v_{i})$ is uniquely determined by $(u_{j}, v_{j})$ for any point $(u_{1}, u_{2}, v_{1}, v_{2})$ on $\Pi^{\mathrm{in}}_{i\sigma}$ or $\Pi^{\mathrm{out}}_{i\sigma}$. In particular, there exists a ${C}^{\infty}$ smooth function $h_{i\sigma}^{\mathrm{in}}$ such that $v_{i} = h_{i\sigma}^{\mathrm{in}}(u_j,v_j)$ on $\Pi_{i\sigma}^{\mathrm{in}}$. More precisely,
\begin{equation}\label{eq56yigkre4r}
\begin{aligned}
h_{1\sigma}^{\mathrm{in}}(u_2,v_2) &= \frac{\kappa\sigma}{\gamma\delta}\,u_{2} v_{2} [1+O(\delta^2)] + O(\delta v_2^2),\\
h_{2\sigma}^{\mathrm{in}}(u_1,v_1) &= \frac{\kappa\gamma\sigma}{\delta}\,u_{1} v_{1} [1+O(\delta^2)] + O(\delta v_1^2).
\end{aligned}
\end{equation}
\end{mylem}

Note also that each homoclinic orbit $\Gamma_{i\sigma}$ ($i=1,2$ and $\sigma = \pm$) intersects $\Pi^{\mathrm{in}}_{i\sigma}$ and $\Pi^{\mathrm{out}}_{i\sigma}$ at $(0,0)$. This is because each loop lies at the intersection of one of the invariant planes with the locally straightened stable and unstable manifolds of the equilibrium $O$.

\begin{figure}[!htb]
\centering
\includegraphics[width=0.9\linewidth]{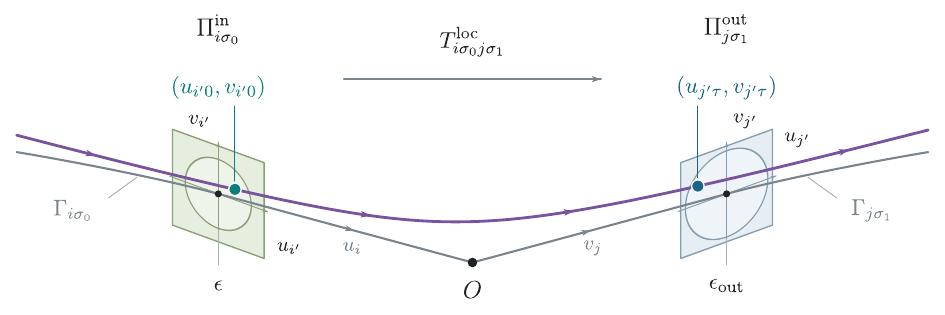}
\caption{Schematic local passage $T^{\mathrm{loc}}_{i\sigma_0j\sigma_1}$.
The purple trajectory enters through the $\epsilon$-disk on
$\Pi^{\mathrm{in}}_{i\sigma_0}$, follows $\Gamma_{i\sigma_0}$ toward $O$,
passes near $O$, and exits through the $\epsilon_{\mathrm{out}}$-disk on
$\Pi^{\mathrm{out}}_{j\sigma_1}$ near $\Gamma_{j\sigma_1}$. The admissible incoming points form the domain $D_{i\sigma_{0} j\sigma_1}$; $i$ and $j$ may coincide. Black dots mark $O$ and the section origins.}
\label{fig:domain_of_local_map}
\end{figure}

Take real $\epsilon$ and $\epsilon_{\mathrm{out}}$ such that $0 < \epsilon \leq \epsilon_{\mathrm{out}} < \delta$, and two integers $i, j \in \{1,2\}$. Let $i^{\prime}, j^{\prime} \in \{1,2\}$ such that $i^{\prime}\neq i$ and $j^{\prime} \neq j$. For any two $\sigma_{0}, \sigma_{1} = \pm$, let $D_{i\sigma_{0} j\sigma_{1}}$ be the set of all points $(u_{i^{\prime}}, v_{i^{\prime}})$ that at time $t = 0$ are in the $\epsilon$-ball of $(0,0)$ in $\Pi^{\mathrm{in}}_{i\sigma_{0}}$, their forward orbits go along $\Gamma_{i\sigma_{0}}$, get close to the equilibrium $O$, then go along $\Gamma_{j\sigma_{1}}$, and at $t=\tau$ (for some $\tau > 0$), they reach $\Pi^{\mathrm{out}}_{j\sigma_{1}}$ at a point $(u_{j^{\prime}}, v_{j^{\prime}})$ which lies in the $\epsilon_{\mathrm{out}}$-ball centered at $(0,0) \in \Pi^{\mathrm{out}}_{j\sigma_{1}}$ (See Figure \ref{fig:domain_of_local_map}). Throughout, we often incorporate the start and end times in the coordinates, too, and say the point $(u_{i^{\prime} 0}, v_{i^{\prime} 0})$ reaches the point $(u_{j^{\prime} \tau}, v_{j^{\prime} \tau})$. If $D_{i\sigma_{0} j\sigma_{1}} \neq \emptyset$, we define the local Poincar\'e map by
\begin{equation*}
\begin{aligned}
T^{\mathrm{loc}}_{i\sigma_{0} j\sigma_{1}}: D_{i\sigma_{0} j\sigma_{1}} &\rightarrow \Pi^{\mathrm{out}}_{j\sigma_{1}}\\
(u_{i^{\prime} 0}, v_{i^{\prime} 0}) &\mapsto (u_{j^{\prime} \tau}, v_{j^{\prime} \tau}).
\end{aligned}
\end{equation*}
The forward orbit of $(u_{j^{\prime}\tau}, v_{j^{\prime}\tau}) \in \Pi^{\mathrm{out}}_{j\sigma_{1}}$ continues going along the homoclinic orbit $\Gamma_{j\sigma_{1}}$ until it intersects $\Pi^{\mathrm{in}}_{j\sigma_{1}}$ at a point $(u_{j^{\prime}}, v_{j^{\prime}})$. We define the global Poincar\'e map $T^{\mathrm{glo}}_{j\sigma_{1}}$ from the $\epsilon_{\mathrm{out}}$-ball around $(0,0)\in \Pi_{j\sigma_{1}}^{\mathrm{out}}$ to $\Pi_{j\sigma_{1}}^{\mathrm{in}}$ by $(u_{j^{\prime}\tau}, v_{j^{\prime}\tau}) \mapsto (u_{j^{\prime}}, v_{j^{\prime}})$. However, symmetries \eqref{eq:symmetry_1} and \eqref{eq:symmetry_2} imply $T^{\mathrm{glo}}_{1-} = T^{\mathrm{glo}}_{1+}$ and $T^{\mathrm{glo}}_{2-} = T^{\mathrm{glo}}_{2+}$. Setting $T^{\mathrm{glo}}_{j} := T^{\mathrm{glo}}_{j\sigma_{1}}$, we then define the Poincar\'e map 
\begin{equation*}
T_{i\sigma_{0} j\sigma_{1}}: = T^{\mathrm{glo}}_{j} \circ T^{\mathrm{loc}}_{i\sigma_{0} j\sigma_{1}}.
\end{equation*}
Note also that the global maps are diffeomorphisms, and Taylor expanding them at $(0,0)$ gives \eqref{eq_Taylor_expansion_global_map}.

For convenience, we combine the domains $D_{i\sigma_{0}j\sigma_{1}}$ into a single domain, and the Poincar\'e maps $T_{i\sigma_{0}j\sigma_{1}}$ into a single map as in the following definition.
\begin{mydefn}\label{dfn:Poincare_map_T}
Let $D$ be the union of all the sets $D_{i\sigma_{0} j\sigma_{1}}$ for $i,j=1,2$ and $\sigma_{0}, \sigma_{1} = \pm$. For any point $z\in D$, there exists a unique map $T_{i\sigma_{0} j\sigma_{1}}$ that has $z$ in its domain. We define the (unified) Poincar\'e map $T : D \rightarrow \mathop{\bigcup}\limits_{i\sigma} \Pi^{\mathrm{in}}_{i\sigma}$ by $T(z) := T_{i\sigma_{0} j\sigma_{1}}(z)$.
\end{mydefn}

\begin{figure}[!tbh]
\centering
\includegraphics[width=\linewidth]{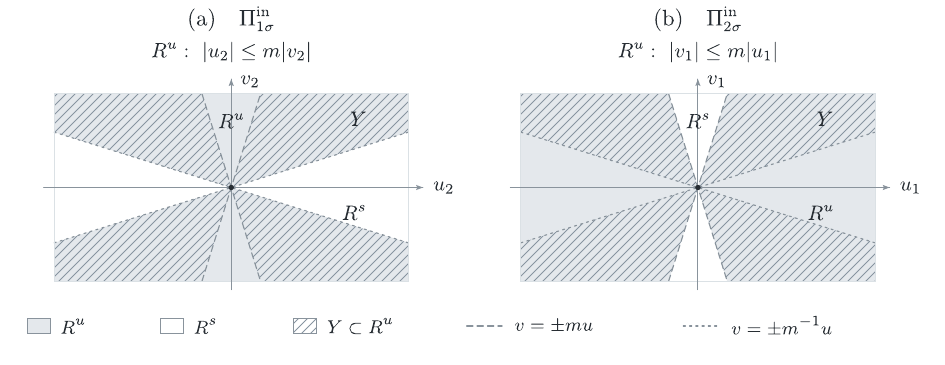}
\caption{The regions $R^u$, $R^s$, and $Y$ on the incoming sections, for either sign $\sigma$. Gray indicates $R^u$, white indicates $R^s$, and hatching indicates $Y\subset R^u$.}
\label{fig:regions_on_cross_sections}
\end{figure}

To analyze the Poincar\'e map $T$ on the domain $D$, it is easier if we divide $D$ into smaller regions, and analyze $T$ on each region separately. Accordingly, we fix a sufficiently large $m >1$ such that
\begin{equation}\label{eq:range_of_m}
\left\{\frac{\left|d_1\right|}{\left|b_1\right|},
\frac{\left|d_1\right|}{\left|c_1\right|},
\frac{\left|d_2\right|}{\left|b_2\right|},
\frac{\left|c_2\right|}{\left|a_2\right|},
\frac{\left|b_2\right|}{\left|a_2\right|}\right\}\subset (m^{-1}, m),
\end{equation}
where $a_{i}$, $b_{i}$, $c_{i}$ and $d_{i}$ are as in \eqref{eq_Taylor_expansion_global_map}, and define
\begin{equation}\label{eq:definition_regions}
\begin{aligned}
R^{u} := \{(u_{2}, v_{2})\in \Pi^{\mathrm{in}}_{1}: \lvert u_{2}\rvert \leq m \lvert v_{2}\rvert\}  &\,\textstyle\bigcup\, \{(u_{1}, v_{1})\in \Pi^{\mathrm{in}}_{2}: \lvert v_{1}\rvert \leq m \lvert u_{1}\rvert\},\\[2pt]
R^{s} := \{(u_{2}, v_{2})\in \Pi^{\mathrm{in}}_{1}: \lvert u_{2}\rvert > m \lvert v_{2}\rvert\} &\,\textstyle\bigcup \, \{(u_{1}, v_{1})\in \Pi^{\mathrm{in}}_{2+} \cup \Pi^{\mathrm{in}}_{2-}: \lvert v_{1}\rvert > m \lvert u_{1}\rvert\},\\[2pt]
Y := \{(u_{2}, v_{2})\in \Pi^{\mathrm{in}}_{1}: m^{-1} \lvert v_{2} \rvert \leq \lvert u_{2}\rvert \leq m \lvert v_{2}\rvert\} &\,\textstyle\bigcup\, \{(u_{1}, v_{1})\in \Pi^{\mathrm{in}}_{2}: m^{-1} \lvert u_{1}\rvert \leq \lvert v_{1}\rvert \leq m \lvert u_{1}\rvert\},
\end{aligned}
\end{equation}
where $\Pi^{\mathrm{in}}_{i} := \Pi^{\mathrm{in}}_{i+} \cup \Pi^{\mathrm{in}}_{i-}$ (see Figure \ref{fig:regions_on_cross_sections}). We further define $D^{u} := D\cap R^{u}$ and $D^{s} := D\cap R^{s}$.

In what follows, we analyze the local maps $T^{\mathrm{loc}}_{i\sigma_{0} j\sigma_{1}}$ for different values of $i,j =1,2$, and $\sigma_{0}, \sigma_{1}= \pm$.

\paragraph{Analysis of $T^{\mathrm{loc}}_{1\sigma_{0} 1\sigma_{1}}$.}
A point $(u_{20}, v_{20}) \in \Pi_{1\sigma_{0}}^{\mathrm{in}}$ reaches $(u_{2\tau}, v_{2\tau}) \in \Pi_{1\sigma_{1}}^{\mathrm{out}}$, where $\sigma_{0},\sigma_{1} = \pm$. Evaluating \eqref{eq_Poincare_map} at $u_{10} = \sigma_{0}\delta$ and $v_{1\tau} = \sigma_{1} \delta$ gives
\begin{align}
v_{10} &= e^{-\lambda_{1}\tau} \sigma_{1} \delta [1 + O(\delta^{2})],\label{eq:34iq5ey1}\\
v_{2\tau} &= e^{\lambda_{2}\tau} v_{20} [1 + O(\delta^{2})] + e^{(\lambda_{2} - 2\lambda_{1})\tau} O(\delta^{4}).\label{eq:34iq5ey2}
\end{align}
Choose $\delta>0$ small enough that both factors $1+O(\delta^2)$
in \eqref{eq:34iq5ey1} and \eqref{eq:34iq5ey2} have absolute value at least $1/2$.
Fix such a $\delta$ and fix $\epsilon_{\mathrm{out}}$. Let
$0<\epsilon\leq\epsilon_{\mathrm{out}}$, and suppose that a point
$(u_{20},v_{20})$ in the $\epsilon$-ball on
$\Pi^{\mathrm{in}}_{1\sigma_0}$ reaches the
$\epsilon_{\mathrm{out}}$-ball on
$\Pi^{\mathrm{out}}_{1\sigma_1}$ at time $\tau>0$.

By \eqref{eq56yigkre4r},
\begin{equation*}
v_{10} = \frac{\kappa}{\gamma}\sigma_0\delta^{-1}u_{20}v_{20}[1+O(\delta^2)] + O(\delta v_{20}^{2}).
\end{equation*}
Since $|u_{20}|,|v_{20}|\leq\epsilon$, there is a constant
$K>0$, independent of $\epsilon$ and $\tau$, such that
\begin{equation*}
|v_{10}|\leq K  \delta^{-1} \epsilon |v_{20}|.
\end{equation*}
On the other hand, \eqref{eq:34iq5ey1}
gives $|v_{10}|\geq \frac{\delta}{2}e^{-\lambda_1\tau}$. Consequently,
\begin{equation*}
e^{\lambda_1\tau}|v_{20}| \geq \frac{\delta^2}{2K\epsilon}.
\end{equation*}

Using \eqref{eq:34iq5ey2}, we therefore
obtain, for some constant $K^{\prime}>0$ independent of $\epsilon$ and
$\tau$,
\begin{align*}
|v_{2\tau}| &\geq \frac{1}{2}e^{\lambda_2\tau}|v_{20}| - K^{\prime}\delta^4e^{(\lambda_2-2\lambda_1)\tau} \\
&\geq \frac{\delta^2}{4 K\epsilon} e^{(\lambda_2-\lambda_1)\tau} - K^{\prime}\delta^4 \\
&\geq \frac{\delta^2}{4K\epsilon}-K^{\prime}\delta^4.
\end{align*}
Here we used $\lambda_1<\lambda_2<2\lambda_1$ and $\tau>0$.
For sufficiently small $\epsilon$, the last expression exceeds
$\epsilon_{\mathrm{out}}$. This contradicts
$\|(u_{2\tau},v_{2\tau})\|_2\leq\epsilon_{\mathrm{out}}$.
Thus no such transition is possible, and we obtain the following
lemma.

\begin{mylem}\label{Lem_565ryftuyg}
For sufficiently small $\epsilon$,
$D_{1\sigma_0 1\sigma_1}=\emptyset$ for all
$\sigma_0,\sigma_1=\pm$.
\end{mylem}

\paragraph{Analysis of $T^{\mathrm{loc}}_{2\sigma_{0} 2\sigma_{1}}$.}

\begin{figure}[!tbh]
\centering
\includegraphics[width=0.9\linewidth]{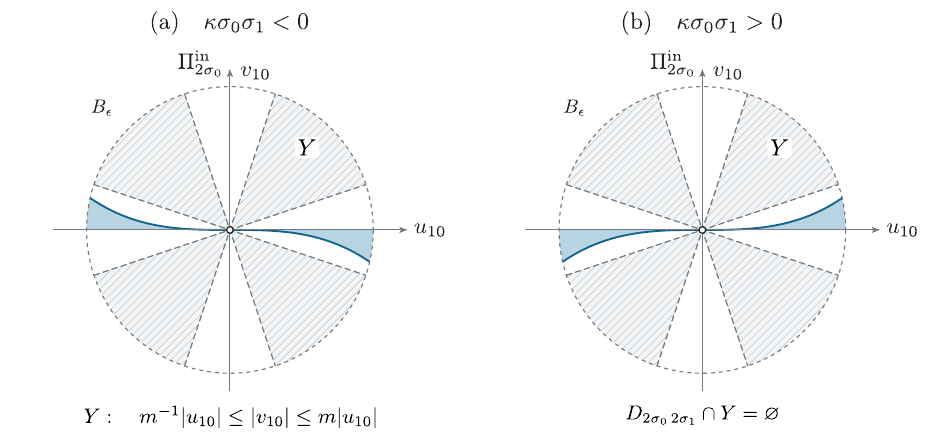}
\caption{Schematic geometry of $D_{2\sigma_0 2\sigma_1}\cap R^u$ for
a sufficiently small entrance radius $\epsilon$, with $\delta$ and
$\epsilon_{\mathrm{out}}$ fixed. The domain is contained in the blue
regions; the gray hatched sectors represent $Y$. The dashed straight
lines are $|v_{10}|=m^{-1}|u_{10}|$ and $|v_{10}|=m|u_{10}|$.
The blue regions lie strictly between the horizontal axis and the
shallow cone boundaries, showing that
$D_{2\sigma_0 2\sigma_1}\cap Y=\varnothing$.}
\label{Fig172yuihgoyvf}
\end{figure}

A point $(u_{10}, v_{10}) \in \Pi_{2\sigma_{0}}^{\mathrm{in}}$ reaches a point $(u_{1\tau}, v_{1\tau}) \in \Pi_{2\sigma_{1}}^{\mathrm{out}}$, where $\sigma_{0},\sigma_{1} = \pm$. Evaluating \eqref{eq_Poincare_map} at $u_{20} = \sigma_{0}\delta$ and $v_{2\tau} = \sigma_{1} \delta$ gives
\begin{align}
u_{1\tau} &= e^{-\lambda_{1}\tau} u_{10} [1 + O(\delta^{2})] + e^{-\lambda_{2} \tau} O(\lvert v_{1\tau}\rvert \delta^{2}),\nonumber\\
v_{1\tau} &= e^{\lambda_{1}\tau} v_{10} [1 + O(\delta^{2})] + e^{(\lambda_{1} -\lambda_{2})\tau} O(\lvert u_{10}\rvert \delta^{2}),\nonumber\\
v_{20} &= e^{-\lambda_{2}\tau} \sigma_{1}\delta [1 + O(\delta^{2})].\label{eq0987yguvjhbj}
\end{align}
Note also that \eqref{eq56yigkre4r} yields
\begin{equation}\label{eq56ygt75iygkjh}
v_{20} = \kappa \gamma \sigma_{0} \delta^{-1} u_{10} v_{10} [1 + O(\delta^{2})] + O\left(\delta v_{10}^{2}\right).
\end{equation}
For our purposes, it is sufficient to investigate $D_{2\sigma_{0} 2\sigma_{1}} \cap R^{u}$. Let $(u_{10}, v_{10}) \in R^{u}$, so $v_{10} = O(\lvert u_{10}\rvert)$. By virtue of this relation, \eqref{eq56ygt75iygkjh} yields $v_{20} = \kappa\gamma \sigma_{0} \delta^{-1} u_{10} v_{10} \left[1 + O\left(\delta^{2}\right)\right]$. Substituting this into \eqref{eq0987yguvjhbj} gives
\begin{equation*}
e^{-\lambda_{2}\tau} = \gamma\delta^{-2}\lvert u_{10} v_{10}\rvert [1 + O(\delta^{2})]
\qquad \mathrm{and} \qquad
e^{-\lambda_{1}\tau} = \gamma^{\gamma} \delta^{-2\gamma}\lvert u_{10} v_{10}\rvert^{\gamma} [1 + O(\delta^{2})].
\end{equation*}
Therefore, for $(u_{10}, v_{10}) \in R^{u}$, we obtain $v_{1\tau} = e^{\lambda_{1}\tau} v_{10}[1 + O(\delta^{2})]$, and so
\begin{equation*}
v_{1\tau} = \gamma^{-\gamma} \delta^{2\gamma}\lvert u_{10}\rvert^{-\gamma} v_{10} \lvert v_{10}\rvert^{-\gamma} [1 + O(\delta^{2})].
\end{equation*}
Fixing $\delta > 0$, and taking $\epsilon$ small enough, we have
\begin{equation*}
\frac{u_{1\tau}}{v_{1\tau}} = e^{-2\lambda_{1}\tau} \frac{u_{10}}{v_{10}}[1 + O(\delta^{2})] + e^{-\lambda_{2} \tau} O(\delta^{2})
= O(\epsilon^{2}).
\end{equation*}
Thus, $u_{1\tau} = O(\lvert v_{1\tau}\rvert \epsilon^{2})$ and so, for $(u_{10}, v_{10}) \in R^{u}$, we have
\begin{equation*}
\|(u_{1\tau}, v_{1\tau})\| = \lvert v_{1\tau}\rvert [1 + O(\epsilon^{2})] = \gamma^{-\gamma} \lvert u_{10}\rvert^{-\gamma} \lvert v_{10}\rvert^{1-\gamma}\delta^{2\gamma}[1+ O(\delta^{2})].
\end{equation*}
Noting $\|(u_{1\tau}, v_{1\tau})\| \leq \epsilon_{\mathrm{out}}$ and the condition $\mathrm{sgn}(u_{10} v_{10}) = \kappa \sigma_{0} \sigma_{1}$ (see \eqref{eq0987yguvjhbj} and \eqref{eq56ygt75iygkjh}), the domain is characterized by curves tangent to the horizontal axis at $(0,0)$ (see the blue region in Figure \ref{Fig172yuihgoyvf}).

\begin{myrem}
Regarding the map $T_{2\sigma_{0} 2\sigma_{1}}$, for $(u_{10}, v_{10})\in D_{2\sigma_{0} 2\sigma_{1}} \cap R^{u}$, by virtue of $u_{1\tau} = O(\lvert v_{1\tau}\rvert \epsilon^{2})$ and \eqref{eq_Taylor_expansion_global_map}, we have
\begin{equation}\label{eq98qq7yewqj}
T_{2\sigma_{0} 2\sigma_{1}}(u_{10}, v_{10}) = \left([b_{2} + O(\epsilon^{2}) + O(\epsilon_{\mathrm{out}}^{2})] v_{1\tau},\, [d_{2} + O(\epsilon^{2}) + O(\epsilon_{\mathrm{out}}^{2})] v_{1\tau}\right).
\end{equation}
\end{myrem}

\paragraph{Analysis of $T^{\mathrm{loc}}_{1\sigma_{0} 2\sigma_{1}}$.}
A point $(u_{20}, v_{20})\in \Pi^{\mathrm{in}}_{1\sigma_{0}}$ reaches $(u_{1\tau}, v_{1\tau})\in \Pi^{\mathrm{out}}_{2\sigma_{1}}$, where $\sigma_{0},\sigma_{1} = \pm$. Evaluating \eqref{eq_Poincare_map} at $u_{10}=\sigma_{0}\delta$ and $v_{2\tau} = \sigma_{1}\delta$ gives
\begin{align}
u_{1\tau} &= e^{-\lambda_{1}\tau} \sigma_{0}\delta [1 + O(\delta^{2})],\label{eq2d04gouhei}\\
v_{10} &= e^{-\lambda_{1} \tau} v_{1\tau} [1 + O(\delta^{2})] + e^{-\lambda_{2} \tau} O(\lvert u_{20}\rvert \delta^{2}) + e^{-2 \lambda_{2} \tau} O(\delta^{4}), \label{eq028uiuktjo}\\
v_{20} &= e^{-\lambda_{2} \tau} \sigma_{1}\delta [1 + O(\delta^{2})]. \label{eq457it7y8}
\end{align}
From \eqref{eq457it7y8}, we obtain
\begin{gather}
\mathrm{sgn}(v_{20}) = \sigma_{1},\label{eq_4576rutiyku}\\
e^{-\lambda_{2}\tau} = \sigma_{1} v_{20} \delta^{-1} [1 + O(\delta^{2})]
\qquad \mathrm{and}\qquad
e^{-\lambda_{1}\tau} = (\sigma_{1} v_{20})^{\gamma}\delta^{-\gamma} [1 + O(\delta^{2})].\label{eq567ug79osw}
\end{gather}
However, \eqref{eq56yigkre4r} yields $v_{10} = \frac{\kappa}{\gamma} \sigma_{0} \delta^{-1} u_{20} v_{20} [1 + O(\delta^{2})] + O(\delta v_{20}^{2})$. This together with \eqref{eq2d04gouhei}, \eqref{eq028uiuktjo} and \eqref{eq567ug79osw} gives
\begin{align}
u_{1\tau} &= \sigma_{0} \lvert v_{20}\rvert^{\gamma} \delta^{1-\gamma} [1 + O(\delta^{2})],\label{eq_45tu9yo}\\
v_{1\tau} &= \kappa \gamma^{-1} \sigma_{0}\sigma_{1} u_{20} \lvert v_{20}\rvert^{1-\gamma} \delta^{\gamma - 1}[1 + O(\delta^{2})] + O(\lvert v_{20}\rvert^{2-\gamma} \delta^{1 + \gamma}).\nonumber
\end{align}

\begin{figure}[!tbh]
\centering
\includegraphics[width=\linewidth]{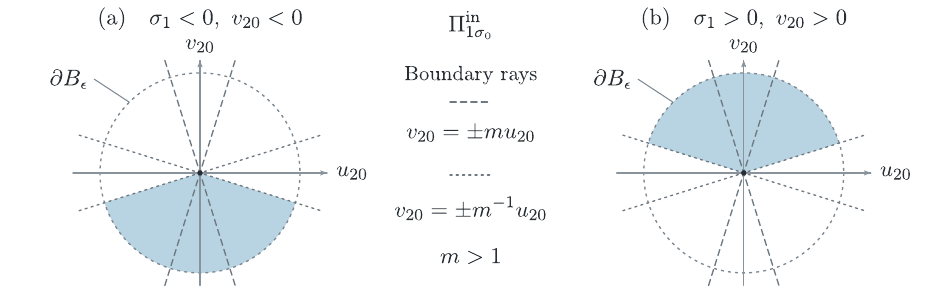}
\caption{The intersection $D_{1\sigma_0 2\sigma_1}\cap R^u$ on $\Pi^{\mathrm{in}}_{1\sigma_0}$ is shown in blue: \textup{(a)} $\sigma_{1}=-$, so $v_{20}<0$, and \textup{(b)} $\sigma_{1}=+$, so $v_{20}>0$. The outer circles are $\partial B_\epsilon$.}
\label{fig:domain_D12}
\end{figure}

Let us now investigate $D_{1\sigma_{0} 2\sigma_{1}} \cap R^{u}$. If $(u_{20}, v_{20}) \in R^{u}$, then $u_{20} = O(\lvert v_{20}\rvert)$. We then have
\begin{equation}\label{eq_67ihbu3y4g4byu}
\frac{v_{1\tau}}{u_{1\tau}} = \frac{\kappa}{\gamma} \sigma_{1} u_{20} \lvert v_{20}\rvert^{1-2\gamma} \delta^{2\gamma - 2} [1 + O(\delta^{2})] + O(\lvert v_{20}\rvert^{2-2\gamma} \delta^{2\gamma}) = O(\epsilon^{2 -2\gamma}),
\end{equation}
and so $v_{1\tau} = O(\epsilon^{2-2\gamma} \lvert u_{1\tau}\rvert)$. Thus, choosing $\epsilon > 0$ small enough, we have $\|(u_{1\tau}, v_{1\tau})\| = \lvert v_{20}\rvert^{\gamma} \delta^{1-\gamma} [1 + O(\delta^{2})]$ for $(u_{20}, v_{20}) \in R^{u}$. Therefore, $\|(u_{1\tau}, v_{1\tau})\| \leq \epsilon_{\mathrm{out}}$ if $\epsilon > 0$ is sufficiently small. Figure \ref{fig:domain_D12} demonstrates $D_{1\sigma_{0} 2\sigma_{1}} \cap R^{u}$ for $\epsilon >0$ small enough.

\begin{myrem}
For later reference, we remark that for $(u_{20}, v_{20})\in Y$, we have
\begin{equation}\label{eq_67iby20iubhv}
v_{1\tau} = \kappa \gamma^{-1} \sigma_{0}\sigma_{1} u_{20} \lvert v_{20}\rvert^{1-\gamma} \delta^{\gamma - 1} [1 + O(\delta^{2})].
\end{equation}
Regarding the map $T_{1\sigma_{0} 2\sigma_{1}}$, by \eqref{eq_Taylor_expansion_global_map} and for when $(u_{20}, v_{20}) \in R^{u}$, we have
\begin{equation}\label{eq096te6dft}
\begin{aligned}
T_{1\sigma_{0} 2\sigma_{1}}(u_{20}, v_{20}) &= \left([a_{2} + O(\epsilon^{2-2\gamma})] u_{1\tau},\, [c_{2} + O(\epsilon^{2-2\gamma})] u_{1\tau}\right)\\
&= \left(a_{2} \sigma_{0} \lvert v_{20}\rvert^{\gamma} \delta^{1-\gamma} [1 + O(\delta^{2})],\, c_{2} \sigma_{0} \lvert v_{20}\rvert^{\gamma} \delta^{1-\gamma} [1 + O(\delta^{2})]\right).
\end{aligned}
\end{equation}
\end{myrem}

\paragraph{Analysis of $T^{\mathrm{loc}}_{2\sigma_{0} 1\sigma_{1}}$.}
A point $(u_{10}, v_{10})\in \Pi^{\mathrm{in}}_{2\sigma_{0}}$ reaches $(u_{2\tau}, v_{2\tau})\in \Pi^{\mathrm{out}}_{1\sigma_{1}}$, where $\sigma_{0}, \sigma_{1} = \pm$. Evaluating \eqref{eq_Poincare_map} at $u_{20}=\sigma_{0}\delta$ and $v_{1\tau} = \sigma_{1}\delta$ gives
\begin{align}
u_{1\tau} &= e^{-\lambda_{1}\tau} u_{10} [1 + O(\delta^{2})] + e^{-2 \lambda_{2}\tau} O(\delta^{4}) + e^{-\lambda_{2} \tau} O(\lvert v_{2\tau}\rvert \delta^{2}),\label{eq6yfigk8ihou}\\
u_{2\tau} &= e^{-\lambda_{2}\tau} \sigma_{0}\delta [1 + O(\delta^{2})],\label{eq43iouy8i}\\
v_{10} &= e^{-\lambda_{1}\tau} \sigma_{1}\delta [1 + O(\delta^{2})],\label{eq024rfghvgyi}\\
v_{20} &= e^{-\lambda_{2}\tau} v_{2\tau}[1 + O(\delta^{2})] + e^{-2 \lambda_{1} \tau} O(u_{10}^{2} \delta^{3}) + e^{-2 \lambda_{2} \tau} O(\delta^{4}) + e^{-3 \lambda_{1} \tau} O(\lvert u_{10}\rvert \delta^{3}).\label{eq23oihi2e}
\end{align}
By virtue of \eqref{eq024rfghvgyi}, we obtain
\begin{gather}
\mathrm{sgn}(v_{10}) = \sigma_{1},\label{eq_09728yugk}\\
e^{-\lambda_{1}\tau} = \sigma_{1} v_{10} \delta^{-1} [1 + O(\delta^{2})]
\qquad \mathrm{and}\qquad
e^{-\lambda_{2}\tau} = (\sigma_{1} v_{10})^{\frac{1}{\gamma}} \delta^{-{\frac{1}{\gamma}}} [1 + O(\delta^{2})].\label{eq9er578iuytry}
\end{gather}
Note also that \eqref{eq56yigkre4r} yields $v_{20} = \kappa\gamma \sigma_{0} \delta^{-1} u_{10} v_{10} [1 + O(\delta^{2})] + O\left(\delta v_{10}^{2}\right)$. This together with \eqref{eq43iouy8i}, \eqref{eq23oihi2e} and \eqref{eq9er578iuytry} gives
\begin{equation}\label{eq_178rtf6y3h23}
\begin{aligned}
u_{2\tau} &= \sigma_{0} \lvert v_{10}\rvert^{\frac{1}{\gamma}} \delta^{1 - \frac{1}{\gamma}} [1 + O(\delta^{2})],\\
v_{2\tau} &= \kappa\gamma\sigma_{0} \sigma_{1} u_{10} \lvert v_{10}\rvert^{1-\frac{1}{\gamma}} \delta^{\frac{1}{\gamma} -1} [1 + O(\delta^{2})] + O(\lvert v_{10}\rvert^{2 -\frac{1}{\gamma}} \delta^{1 + \frac{1}{\gamma}}).
\end{aligned}
\end{equation}

\begin{figure}[!tbh]
\centering
\includegraphics[width=\linewidth]{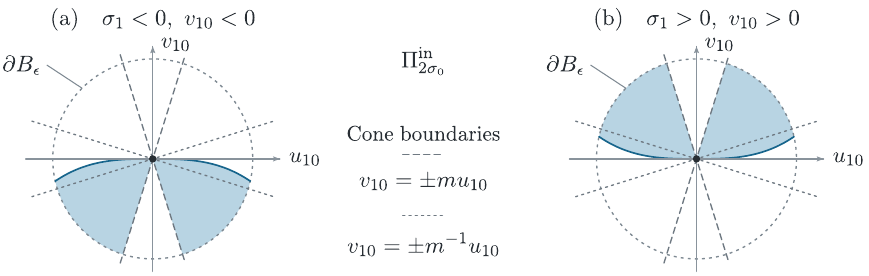}
\caption{Schematic of the domains $D_{2\sigma_0 1\sigma_1}\cap R^u$ on $\Pi^{\mathrm{in}}_{2\sigma_0}$. The blue region lies in $v_{10}<0$ in \textup{(a)} ($\sigma_1=-$) and in $v_{10}>0$ in \textup{(b)} ($\sigma_1=+$). The circle $\partial B_\epsilon$ bounds the entrance disk, and the solid blue curves mark the entrance points whose exit coordinates satisfy $\lVert(u_{2\tau},v_{2\tau})\rVert=\epsilon_{\mathrm{out}}$.}
\label{fig:domain_D21}
\end{figure}

Let us now investigate $D_{2\sigma_{0} 1\sigma_{1}} \cap R^{u}$. When $(u_{10}, v_{10})\in R^{u}$, we have $v_{10} = O(\lvert u_{10}\rvert)$, and so
\begin{equation}\label{eq_45tuvkb7taguer}
v_{2\tau} = \kappa\gamma\sigma_{0} \sigma_{1} u_{10} \lvert v_{10}\rvert^{1-\frac{1}{\gamma}} \delta^{\frac{1}{\gamma} -1} [1 + O(\delta^{2})].
\end{equation}
This gives
\begin{equation}\label{eq_98uriygvjbmdjer}
\frac{u_{2\tau}}{v_{2\tau}} = O\left(\epsilon^{\frac{2}{\gamma} -2} \delta^{2 - \frac{2}{\gamma}}\right).
\end{equation}
Thus, choosing $\epsilon>0$ small enough, we have $\|(u_{2\tau}, v_{2\tau})\| = \gamma\lvert u_{10}\rvert \lvert v_{10}\rvert^{1-\frac{1}{\gamma}} \delta^{\frac{1}{\gamma} -1} [1 + O(\delta^{2})]$. However, we must have $\|(u_{2\tau}, v_{2\tau})\| \leq \epsilon_{\mathrm{out}}$ and $\mathrm{sgn}(v_{10}) = \sigma_{1}$. Figure \ref{fig:domain_D21} demonstrates $D_{2\sigma_{0} 1\sigma_{1}} \cap R^{u}$ for $\epsilon > 0$ small enough.

For later reference, we remark that when $(u_{10}, v_{10}) \in D_{2\sigma_{0} 1\sigma_{1}} \cap R^{u}$, the map $T_{2\sigma_{0} 1\sigma_{1}}$ takes the form
\begin{equation}\label{eq017y3te569n}
\begin{aligned}
T_{2\sigma_{0} 1\sigma_{1}} (u_{10}, v_{10}) &= \left([b_{1} + O(\epsilon^{\frac{2}{\gamma}-2}) + O(\epsilon_{\mathrm{out}}^{2})] v_{2\tau}, [d_{1} + O(\epsilon^{\frac{2}{\gamma}-2}) + O(\epsilon_{\mathrm{out}}^{2})] v_{2\tau}\right)\\
&= \Big(b_{1} \kappa\gamma\sigma_{0} \sigma_{1} u_{10} \lvert v_{10}\rvert^{1-\frac{1}{\gamma}} \delta^{\frac{1}{\gamma} -1} [1 + O(\delta^{2})],\\
&\qquad d_{1} \kappa\gamma\sigma_{0} \sigma_{1} u_{10} \lvert v_{10}\rvert^{1-\frac{1}{\gamma}} \delta^{\frac{1}{\gamma} -1} [1 + O(\delta^{2})]\Big).
\end{aligned}
\end{equation}

\section{Dynamics and switching near the homoclinic loops}\label{sec:Dynamics_near_four_loops}

In this section, we study the dynamics near the homoclinic orbits and prove Theorems \ref{thm:dynamic_near_loops}, \ref{thm:coding_unstable} and \ref{thm:coding_stable}. This is done by analyzing the Poincar\'e maps on the cross-sections. To this end, we first investigate the dynamics of the Poincar\'e maps on the regions $R^{u}$ and $R^{s}$, separately. In particular, we show that the forward (resp. backward) iterations of the Poincar\'e maps keep $R^{u}$ (resp. $R^{s}$) invariant and eventually expand it. We then investigate the dynamics on the whole domain of the Poincar\'e maps, and prove the Theorems.

To start with, it would be beneficial if we introduce a uniform coordinate system on the cross-sections. Consider a point $(u_{10}, u_{20}, v_{10}, v_{20})$ on one of the four cross-sections. We identify this point by $(x,y)$, where $x= u_{20}$ and $y = v_{20}$ if the point belongs to $\Pi^{\mathrm{in}}_{1-} \cup \Pi^{\mathrm{in}}_{1+}$, and $x= u_{10}$ and $y = v_{10}$ if the point belongs to $\Pi^{\mathrm{in}}_{2-} \cup \Pi^{\mathrm{in}}_{2+}$.

Let $(x_{0}, y_{0}) \in D \cap \Pi^{\mathrm{in}}_{i_{0} \sigma_{0}}$, where $i_{0}\in\{1,2\}$ and $\sigma_{0} \in \{-,+\}$. There exists a Poincar\'e map $T_{i_{0}\sigma_{0} i_{1}\sigma_{1}}$ that takes $(x_{0}, y_{0})$ to $\Pi^{\mathrm{in}}_{i_{1} \sigma_{1}}$ for some $i_{1}\in\{1,2\}$ and $\sigma_{1} \in \{-,+\}$. Let $(x_{1}, y_{1}) := T_{i_{0}\sigma_{0} i_{1}\sigma_{1}}(x_{0}, y_{0})$. Similarly, if $(x_{1}, y_{1})\in D$, i.e., it lies in the domain of one of the Poincar\'e maps $T_{i_{1}\sigma_{1} i_{2}\sigma_{2}}$ for some $i_{2}$ and $\sigma_{2}$, we let $(x_{2}, y_{2}) := T_{i_{1}\sigma_{1} i_{2}\sigma_{2}}(x_{1}, y_{1}) \in \Pi^{\mathrm{in}}_{i_{2} \sigma_{2}}$. This point can be considered as the second iteration of $(x_{0}, y_{0})$. As long as the point $(x_{i}, y_{i})$ lies in $D$, we can continue this procedure and define the next iteration $(x_{i+1}, y_{i+1})$. Analogously, we can define the backward iterations of a point, too. We say that $(x_{-1}, y_{-1}) \in D \cap \Pi^{\mathrm{in}}_{i_{-1} \sigma_{-1}}$ for some $i_{-1}\in\{1,2\}$ and $\sigma_{-1} \in \{-,+\}$ is the backward iteration of $(x_{0}, y_{0}) \in \Pi^{\mathrm{in}}_{i_{0} \sigma_{0}}$ if $T_{i_{-1} \sigma_{-1} i_{0} \sigma_{0}}(x_{-1}, y_{-1}) = (x_{0}, y_{0})$.

We start our analysis with the following lemma. It states that $R^u$ is invariant under the forward dynamics.
\begin{mylem}\label{Lem_57fugip9}
If $(x_{0}, y_{0}) \in D^{u}$, then $(x_{1}, y_{1})\in Y \subset R^{u}$.
\end{mylem}

\begin{proof}
Let $(x_{0}, y_{0}) \in \Pi^{\mathrm{in}}_{1\sigma_{0}}$ for some $\sigma_{0}$. By Lemma \ref{Lem_565ryftuyg}, we must have $(x_{1}, y_{1}) \in \Pi^{\mathrm{in}}_{2\sigma_{1}}$ for some $\sigma_{1}$. Thus, $(x_{1}, y_{1}) = T_{1\sigma_{0} 2\sigma_{1}}(x_{0}, y_{0})$, and so by \eqref{eq096te6dft}, we have that $\frac{y_{1}}{x_{1}} = \frac{c_{2}}{a_{2}} + O(\epsilon^{2-2\gamma})$. Taking $\epsilon$ small enough, we have $(x_{1}, y_{1})\in Y$.

Now, suppose $(x_{0}, y_{0}) \in \Pi^{\mathrm{in}}_{2\sigma_{0}}$ for some $\sigma_{0}$. Then, there are two cases:
\begin{equation*}
(x_{1}, y_{1}) = T_{2\sigma_{0} 1\sigma_{1}}(x_{0}, y_{0})
\qquad \mathrm{or}\qquad
(x_{1}, y_{1}) = T_{2\sigma_{0} 2\sigma_{1}}(x_{0}, y_{0})
\end{equation*}
for some $\sigma_{1}$. In the first case, by \eqref{eq017y3te569n}, we have $\frac{y_{1}}{x_{1}} = \frac{d_{1}}{b_{1}} + O(\epsilon^{\frac{2}{\gamma} - 2}) + O(\epsilon_{\mathrm{out}}^{2})$, and in the second case, \eqref{eq98qq7yewqj} yields $\frac{y_{1}}{x_{1}} = \frac{d_{2}}{b_{2}} + O(\epsilon^{2}) + O(\epsilon_{\mathrm{out}}^{2})$. In both cases, noting \eqref{eq:range_of_m} and taking $\epsilon$ and $\epsilon_{\mathrm{out}}$ sufficiently small, we have $(x_{1}, y_{1})\in Y$, as desired.
\end{proof}

\begin{mycor}\label{cor:inavriance_of_Rs}
As a direct consequence of Lemma \ref{Lem_57fugip9}, we have that $R^s$ is invariant under the backward iterations. More precisely, if $(x_{0}, y_{0}) \in D$ and $(x_{1}, y_{1}) \in R^{s}$, then $(x_{0}, y_{0})\in D^{s} \subset R^s$.
\end{mycor}

Let us now explore the forward dynamics on $Y$ further. Since $Y \subset R^{u}$, Lemma \ref{Lem_57fugip9} implies that $Y$ is invariant under the forward iteration. Note also that the domain $D_{2\sigma_{0} 2\sigma_{1}}$ has an empty intersection with $Y$ (see Figure \ref{Fig172yuihgoyvf}). This means that if we start from $(x_{0}, y_{0}) \in Y \cap (\Pi^{\mathrm{in}}_{2-} \cup \Pi^{\mathrm{in}}_{2+})$, then $(x_{1}, y_{1})$ must lie in $Y \cap (\Pi^{\mathrm{in}}_{1-} \cup \Pi^{\mathrm{in}}_{1+})$. On the other hand, by Lemma \ref{Lem_565ryftuyg}, there is no possibility that $(x_{2}, y_{2})$ remains in $\Pi^{\mathrm{in}}_{1-} \cup \Pi^{\mathrm{in}}_{1+}$, and so it must return to $Y \cap (\Pi^{\mathrm{in}}_{2-} \cup \Pi^{\mathrm{in}}_{2+})$. Therefore, as long as the forward iterations of a point in $Y$ are defined, they switch between $\Pi^{\mathrm{in}}_{1-} \cup \Pi^{\mathrm{in}}_{1+}$ and $\Pi^{\mathrm{in}}_{2-} \cup \Pi^{\mathrm{in}}_{2+}$ alternately. We have
\begin{mylem}\label{lem:75rf2}
Let $(x_{0}, y_{0}) \in Y \cap D$. Then, $(x_{1}, y_{1}) \in Y$. Moreover, if $(x_{0}, y_{0}) \in \Pi^{\mathrm{in}}_{i\sigma_{0}}$ and $(x_{1}, y_{1}) \in \Pi^{\mathrm{in}}_{j\sigma_{1}}$, where $i,j\in\{1,2\}$ and $\sigma_{0}, \sigma_{1} = \pm$, then $i\neq j$.
\end{mylem}

Take a point $(x_{0}, y_{0}) \in Y$ and let the next few forward iterations of it exist (so, all the iterations lie in $Y$). Suppose $(x_{0}, y_{0})\in \Pi^{\mathrm{in}}_{1\sigma_{0}}$ for some $\sigma_{0} = \pm$. Then, for some $\sigma_{1}, \sigma_{2}, \sigma_{3}, \sigma_{4} = \pm$, we have
\begin{equation*}
(x_{0}, y_{0}) \xmapsto{T_{1\sigma_{0} 2\sigma_{1}}} (x_{1}, y_{1})
\xmapsto{T_{2\sigma_{1} 1\sigma_{2}}} (x_{2}, y_{2})
\xmapsto{T_{1\sigma_{2} 2\sigma_{3}}} (x_{3}, y_{3})
\xmapsto{T_{2\sigma_{3} 1\sigma_{4}}} (x_{4}, y_{4}).
\end{equation*}
Let us investigate how $\sigma_{i}$s evolve. Noting \eqref{eq_4576rutiyku} and \eqref{eq_09728yugk}, we obtain $\sigma_{1} = \mathrm{sgn}(y_{0})$, $\sigma_{2} = \mathrm{sgn}(y_{1})$, $\sigma_{3} = \mathrm{sgn}(y_{2})$, and $\sigma_{4} = \mathrm{sgn}(y_{3})$. We then have
\begin{enumerate}
\item By \eqref{eq096te6dft}, we have $\mathrm{sgn}(y_{1}) = \sigma_{0}\mathrm{sgn}(c_{2})$ and $\mathrm{sgn}(x_{1}) = \sigma_{0}\mathrm{sgn}(a_{2})$. So, $\sigma_{2} = \sigma_{0}\mathrm{sgn}(c_{2})$.
\item By \eqref{eq017y3te569n}, we have $\mathrm{sgn}(y_{2}) = \kappa \sigma_{1}\sigma_{2} \mathrm{sgn}(d_{1} x_{1}) = \kappa \mathrm{sgn}(d_{1} a_{2} c_{2} y_{0})$. So, $\sigma_{3} = \kappa \mathrm{sgn}(d_{1} a_{2} c_{2} y_{0})$.
\item By \eqref{eq096te6dft}, we have $\mathrm{sgn}(y_{3}) = \sigma_{2}\mathrm{sgn}(c_{2}) = \sigma_{0}$ and $\mathrm{sgn}(x_{3}) = \sigma_{2}\mathrm{sgn}(a_{2}) = \sigma_{0}\mathrm{sgn}(a_{2} c_{2})$. So, $\sigma_{4} = \sigma_{0}$.
\item By \eqref{eq017y3te569n}, we have $\mathrm{sgn}(y_{4}) = \kappa \sigma_{3} \sigma_{4} \mathrm{sgn}(d_{1} x_{3}) = \mathrm{sgn}(y_{0})$.
\end{enumerate}
In summary, $\sigma_{i}$s evolve as follows
\begin{equation}\label{eq_evolution_sigma}
\sigma_{0} \,\longrightarrow\,
\sigma_{1} = \mathrm{sgn}(y_{0}) \,\longrightarrow\,
\sigma_{2} = \sigma_{0}\mathrm{sgn}(c_{2}) \,\longrightarrow\,
\sigma_{3} = \kappa \mathrm{sgn}(d_{1} a_{2} c_{2} y_{0}) \,\longrightarrow\, \sigma_{4} =\sigma_{0}.
\end{equation}
This shows that starting from $Y\cap \Pi^{\mathrm{in}}_{1\sigma_{0}}$, after at most four iterations, we return to $Y\cap \Pi^{\mathrm{in}}_{1\sigma_{0}}$.

Throughout the paper, sign variables take values in $\{-1,1\}$. In labels and tables, we write $+$ and $-$ for $+1$ and $-1$, respectively. Products of signs are interpreted as ordinary multiplication.

Consider a point $(x_{0},y_{0}) \in Y$ and let some forward iterations of it exist. To each iteration $(x_{i}, y_{i})$, we assign the symbol ${k_{i}}_{\sigma_{i}}^{+}$ if $(x_{i}, y_{i})$ lies in the upper-half of $Y \cap \Pi^{\mathrm{in}}_{k_{i} \sigma_{i}}$, and ${k_{i}}_{\sigma_{i}}^{-}$ if it lies in the lower-half of $Y \cap \Pi^{\mathrm{in}}_{k_{i} \sigma_{i}}$. Putting these symbols in chronological order, we can then associate a string with (the forward orbit of) the point $(x_{0},y_{0})$. Indeed, suppose $(x_{0},y_{0}) \in Y \cap \Pi^{\mathrm{in}}_{1 \sigma_{0}}$, and let $A = \sgn(y_{0})$, $B = \sgn(c_{2})$ and $C = \sgn(\kappa d_{1} a_{2})$. By the analysis above, the string associated with $(x_{0}, y_{0})$ is
\begin{equation}\label{eq_234ijvtmbwjr}
{1}_{\sigma_{0}}^{A} \,\,
{2}_{\sigma_{1}}^{\sigma_{0} B} \,\,
{1}_{\sigma_{2}}^{ABC} \,\,
{2}_{\sigma_{3}}^{\sigma_{0}} \,\,
{1}_{\sigma_{4}}^{A}
\,\, \ldots,
\end{equation}
where, by \eqref{eq_evolution_sigma}, we have $\sigma_{1} = A$, $\sigma_{2} = \sigma_{0} B$, $\sigma_{3} = ABC$, and $\sigma_{4} = \sigma_{0}$. Accordingly, we can define a function that assigns to the symbol associated with a point in $Y$ the symbol associated with the forward iteration of that point.
\begin{mydefn}\label{dfn_evolution_function_Su}
Let $E := \{{1}_{+}^{+}, {1}_{+}^{-}, {1}_{-}^{+}, {1}_{-}^{-}, {2}_{+}^{+}, {2}_{+}^{-}, {2}_{-}^{+}, {2}_{-}^{-}\}$ be the set of all symbols. Given $B = \sgn(c_{2})$ and $C = \sgn(\kappa d_{1} a_{2})$, we define the bijection $S_{u} = S_{u; B,C}: E \rightarrow E$ as in Table \ref{table:bijection_Su}.
\end{mydefn}

\begin{table}[!htb]
\centering
\begingroup
\normalfont\normalsize
\setlength{\tabcolsep}{6pt}
\renewcommand{\arraystretch}{1.12}
\begin{tabular}{@{}cc@{\hspace{1.4em}}cccc@{\hspace{1.4em}}cccc@{}}
\toprule
\multicolumn{2}{c}{Signs}
& \multicolumn{4}{c}{Input $x={1}_{\sigma}^{q}$}
& \multicolumn{4}{c}{Input $x={2}_{\sigma}^{q}$} \\
\cmidrule(lr){1-2}\cmidrule(lr){3-6}\cmidrule(lr){7-10}
$B$ & $C$
& ${1}_{+}^{+}$ & ${1}_{+}^{-}$ & ${1}_{-}^{+}$ & ${1}_{-}^{-}$
& ${2}_{+}^{+}$ & ${2}_{+}^{-}$ & ${2}_{-}^{+}$ & ${2}_{-}^{-}$ \\
\midrule
$+$ & $+$
& ${2}_{+}^{+}$ & ${2}_{-}^{+}$ & ${2}_{+}^{-}$ & ${2}_{-}^{-}$
& ${1}_{+}^{+}$ & ${1}_{-}^{+}$ & ${1}_{+}^{-}$ & ${1}_{-}^{-}$ \\
$+$ & $-$
& ${2}_{+}^{+}$ & ${2}_{-}^{+}$ & ${2}_{+}^{-}$ & ${2}_{-}^{-}$
& ${1}_{+}^{-}$ & ${1}_{-}^{-}$ & ${1}_{+}^{+}$ & ${1}_{-}^{+}$ \\
\midrule
$-$ & $+$
& ${2}_{+}^{-}$ & ${2}_{-}^{-}$ & ${2}_{+}^{+}$ & ${2}_{-}^{+}$
& ${1}_{+}^{-}$ & ${1}_{-}^{-}$ & ${1}_{+}^{+}$ & ${1}_{-}^{+}$ \\
$-$ & $-$
& ${2}_{+}^{-}$ & ${2}_{-}^{-}$ & ${2}_{+}^{+}$ & ${2}_{-}^{+}$
& ${1}_{+}^{+}$ & ${1}_{-}^{+}$ & ${1}_{+}^{-}$ & ${1}_{-}^{-}$ \\
\bottomrule
\end{tabular}
\endgroup
\caption{The bijection $S_{u;B,C}:E\to E$.
For each sign pair $(B,C)$, the entry in column $x$ is $S_{u;B,C}(x)$.}
\label{table:bijection_Su}
\end{table}

The symbolic pattern in \eqref{eq_234ijvtmbwjr} can now be written as
\begin{equation*}
{1}_{\sigma_0}^{A} \,\,\, S_{u}({1}_{\sigma_0}^{A})\,\,\, S_{u}^{2}({1}_{\sigma_0}^{A})\,\,\, S_{u}^{3}({1}_{\sigma_0}^{A}) \,\,\, S_{u}^{4}({1}_{\sigma_0}^{A}) \,\, \ldots .
\end{equation*}
Furthermore, observe that
\begin{equation*}
S_{u}^{4} = \mathrm{id}_E.
\end{equation*}
In other words, the fifth and first elements of string \eqref{eq_234ijvtmbwjr} are the same which means that if we consider the block consisting of the first four elements of the string, then the whole string is made of repeating this block. We call the minimal block whose repetitions generate the string of the point $(x_{0},y_{0})$ the code of $(x_{0},y_{0})$.

So far, we have seen that the forward iterations in $Y$ are governed by the two types of the maps $T_{1\sigma_{0} 2\sigma_{1}}$ and $T_{2\sigma_{0} 1\sigma_{1}}$, for some $\sigma_{0}, \sigma_{1} = \pm$, alternately. The following lemma states that these two maps expand $R^{u}$, and so expand $Y$, too.
\begin{mylem}\label{Lem_66fiytfi5}
Let $(x_{0}, y_{0}) \in D^{u} \subset R^{u}$; so $(x_{0}, y_{0}) \in \Pi^{\mathrm{in}}_{1\sigma_{0}}$ or $(x_{0}, y_{0}) \in \Pi^{\mathrm{in}}_{2\sigma_{0}}$ for some $\sigma_{0} = \pm$.
\begin{enumerate}
\item Suppose $(x_{0}, y_{0}) \in \Pi^{\mathrm{in}}_{1\sigma_{0}}$ and let $(x_{1}, y_{1}) := T_{1\sigma_{0} 2\sigma_{1}}(x_{0}, y_{0})$ for some $\sigma_{1} = \pm$. Then as $\epsilon \rightarrow 0$, we have $\|(x_{0}, y_{0})\| = o(\|(x_{1}, y_{1})\|)$ and $\frac{y_{1}}{x_{1}} \rightarrow \frac{c_{2}}{a_{2}}$.

\item Suppose $(x_{0}, y_{0}) \in \Pi^{\mathrm{in}}_{2\sigma_{0}}$ and let $(x_{1}, y_{1}) := T_{2\sigma_{0} 1\sigma_{1}}(x_{0}, y_{0})$ for some $\sigma_{1} = \pm$. Then as $(\epsilon, \epsilon_{\mathrm{out}}) \rightarrow (0,0)$, we have $\|(x_{0}, y_{0})\| = o(\|(x_{1}, y_{1})\|)$ and $\frac{y_{1}}{x_{1}} \rightarrow \frac{d_{1}}{b_{1}}$.
\end{enumerate}
\end{mylem}

\begin{proof}
To prove the first part, note that since $(x_{0}, y_{0}) \in D^{u} \cap\Pi^{\mathrm{in}}_{1\sigma_{0}} \subset R^{u}$, by definition, we have $\lvert x_{0}\rvert \leq m\lvert y_{0}\rvert$, and so $\|(x_{0}, y_{0})\| \leq \sqrt{m^{2} + 1} \lvert y_{0}\rvert$. Thus, by \eqref{eq096te6dft}, we have
\begin{equation*}
\frac{\|(x_{0}, y_{0})\|}{\|(x_{1}, y_{1})\|} \leq \lvert y_{0}\rvert^{1-\gamma} \delta^{\gamma -1} \sqrt{\frac{m^{2} + 1}{[a_{2}^{2} + c_{2}^{2}] [1+ O(\delta^2)]}}.
\end{equation*}
Taking into account that $\lvert y_{0}\rvert \leq \epsilon$, this relation yields $\|(x_{0}, y_{0})\| = o(\|(x_{1}, y_{1})\|)$ as $\epsilon \rightarrow 0$. It also follows from \eqref{eq096te6dft} that $\frac{y_{1}}{x_{1}} = \frac{c_{2}}{a_{2}} + O(\epsilon^{2-2\gamma})$ as $\epsilon \rightarrow 0$. This proves the first part of the lemma.

Regarding the second part, by the definition of $R^u$, relation $(x_{0}, y_{0}) \in D^{u} \cap\Pi^{\mathrm{in}}_{2\sigma_{0}} \subset R^{u}$ yields $\lvert y_{0}\rvert \leq m \lvert x_{0}\rvert$, and so $\|(x_{0}, y_{0})\| \leq \sqrt{m^{2} + 1} \lvert x_{0}\rvert$. Then, by \eqref{eq017y3te569n}, we have
\begin{equation*}
\frac{\|(x_{0}, y_{0})\|}{\|(x_{1}, y_{1})\|} \leq
\gamma^{-1} |y_{0}|^{\frac{1}{\gamma} - 1} \delta^{1 - \frac{1}{\gamma}} \sqrt{\frac{1 + m^2}{[b_{1}^{2} + d_{1}^{2}] [1 + O(\delta^{2})]}}.
\end{equation*}
However, $\lvert y_{0}\rvert \leq \epsilon$, and so $\|(x_{0}, y_{0})\| = o(\|(x_{1}, y_{1})\|)$ as $\epsilon \rightarrow 0$. Note also that $\frac{y_{1}}{x_{1}} = \frac{d_{1}}{b_{1}} + O(\epsilon^{\frac{2}{\gamma} - 2}) + O(\epsilon_{\mathrm{out}}^{2})$ by \eqref{eq017y3te569n}. This ends the proof.
\end{proof}

Let us now have a look at the dynamics in $R^s$. By Corollary \ref{cor:inavriance_of_Rs}, we have that $R^{s}$ is invariant under the backward iterations. The following lemma is an analogue of Lemma \ref{Lem_66fiytfi5}, and states that $T^{-1}_{1\sigma_{0} 2\sigma_{1}}$ and $T^{-1}_{2\sigma_{0} 1\sigma_{1}}$ expand $R^{s}$. Later, in Proposition \ref{prp:switching_backward_dynamics}, we show that the backward dynamics in $R^s$ is governed only by these two types of the maps.

\begin{mylem}\label{Lem_989giwykd}
Let $(x_{0}, y_{0}) \in D$ for which $(x_{1}, y_{1}) \in R^{s}$. Then, $(x_{0}, y_{0})$ is in $R^{s}$ as well, i.e. $(x_{0}, y_{0}) \in D^{s} \subset R^{s}$. There are two possibilities: $(x_{1}, y_{1}) \in \Pi^{\mathrm{in}}_{1\sigma_{1}}$ or $(x_{1}, y_{1}) \in \Pi^{\mathrm{in}}_{2\sigma_{1}}$. Then:
\begin{enumerate}
\item Suppose $(x_{1}, y_{1}) \in \Pi^{\mathrm{in}}_{1\sigma_{1}}$ and let $(x_{0}, y_{0}) = T^{-1}_{2\sigma_{0} 1\sigma_{1}}(x_{1}, y_{1})$. Then as $(\epsilon, \epsilon_\mathrm{out})\rightarrow (0,0)$, we have $\|(x_{1}, y_{1})\| = o(\|(x_{0}, y_{0})\|)$ and $|\frac{y_{0}}{x_{0}}| \rightarrow \infty$.
\item Suppose $(x_{1}, y_{1}) \in \Pi^{\mathrm{in}}_{2\sigma_{1}}$ and let $(x_{0}, y_{0}) = T^{-1}_{1\sigma_{0} 2\sigma_{1}}(x_{1}, y_{1})$. Then as $(\epsilon, \epsilon_\mathrm{out})\rightarrow (0,0)$, we have $\|(x_{1}, y_{1})\| = o(\|(x_{0}, y_{0})\|)$ and $\frac{y_{0}}{x_{0}} \rightarrow 0$.
\end{enumerate}
\end{mylem}

\begin{proof}
Take $i,j\in \{1,2\}$ such that $i\neq j$. Consider system \eqref{eq:normal_form} and let us call it system I. Take a point $z_{0}\in \Pi^{\mathrm{out}}_{i\sigma}$, and let $z_{1}$, $z_{2}$, and $z_{3}$ be the first, second, and third intersection point of the forward orbit of $z_{0}$ with any of the eight cross-sections, respectively. Thus,
\begin{equation*}
z_{0}\in \Pi^{\mathrm{out}}_{i\sigma} \longrightarrow
z_{1}\in \Pi^{\mathrm{in}}_{i\sigma} \longrightarrow
z_{2}\in \Pi^{\mathrm{out}}_{j\sigma^{\prime}}\longrightarrow
z_{3}\in \Pi^{\mathrm{in}}_{j\sigma^{\prime}} \qquad \mathrm{for\,\, some\,\,} \sigma,\sigma^{\prime} = \pm.
\end{equation*}
So, $z_{3} = T_{i\sigma j\sigma^{\prime}}(z_{1})$. Call system II the system obtained from system I by reversing time. For system II, we have
\begin{equation*}
z_{3}\in \Pi^{\mathrm{in}}_{j\sigma^{\prime}}\longrightarrow
z_{2}\in \Pi^{\mathrm{out}}_{j\sigma^{\prime}}\longrightarrow
z_{1}\in \Pi^{\mathrm{in}}_{i\sigma} \longrightarrow
z_{0}\in \Pi^{\mathrm{out}}_{i\sigma}.
\end{equation*}
The homoclinic orbits in system II leave and approach $O$ along $u$-axes and $v$-axes, respectively, so the Poincar\'e map $\overline{T}_{j\sigma^{\prime} i\sigma}$ of this system is defined on a subset of $\Pi^{\mathrm{out}}_{j\sigma^{\prime}}$ into $\Pi^{\mathrm{out}}_{i\sigma}$. In particular, $\overline{T}_{j\sigma^{\prime} i\sigma}(z_{2}) = z_{0}$. Noting that
\begin{equation}\label{eq:98uyhhyguh}
z_{2} = (T_{j}^{\mathrm{glo}})^{-1}(z_{3})
\qquad \text{and}\qquad
z_{0} = (T_{i}^{\mathrm{glo}})^{-1}(z_{1}),
\end{equation}
we have
\begin{equation*}
T_{i}^{\mathrm{glo}} \circ \overline{T}_{j\sigma^{\prime} i\sigma} \circ (T_{j}^{\mathrm{glo}})^{-1} (z_{3}) = z_{1}.
\end{equation*}
In other words, $T_{i}^{\mathrm{glo}} \circ \overline{T}_{j\sigma^{\prime} i\sigma} \circ (T_{j}^{\mathrm{glo}})^{-1} = T_{i\sigma j\sigma^{\prime}}^{-1}$.

Applying the change of coordinates $(\widetilde{u}_{1}, \widetilde{u}_{2}, \widetilde{v}_{1}, \widetilde{v}_{2}) = (v_{1}, v_{2}, u_{1}, u_{2})$ on system II, i.e. swapping the stable and unstable components, we obtain a new system which we call system III. Consider also the cross-sections for system III:
\begin{equation*}
\Sigma_{k\sigma}^{\mathrm{in}} = \{\widetilde{u}_{k} = \sigma\delta\} \cap \{H=0\}
\qquad \mathrm{and}\qquad
\Sigma_{k\sigma}^{\mathrm{out}} = \{\widetilde{v}_{k} = \sigma\delta\} \cap \{H=0\}
\end{equation*}
for any $k=1,2$ and $\sigma = \pm$. Then,
\begin{equation*}
Jz_{3} \in \Sigma_{j\sigma^{\prime}}^{\mathrm{out}} \longrightarrow
Jz_{2} \in \Sigma_{j\sigma^{\prime}}^{\mathrm{in}} \longrightarrow
Jz_{1} \in \Sigma_{i\sigma}^{\mathrm{out}} \longrightarrow
Jz_{0} \in \Sigma_{i\sigma}^{\mathrm{in}},
\end{equation*}
where $J = \left(\begin{smallmatrix}
0 & 1\\ 1 & 0
\end{smallmatrix}\right)$. The Poincar\'e map of system III maps $Jz_{2}$ to $Jz_{0}$. The global maps of system III are $Jz_{3} \mapsto Jz_{2}$ and $Jz_{1} \mapsto Jz_{0}$ given by $J (T^{\mathrm{glo}}_{j})^{-1}J^{-1}$ and $J (T^{\mathrm{glo}}_{i})^{-1}J^{-1}$, respectively. Note also that the differential of the global maps at $(0,0)$ are
\begin{equation}\label{eq:79ugvjfykuci}
\frac{1}{a_{1} d_{1} - b_{1} c_{1}} \left(\begin{matrix}
a_{1} & -c_{1}\\
-b_{1} & d_{1}
\end{matrix}\right)
\qquad \text{and}\qquad
\frac{1}{a_{2} d_{2} - b_{2} c_{2}} \left(\begin{matrix}
a_{2} & -c_{2}\\
-b_{2} & d_{2}
\end{matrix}\right),
\end{equation}
respectively.

The key observation is that system III is of the form of system \eqref{eq:normal_form}, and so all the statements that we have proved for system \eqref{eq:normal_form} hold for system III as well. This resembles time reversibility, even though system \eqref{eq:normal_form} is not necessarily time-reversible. Following this observation, if we define
\begin{equation*}
\begin{aligned}
{\mathcal{R}}^{u} &:= \{(\widetilde{u}_{2}, \widetilde{v}_{2})\in \Sigma^{\mathrm{in}}_{1+} \cup \Sigma^{\mathrm{in}}_{1-}: \lvert \widetilde{u}_{2}\rvert \leq m \lvert \widetilde{v}_{2}\rvert\} \cup \{(\widetilde{u}_{1}, \widetilde{v}_{1})\in \Sigma^{\mathrm{in}}_{2+} \cup \Sigma^{\mathrm{in}}_{2-}: \lvert \widetilde{v}_{1}\rvert \leq m \lvert \widetilde{u}_{1}\rvert\},\\
\mathcal{Y} &:= \{(\widetilde{u}_{2}, \widetilde{v}_{2})\in \Sigma^{\mathrm{in}}_{1+} \cup \Sigma^{\mathrm{in}}_{1-}: m^{-1} \lvert \widetilde{v}_{2}\rvert \leq \lvert \widetilde{u}_{2}\rvert \leq m \lvert \widetilde{v}_{2}\rvert\}\\
&\qquad \cup \{(\widetilde{u}_{1}, \widetilde{v}_{1})\in \Sigma^{\mathrm{in}}_{2+} \cup \Sigma^{\mathrm{in}}_{2-}: m^{-1} \lvert \widetilde{u}_{1}\rvert \leq \lvert \widetilde{v}_{1}\rvert \leq m \lvert \widetilde{u}_{1}\rvert\},
\end{aligned}
\end{equation*}
which are analogous to $R^{u}$ and $Y$ in \eqref{eq:definition_regions}, but defined for system III, applying Lemma \ref{Lem_57fugip9} on system III yields that if $Jz_{2} \in {\mathcal{R}}^{u}$, then $Jz_{0} \in \mathcal{Y}$. Moreover, writing $J z_{0} = (\widetilde{u}_{j}, \widetilde{v}_{j})$ and taking \eqref{eq:79ugvjfykuci} into account, Lemma \ref{Lem_66fiytfi5} implies that as $(\epsilon, \epsilon_\mathrm{out})\rightarrow (0,0)$, we have $\|z_{2}\| = o(\|z_{0}\|)$ and $\frac{\widetilde{v}_{j}}{\widetilde{u}_{j}} \rightarrow s$, where $s= \frac{-b_{2}}{a_{2}}$ if $j=1$, and $s= \frac{-d_{1}}{c_{1}}$ if $j=2$. Consequently, for system II, if $z_{2}\in J^{-1}\mathcal{R}^{u}$ then $z_{0}\in J^{-1}\mathcal{Y}$. Moreover, writing $z_{0} = (u_{j}, v_{j})$, we have $\frac{v_{j}}{u_{j}} \rightarrow \frac{1}{s}$ as $(\epsilon, \epsilon_\mathrm{out})\rightarrow (0,0)$.

Repeating the same argument and considering \eqref{eq:98uyhhyguh}, we have that 
\begin{equation*}
\text{if}\quad
z_{3} = T^{\mathrm{glo}}_{j}(z_{2}) \in T^{\mathrm{glo}}_{j}(J^{-1}\mathcal{R}^{u})
\quad \text{then}\quad
z_{1} = T^{\mathrm{glo}}_{i}(z_{0}) \in T^{\mathrm{glo}}_{i}(J^{-1}\mathcal{Y}).
\end{equation*}
Moreover, since $T^{\mathrm{glo}}_{i}$ and $T^{\mathrm{glo}}_{j}$ are diffeomorphisms, as $(\epsilon, \epsilon_\mathrm{out})\rightarrow (0,0)$, we also have $\|z_{3}\| = o(\|z_{1}\|)$, and if we write $z_{1} = (u_{j}, v_{j})$, then $\frac{v_{j}}{u_{j}} \rightarrow p$, where $p = \pm\infty$ if $j=1$, and $p=0$ if $j=2$ (the latter follows from applying $DT^{\mathrm{glo}}_{i}$ on the line passing through $(0,0)$ with the slope $\frac{1}{s}$), which implies $z_{1} \in R^{s}$, and so $z_{1} \in D^{s} \subset R^{s}$. This means that if $R^{s} \subset T^{\mathrm{glo}}_{j}(J^{-1}\mathcal{R}^{u})$, then we are done with the proof of Lemma \ref{Lem_989giwykd}.

Observe that
\begin{equation*}
J^{-1} \mathcal{R}^{u} = \{(u_{2}, v_{2})\in \Pi^{\mathrm{out}}_{1-} \cup \Pi^{\mathrm{out}}_{1+}: \lvert v_{2}\rvert \leq m \lvert u_{2}\rvert\} \cup \{(u_{1}, v_{1})\in \Pi^{\mathrm{out}}_{2-} \cup \Pi^{\mathrm{out}}_{2+}: \lvert u_{1}\rvert \leq m \lvert v_{1}\rvert\}.
\end{equation*}
Choosing $m > 1$ sufficiently large, the set $J^{-1} \mathcal{R}^{u}$ covers $\Pi^{\mathrm{out}}_{1-} \cup \Pi^{\mathrm{out}}_{1+} \cup \Pi^{\mathrm{out}}_{2-} \cup \Pi^{\mathrm{out}}_{2+}$ except for two narrow cones containing the vertical axes on $\Pi^{\mathrm{out}}_{1-}$ and $\Pi^{\mathrm{out}}_{1+}$, and two narrow cones containing the horizontal axes on $\Pi^{\mathrm{out}}_{2-}$ and $\Pi^{\mathrm{out}}_{2+}$. Let $L:= D T_{i}^{\mathrm{glo}}(0,0)$. Then, $L$ for $i=1$ (resp. $i=2$) takes the vertical (resp. horizontal) cones on $\Pi^{\mathrm{out}}_{i-}$ and $\Pi^{\mathrm{out}}_{i+}$ to two narrow cones on $\Pi^{\mathrm{in}}_{i-}$ and $\Pi^{\mathrm{in}}_{i+}$ containing the straight lines with the slope $\frac{d_{1}}{b_{1}}$ (resp. $\frac{c_{2}}{a_{2}}$) that pass through $(0,0)$. This gives $R^{s} \subset L(J^{-1}\mathcal{R}^{u})$. However, $T^{\mathrm{glo}}_{j}$ is well-approximated by $L$ near $(0,0)$, and so $R^{s} \subset T^{\mathrm{glo}}_{j}(J^{-1}\mathcal{R}^{u})$. This ends the proof of Lemma \ref{Lem_989giwykd}.
\end{proof}

Earlier, we associated a symbol to any point in $Y$, and studied how these symbols evolve by time. Let us now do the same with the points in $R^s$. We associate the symbol ${1}_{\sigma}^{+}$ (resp. ${1}_{\sigma}^{-}$) to the points on $R^{s} \cap \Pi^{\mathrm{in}}_{1\sigma}$ that lie in the right (resp. left) half-plane. We also associate the symbol ${2}_{\sigma}^{+}$ (resp. ${2}_{\sigma}^{-}$) to the points on $R^{s} \cap \Pi^{\mathrm{in}}_{2\sigma}$ that lie in the upper (resp. lower) half-plane. We have shown that $R^{s}$ is invariant under the backward iterations of the Poincar\'e map. Let us now investigate how the symbols associated to the points in $R^s$ evolve under the backward dynamics. To do this, we follow the discussions on systems I, II, and III from the proof of Lemma \ref{Lem_989giwykd}.

Consider system I in the proof of Lemma \ref{Lem_989giwykd}. Take a point $z_{0} \in \Pi^{\mathrm{in}}_{1\sigma_{0}} \cap R^{s}$ for some $\sigma_{0} = \pm$. Let $z_{1}$,..., $z_{7}$ be the intersection points of the backward orbit of $z_{0}$ with any of the eight cross-sections (in the order encountered backward). In summary, we have
\begin{equation*}
\begin{tikzcd}[row sep=-1.1ex, column sep=1.6ex]
z_{7} & z_{6} & z_{5} & z_{4} & z_{3} & z_{2} & z_{1} & z_{0}\\
\indown & \indown & \indown & \indown & \indown & \indown & \indown & \indown\\
\Pi^{\mathrm{out}}_{i_{3}\sigma_{3}} 
& \Pi^{\mathrm{in}}_{i_{3}\sigma_{3}} \cap R^{s} 
& \Pi^{\mathrm{out}}_{i_{2}\sigma_{2}} 
& \Pi^{\mathrm{in}}_{i_{2}\sigma_{2}} \cap R^{s} 
& \Pi^{\mathrm{out}}_{i_{1}\sigma_{1}} 
& \Pi^{\mathrm{in}}_{i_{1}\sigma_{1}} \cap R^{s} 
& \Pi^{\mathrm{out}}_{1\sigma_{0}} 
& \Pi^{\mathrm{in}}_{1\sigma_{0}} \cap R^{s}
\arrow[from=1-1, to=1-2] \arrow[from=1-2, to=1-3] \arrow[from=1-3, to=1-4] \arrow[from=1-4, to=1-5] \arrow[from=1-5, to=1-6] \arrow[from=1-6, to=1-7] \arrow[from=1-7, to=1-8]
\end{tikzcd}
\end{equation*}
for some $\sigma_{0},\sigma_{1},\sigma_{2},\sigma_{3} = \pm$ and $i_{1}, i_{2}, i_{3} = 1,2$. Here, the arrow directions show how points move in forward time. Note also that the property $z_{k}\in R^{s}$ for $k=2,4,6$ is a consequence of Corollary \ref{cor:inavriance_of_Rs} and assumption $z_{0} \in R^{s}$. 

Let the symbols of $z_{0}$, $z_{2}$, $z_{4}$, and $z_{6}$ be ${1}_{\sigma_0}^{p_0}$, ${i_{1}}_{\sigma_1}^{p_1}$, ${i_{2}}_{\sigma_2}^{p_2}$, and ${i_{3}}_{\sigma_3}^{p_3}$, respectively, for some $p_{0}, p_{1}, p_{2}, p_{3} = \pm$. Our aim is to investigate the evolution ${i_{3}}_{\sigma_3}^{p_3} \mapsto {i_{2}}_{\sigma_2}^{p_2} \mapsto {i_{1}}_{\sigma_1}^{p_1} \mapsto {1}_{\sigma_0}^{p_0}$.

Reversing time, we obtain system II from system I. For system II, we have
\begin{equation*}
\begin{tikzcd}[row sep=-1.1ex, column sep=1.6ex]
z_{7} & z_{6} & z_{5} & z_{4} & z_{3} & z_{2} & z_{1} & z_{0}.\\
\indown & \indown & \indown & \indown & \indown & \indown & \indown & \indown\\
\Pi^{\mathrm{out}}_{i_{3}\sigma_{3}} 
& \Pi^{\mathrm{in}}_{i_{3}\sigma_{3}} \cap R^{s} 
& \Pi^{\mathrm{out}}_{i_{2}\sigma_{2}} 
& \Pi^{\mathrm{in}}_{i_{2}\sigma_{2}} \cap R^{s} 
& \Pi^{\mathrm{out}}_{i_{1}\sigma_{1}} 
& \Pi^{\mathrm{in}}_{i_{1}\sigma_{1}} \cap R^{s} 
& \Pi^{\mathrm{out}}_{1\sigma_{0}} 
& \Pi^{\mathrm{in}}_{1\sigma_{0}} \cap R^{s}
\arrow[from=1-2, to=1-1] \arrow[from=1-3, to=1-2] \arrow[from=1-4, to=1-3] \arrow[from=1-5, to=1-4] \arrow[from=1-6, to=1-5] \arrow[from=1-7, to=1-6] \arrow[from=1-8, to=1-7]
\end{tikzcd}
\end{equation*}
Then, for system III, we have
\begin{equation*}
\begin{tikzcd}[row sep=-1.1ex, column sep=1.4ex]
Jz_{7} & Jz_{6} & Jz_{5} & Jz_{4} & Jz_{3} & Jz_{2} & Jz_{1} & Jz_{0},\\
\indown & \indown & \indown & \indown & \indown & \indown & \indown & \indown\\
\Sigma^{\mathrm{in}}_{i_{3}\sigma_{3}}
& \Sigma^{\mathrm{out}}_{i_{3}\sigma_{3}}
& \Sigma^{\mathrm{in}}_{i_{2}\sigma_{2}}
& \Sigma^{\mathrm{out}}_{i_{2}\sigma_{2}}  
& \Sigma^{\mathrm{in}}_{i_{1}\sigma_{1}}
& \Sigma^{\mathrm{out}}_{i_{1}\sigma_{1}}
& \Sigma^{\mathrm{in}}_{1\sigma_{0}}
& \Sigma^{\mathrm{out}}_{1\sigma_{0}}
\arrow[from=1-2, to=1-1] \arrow[from=1-3, to=1-2] \arrow[from=1-4, to=1-3] \arrow[from=1-5, to=1-4] \arrow[from=1-6, to=1-5] \arrow[from=1-7, to=1-6] \arrow[from=1-8, to=1-7]
\end{tikzcd}
\end{equation*}
where $J = \left(\begin{smallmatrix} 0 & 1\\ 1 & 0 \end{smallmatrix}\right)$. Note that $Jz_{1} = J (T_{1}^{\mathrm{glo}})^{-1}(z_{0})$. For system II, the point $z_{0}$ lies in $\Pi^{\mathrm{in}}_{1\sigma_{0}} \cap R^{s}$, however, $\Pi^{\mathrm{in}}_{1\sigma_{0}} \cap R^{s}$ is a narrow cone containing the horizontal axis on $\Pi^{\mathrm{in}}_{1\sigma_{0}}$ (see Figure \ref{fig:regions_on_cross_sections}). This implies that $Jz_{1}$ for system III lies in $\mathcal{Y}$ (see also \eqref{eq:876yvbfprs}). Since $\mathcal{Y}$ is invariant under the Poincar\'e map (see Lemma \ref{Lem_57fugip9}) and $Jz_{1} \in \mathcal{Y}$, we have that $Jz_{3}$, $Jz_{5}$, and $Jz_{7}$ lie in $\mathcal{Y}$ as well. In particular, Lemma \ref{lem:75rf2} yields
\begin{equation}\label{eq:64gyeuj}
i_{1} = 2, \quad i_{2} = 1, \quad i_{3} = 2.
\end{equation}
In summary,
\begin{equation*}
\begin{tikzcd}[row sep=-1.1ex, column sep=1.4ex]
Jz_{7} & Jz_{6} & Jz_{5} & Jz_{4} & Jz_{3} & Jz_{2} & Jz_{1} & Jz_{0}.\\
\indown & \indown & \indown & \indown & \indown & \indown & \indown & \indown\\
\Sigma^{\mathrm{in}}_{2\sigma_{3}} \cap\mathcal{Y}
& \Sigma^{\mathrm{out}}_{2\sigma_{3}}
& \Sigma^{\mathrm{in}}_{1\sigma_{2}} \cap\mathcal{Y}
& \Sigma^{\mathrm{out}}_{1\sigma_{2}}  
& \Sigma^{\mathrm{in}}_{2\sigma_{1}} \cap\mathcal{Y}
& \Sigma^{\mathrm{out}}_{2\sigma_{1}}
& \Sigma^{\mathrm{in}}_{1\sigma_{0}} \cap\mathcal{Y}
& \Sigma^{\mathrm{out}}_{1\sigma_{0}}
\arrow[from=1-2, to=1-1] \arrow[from=1-3, to=1-2] \arrow[from=1-4, to=1-3] \arrow[from=1-5, to=1-4] \arrow[from=1-6, to=1-5] \arrow[from=1-7, to=1-6] \arrow[from=1-8, to=1-7]
\end{tikzcd}
\end{equation*}

Relation \eqref{eq:64gyeuj} gives the following proposition.
\begin{myprop}\label{prp:switching_backward_dynamics}
The backward dynamics in $R^s$ switches between type-1 and type-2 cross-sections. More precisely, if the backward iteration of $(x_{1}, y_{1}) \in R^{s} \cap \Pi^{\mathrm{in}}_{i\sigma_{1}}$ under the Poincar\'e map is $(x_{0}, y_{0}) \in R^{s} \cap \Pi^{\mathrm{in}}_{j\sigma_{0}}$ for some $i,j\in\{1,2\}$, then $i\neq j$.
\end{myprop}

The evolution of the symbols of the points in $\mathcal{Y}$ is described by the function $S_{u; B^{\prime}, C^{\prime}}$ (see \eqref{eq_234ijvtmbwjr} and Definition \ref{dfn_evolution_function_Su}), where $B^{\prime} = \sgn(\frac{-b_{2}}{\Delta_2})$ and $C^{\prime} = \sgn(\frac{\kappa d_{1} a_{2}}{\Delta_{1} \Delta_{2}})$ for $\Delta_{i} := a_{i} d_{i} - b_{i} c_{i}$ (see \eqref{eq:79ugvjfykuci}, and note that system I and III have the same first integral, and so the same value of $\kappa$). Indeed, if the symbol of $Jz_1$ in system III is ${1}_{\sigma_0}^{A}$, then 
\begin{equation}\label{eq:78rufjhv}
\begin{gathered}
\mathrm{symbol}(Jz_1) = {1}_{\sigma_0}^{A} \longmapsto
\mathrm{symbol}(Jz_3) = {2}_{\sigma_1}^{\sigma_{0} B^{\prime}}\\ \longmapsto
\mathrm{symbol}(Jz_5) = {1}_{\sigma_2}^{A B^{\prime} C^{\prime}} \longmapsto
\mathrm{symbol}(Jz_7) = {2}_{\sigma_3}^{\sigma_{0}},
\end{gathered}
\end{equation}
where $\sigma_{1} = A$, $\sigma_{2} = \sigma_{0} B^{\prime}$ and $\sigma_{3} = AB^{\prime} C^{\prime}$.

\begin{mylem}\label{lem:orientaion_preserving_global_maps}
The global map $T_{i}^{\mathrm{glo}}$ for $i=1,2$ is orientation-preserving, i.e. $\Delta_{i} > 0$.
\end{mylem}

\begin{proof}
Write the normal form in the form \eqref{eq56guyhw}. Along the global segment $p(t)$ of $\Gamma_{1\sigma}$, we have $u_{2} = v_{2} = 0$. Symmetry \eqref{eq:symmetry_2} fixes the $(u_1,v_1)$-plane and preserves its transverse variational subspace. In the ordered transverse coordinates $(u_2,v_2)$ the variational equation is therefore
\begin{equation}\label{eq:variational_equation}
\dot{\xi} = A(t)\xi,
\qquad
A(t)= \begin{pmatrix}
-\lambda_{2} & 0 \\ 0 & \lambda_{2}
\end{pmatrix}
+ \frac{\partial (F_{2}, G_{2})}{\partial (u_{2}, v_{2})}\Big|_{p(t)}.
\end{equation}
If $t_*$ is the travel time for $(0,0)$ on $\Pi^{\mathrm{out}}_{1\sigma}$ to reach $(0,0)$ on $\Pi^{\mathrm{in}}_{1\sigma}$, then Liouville's formula gives
\begin{equation*}
\det\Phi(t_*)=
\exp\!\left(\int_0^{t_*}\operatorname{tr} A(t)\,dt\right)>0,
\end{equation*}
where $\Phi(t)$ is the fundamental matrix of \eqref{eq:variational_equation}. Regarding the determinant of the global map $T_{1}^{\mathrm{glo}}$, note that differentiating the variable traveling time (the time needed to arrive at $\Pi^{\mathrm{in}}_{1\sigma}$) adds a multiple of the vector field at $(0,0)$ on $\Pi^{\mathrm{in}}_{1\sigma}$. This term has zero transverse projection since the endpoint belongs to the invariant $(u_{1},v_{1})$-plane. Therefore this traveling time term contributes nothing to the transverse derivative at $\Gamma_{1\sigma}$. This proves $\Delta_{1} >0$. The proof of $\Delta_{2} >0$ is similar.
\end{proof}

As a consequence of Lemma \ref{lem:orientaion_preserving_global_maps}, we have $B^{\prime} = -\sgn(b_2)$ and $C^{\prime} = \sgn(\kappa d_{1} a_{2})$.

For any even $k$, let $z_{k} \in \Pi^{\mathrm{in}}_{i\sigma} \cap R^{s}$  for $i=1,2$ be a point of system I. Then, for the point $Jz_{k+1}\in \Sigma_{i\sigma} \cap \mathcal{Y}$ of system III, we have $Jz_{k+1} = J (T^{\mathrm{glo}}_{i})^{-1}(z_{k})$, where $T^{\mathrm{glo}}_{i}$ is the global map along $\Gamma_{i\pm}$ of system I. Take two sample points of system I: $z_{k} = \left(\begin{smallmatrix} q \\ 0 \end{smallmatrix}\right) \in \Pi^{\mathrm{in}}_{1\sigma} \cap R^{s}$ and $z_{k} = \left(\begin{smallmatrix} 0 \\ q \end{smallmatrix}\right) \in \Pi^{\mathrm{in}}_{2\sigma} \cap R^{s}$, where $\lvert q\rvert$ is small. Setting $\Delta_{i} := a_{i} d_{i} - b_{i} c_{i}$, by Taylor's theorem, we have
\begin{equation}\label{eq:876yvbfprs}
\begin{aligned}
Jz_{k+1} &= \frac{1}{\Delta_{1}} \left(\begin{matrix}
0 & 1\\ 1& 0
\end{matrix}\right)\left(\begin{matrix}
d_{1} & -b_{1} \\ -c_{1} & a_{1}
\end{matrix}\right) \left(\begin{matrix} q \\ 0 \end{matrix}\right) + \left(\begin{matrix} o(\lvert q\rvert) \\ o(\lvert q\rvert) \end{matrix}\right)\\
&= \frac{q}{\Delta_{1}} \left(\begin{matrix}
-c_{1} \\ d_{1}
\end{matrix}\right) + \left(\begin{matrix} o(\lvert q\rvert) \\ o(\lvert q\rvert) \end{matrix}\right)
\end{aligned}
\end{equation}
and for the second sample point, we have
\begin{equation*}
\begin{aligned}
Jz_{k+1} &= \frac{1}{\Delta_{2}}
\begin{pmatrix} 0 & 1\\ 1& 0 \end{pmatrix}
\begin{pmatrix} d_{2} & -b_{2} \\ -c_{2} & a_{2} \end{pmatrix} \begin{pmatrix} 0 \\ q \end{pmatrix} + \begin{pmatrix} o(\lvert q\rvert) \\ o(\lvert q\rvert) \end{pmatrix}\\
&= \frac{q}{\Delta_{2}} \begin{pmatrix}
a_{2} \\ -b_{2}
\end{pmatrix} + \begin{pmatrix} o(\lvert q\rvert) \\ o(\lvert q\rvert) \end{pmatrix}.
\end{aligned}
\end{equation*}
Therefore, for even $k$ and $q = \pm$, we have
\begin{itemize}
\setlength{\itemsep}{-2pt}
\item $z_{k}$ in system I is assigned to ${1}_{\sigma}^{q}$ if and only if $Jz_{k+1}$ in system III is assigned to ${1}_{\sigma}^{\sgn(q d_{1})}$.
\item $z_{k}$ in system I is assigned to ${2}_{\sigma}^{q}$ if and only if $Jz_{k+1}$ in system III is assigned to ${2}_{\sigma}^{q B^{\prime}}$.
\end{itemize}
Let $A^{\prime} = \sgn(\frac{d_{1}}{\Delta_{1}}) = \sgn(d_1)$. Setting $A = q A^{\prime}$ in \eqref{eq:78rufjhv}, the symbols in system I evolve as
\begin{equation}\label{eq:0oh8ygi34w}
\begin{gathered}
\mathrm{symbol}(z_{0}) = {1}_{\sigma_0}^{q} \longmapsto
\mathrm{symbol}(z_{2}) = {2}_{\sigma_1}^{\sigma_{0}} \longmapsto\\
\mathrm{symbol}(z_{4}) = {1}_{\sigma_2}^{q B^{\prime} C^{\prime}} \longmapsto
\mathrm{symbol}(z_{6}) = {2}_{\sigma_3}^{\sigma_{0} B^{\prime}} \longmapsto  \, \ldots,
\end{gathered}
\end{equation}
where $\sigma_{1} = qA^{\prime}$, $\sigma_{2} = \sigma_{0} B^{\prime}$ and $\sigma_{3} = q A ^{\prime} B^{\prime} C^{\prime}$.

\begin{mydefn}
Let $E := \{{1}_{+}^{+}, {1}_{+}^{-}, {1}_{-}^{+}, {1}_{-}^{-}, {2}_{+}^{+}, {2}_{+}^{-}, {2}_{-}^{+}, {2}_{-}^{-}\}$ be the set of all symbols. Given
\begin{equation*}
A^{\prime} = \sgn(d_{1}),\quad
B^{\prime} = \sgn(-b_{2}),\quad
C^{\prime} = \sgn(\kappa d_{1} a_{2}),
\end{equation*}
we define the bijection $S_{s} = S_{s; A^{\prime},B^{\prime},C^{\prime}}: E \rightarrow E$ as in Table \ref{table:bijection_Ss}.
\end{mydefn}

\begin{table}[!htb]
\centering
\begingroup
\normalfont\normalsize
\setlength{\tabcolsep}{6pt}
\renewcommand{\arraystretch}{1.12}
\begin{tabular}{@{}ccc@{\hspace{1.4em}}cccc@{\hspace{1.4em}}cccc@{}}
\toprule
\multicolumn{3}{c}{Signs}
& \multicolumn{4}{c}{Input $x={1}_{\sigma}^{q}$}
& \multicolumn{4}{c}{Input $x={2}_{\sigma}^{q}$} \\
\cmidrule(lr){1-3}\cmidrule(lr){4-7}\cmidrule(lr){8-11}
$A^{\prime}$ & $B^{\prime}$ & $C^{\prime}$
& ${1}_{+}^{+}$ & ${1}_{+}^{-}$ & ${1}_{-}^{+}$ & ${1}_{-}^{-}$
& ${2}_{+}^{+}$ & ${2}_{+}^{-}$ & ${2}_{-}^{+}$ & ${2}_{-}^{-}$ \\
\midrule
$+$ & $+$ & $+$
& ${2}_{+}^{+}$ & ${2}_{-}^{+}$ & ${2}_{+}^{-}$ & ${2}_{-}^{-}$
& ${1}_{+}^{+}$ & ${1}_{-}^{+}$ & ${1}_{+}^{-}$ & ${1}_{-}^{-}$ \\
$+$ & $+$ & $-$
& ${2}_{-}^{+}$ & ${2}_{+}^{+}$ & ${2}_{-}^{-}$ & ${2}_{+}^{-}$
& ${1}_{+}^{+}$ & ${1}_{-}^{+}$ & ${1}_{+}^{-}$ & ${1}_{-}^{-}$ \\
\addlinespace[2pt]
$+$ & $-$ & $+$
& ${2}_{-}^{-}$ & ${2}_{+}^{-}$ & ${2}_{-}^{+}$ & ${2}_{+}^{+}$
& ${1}_{+}^{+}$ & ${1}_{-}^{+}$ & ${1}_{+}^{-}$ & ${1}_{-}^{-}$ \\
$+$ & $-$ & $-$
& ${2}_{+}^{-}$ & ${2}_{-}^{-}$ & ${2}_{+}^{+}$ & ${2}_{-}^{+}$
& ${1}_{+}^{+}$ & ${1}_{-}^{+}$ & ${1}_{+}^{-}$ & ${1}_{-}^{-}$ \\
\midrule
$-$ & $+$ & $+$
& ${2}_{-}^{+}$ & ${2}_{+}^{+}$ & ${2}_{-}^{-}$ & ${2}_{+}^{-}$
& ${1}_{+}^{-}$ & ${1}_{-}^{-}$ & ${1}_{+}^{+}$ & ${1}_{-}^{+}$ \\
$-$ & $+$ & $-$
& ${2}_{+}^{+}$ & ${2}_{-}^{+}$ & ${2}_{+}^{-}$ & ${2}_{-}^{-}$
& ${1}_{+}^{-}$ & ${1}_{-}^{-}$ & ${1}_{+}^{+}$ & ${1}_{-}^{+}$ \\
\addlinespace[2pt]
$-$ & $-$ & $+$
& ${2}_{+}^{-}$ & ${2}_{-}^{-}$ & ${2}_{+}^{+}$ & ${2}_{-}^{+}$
& ${1}_{+}^{-}$ & ${1}_{-}^{-}$ & ${1}_{+}^{+}$ & ${1}_{-}^{+}$ \\
$-$ & $-$ & $-$
& ${2}_{-}^{-}$ & ${2}_{+}^{-}$ & ${2}_{-}^{+}$ & ${2}_{+}^{+}$
& ${1}_{+}^{-}$ & ${1}_{-}^{-}$ & ${1}_{+}^{+}$ & ${1}_{-}^{+}$ \\
\bottomrule
\end{tabular}
\endgroup
\caption{The symbol-transition bijection
$S_{s; A^{\prime},B^{\prime},C^{\prime}}:E\to E$. For each sign triple $(A^{\prime},B^{\prime},C^{\prime})$, the entry in column $x$ is $S_{s; A^{\prime},B^{\prime},C^{\prime}}(x)$, the symbol at the next incoming-section intersection in forward time, whenever that intersection is defined.}
\label{table:bijection_Ss}
\end{table}

The bijection $S_s$ describes how symbols evolve as we move in $R^s$ in forward time. Accordingly, evolution \eqref{eq:0oh8ygi34w} is formulated as
\begin{equation*}
{1}_{\sigma_0}^{q} \,\, S_s^{-1}({1}_{\sigma_0}^{q})\,\, S_s^{-2}({1}_{\sigma_0}^{q})\,\, S_s^{-3}({1}_{\sigma_0}^{q})\,\, S_s^{-4}({1}_{\sigma_0}^{q}) \, \ldots.
\end{equation*}
We further have that $S_s^{4} = \mathrm{id}_{E}$.

\subsection{Proofs of Theorems \ref{thm:dynamic_near_loops}, \ref{thm:coding_unstable} and \ref{thm:coding_stable}}

We prove Theorems \ref{thm:dynamic_near_loops}, \ref{thm:coding_unstable} and \ref{thm:coding_stable} in this section. The main technical step, the construction of the invariant manifolds in Lemmas \ref{Lem_909tfud34} and \ref{Lem_0ouihiu}, is proved in Section \ref{sec:existence_invariant_manifolds}.

The next lemma describes the set of the points $(x_{0}, y_{0})\in Y$ for which all the backward iterations exist and lie in $Y$, i.e. $(x_{i}, y_{i}) \in Y$ for all $i\leq 0$.

\begin{mylem}\label{Lem_909tfud34}
Let $B_{\epsilon}$ be the union of the $\epsilon$-balls centered at $(0,0)$ on all the four cross-sections $\Pi^{\mathrm{in}}_{i\sigma}$. Then, for a sufficiently small $\epsilon > 0$, and every $i=1,2$ and $\sigma = \pm$, there exists a unique $C^{1}$ curve $G^{u}_{i\sigma}$ in $B_{\epsilon} \cap \Pi^{\mathrm{in}}_{i\sigma}$ passing through $(0,0)$ such that the following hold.
\begin{enumerate}[(i)]
\item\label{item:unstable_curve_1} Each curve lies entirely in $Y$. Moreover, the slope of the tangent line to any of the two curves $G^{u}_{1-}$ and $G^{u}_{1+}$ (resp. $G^{u}_{2-}$ and $G^{u}_{2+}$) at $(0,0)$ is $\frac{d_{1}}{b_{1}}$ (resp. $\frac{c_{2}}{a_{2}}$).

\item\label{item:unstable_curve_2} Let $G^{u}$ be the union of the four curves, i.e. $G^{u} := G^{u}_{1+} \cup G^{u}_{1-} \cup G^{u}_{2+} \cup G^{u}_{2-}$. Then, $G^{u} \setminus \{(0,0)\}$ is invariant under the Poincar\'e map. To be more precise, note that by the previous item, each curve $G^{u}_{i\sigma}$ is the union of the point $(0,0)$, and two pieces $G^{u}_{{i}_{\sigma}^{+}}$ and $G^{u}_{{i}_{\sigma}^{-}}$ lying in the upper and lower half of $\Pi^{\mathrm{in}}_{i\sigma}$, respectively. Then,
\begin{equation*}
T\left(G^{u}_{{i}_{\sigma}^{p}} \setminus \{(0,0)\}\right)\Big\vert_{B_{\epsilon}} = G^{u}_{S_{u}({i}_{\sigma}^{p})} \setminus \{(0,0)\},
\end{equation*}
where $T$ is the Poincar\'e map (see Definition \ref{dfn:Poincare_map_T}), $p = \pm$, and $S_u$ is the evolution rule introduced in Definition \ref{dfn_evolution_function_Su}.

\item\label{item:unstable_curve_3} All the backward iterations of a point $(x_{0}, y_{0})$ in $Y_{\epsilon} := Y \cap B_{\epsilon}$ exist and lie in $Y_{\epsilon}$, i.e. $(x_{i}, y_{i}) \in Y_{\epsilon}$ for all $i\leq 0$, if and only if $(x_{0}, y_{0}) \in G^{u} \setminus \{(0,0)\}$.
\end{enumerate}
\end{mylem}

We prove this lemma by applying the Afraimovich--Shilnikov theorem on the Poincar\'e maps. Section \ref{sec:existence_invariant_manifolds} is dedicated to the proof of this lemma. By reversing time and applying the exact same analysis as in the proof of Lemma \ref{Lem_989giwykd}, we obtain the following lemma.

\begin{mylem}\label{Lem_0ouihiu}
Recall $B_{\epsilon}$ from Lemma \ref{Lem_909tfud34}. Then, for a sufficiently small $\epsilon > 0$, and every $i=1,2$ and $\sigma = \pm$, there exists a unique $C^{1}$ curve $G^{s}_{i\sigma}$ in $B_{\epsilon} \cap \Pi^{\mathrm{in}}_{i\sigma}$ passing through $(0,0)$ such that the following hold.
\begin{enumerate}[(i)]
\item Each curve lies entirely in $R^s \cup \{(0,0)\}$. Moreover, the curves $G^{s}_{1-}$ and $G^{s}_{1+}$ (resp. $G^{s}_{2-}$ and $G^{s}_{2+}$) are tangent to the horizontal (resp. vertical) axis at $(0,0)$.

\item Let $G^{s}$ be the union of the four curves, i.e. $G^{s} := G^{s}_{1+} \cup G^{s}_{1-} \cup G^{s}_{2+} \cup G^{s}_{2-}$. Then, $G^{s} \setminus \{(0,0)\}$ is invariant under the inverse of the Poincar\'e map. To be more precise, note that by the previous item, each curve $G^{s}_{i\sigma}$ is the union of the point $(0,0)$, and two pieces $G^{s}_{{i}_{\sigma}^{+}}$ and $G^{s}_{{i}_{\sigma}^{-}}$ lying in the right (resp. upper) and left (resp. lower) half of $\Pi^{\mathrm{in}}_{i\sigma}$, respectively, when $i= 1$ (resp. $i= 2$). Then,
\begin{equation*}
T^{-1}\left(G^{s}_{{i}_{\sigma}^{p}} \setminus \{(0,0)\}\right)\Big\vert_{B_{\epsilon}} = G^{s}_{S_s^{-1}({i}_{\sigma}^p)} \setminus \{(0,0)\},
\end{equation*}
where $T$ is the Poincar\'e map (see Definition \ref{dfn:Poincare_map_T}), $p = \pm$ , and $S_s$ is the evolution rule given by Table \ref{table:bijection_Ss}.

\item All the forward iterations of a point $(x_{0}, y_{0})$ in $R^{s}_{\epsilon} := R^{s} \cap B_{\epsilon}$ exist and lie in $R^{s}_{\epsilon}$, i.e. $(x_{i}, y_{i}) \in R^{s}_{\epsilon}$ for all $i\geq 0$, if and only if $(x_{0}, y_{0}) \in G^{s} \setminus \{(0,0)\}$.
\end{enumerate}
\end{mylem}

\begin{myrem}
Following Lemma \ref{Lem_66fiytfi5} (resp. Lemma \ref{Lem_989giwykd}), the points on the curves $G_{i\sigma}^{u}$ (resp. $G_{i\sigma}^{s}$) converge to $(0,0)$ by the backward (resp. forward) iteration of $T$. So, the curves $G_{i\sigma}^{u}$ (resp. $G_{i\sigma}^{s}$) can be considered as the unstable (resp. stable) invariant manifolds of the points $(0,0)$.
\end{myrem}

\begin{myrem}
Following symmetries \eqref{eq:symmetry_1} and \eqref{eq:symmetry_2}, the curves $G^{u}_{1-}$ and $G^{u}_{1+}$ have to be identical. The same is true for $G^{u}_{2-}$ and $G^{u}_{2+}$.
\end{myrem}

We are now ready to prove Theorems \ref{thm:dynamic_near_loops}, \ref{thm:coding_unstable} and \ref{thm:coding_stable}.

\begin{proof}[Proof of Theorem \ref{thm:dynamic_near_loops}]
Let $\Lambda^{u}$ be the set of all points in $D$ that stay in $D$ for all the backward iterations of $T$. Take $z_{0} \in \Lambda^{u}$, and let  $z_{-k} = T^{-k}(z_{0})$ for all $k \geq 0$. By definition, $z_{-k} \in \Lambda^{u}$ for all $k\geq 0$. If $z_{0} \in R^{s}$, then according to Lemma \ref{Lem_989giwykd}, the entire sequence $\{z_{-k}:\, k\geq 0\}$ lies in $R^{s}$ and exponentially expands as $k\rightarrow \infty$. Thus, it eventually leaves $B_{\epsilon}$. This implies that $\{z_{-k}:\, k\geq 0\}$ has no intersection with $R^{s}$. Therefore, $\{z_{-k}:\, k\geq 0\} \subset R^{u}$.

Now, suppose $\{z_{-k}:\, k\geq 0\} \subset R^{u}$. For any $k\geq 0$, since $z_{-k-1} \in D^{u}$, by Lemma \ref{Lem_57fugip9}, we have $T(z_{-k-1}) = z_{-k} \in Y$. This implies that $\Lambda^{u}$ is indeed the set of the points in $Y$ that stay in $Y$ for all the backward iterations of $T$. However, by Lemma \ref{Lem_909tfud34}, we have $\Lambda^{u} = G^{u}\setminus \{(0,0)\}$. Applying the flow of the system locally to $G^{u}$, including its section origins, gives $W^{u}_{\mathrm{loc}}(B)$. The stable manifold is obtained analogously by applying the flow locally to $G^{s}$. This completes the proof of Theorem \ref{thm:dynamic_near_loops}.
\end{proof}

\begin{proof}[Proof of Theorems \ref{thm:coding_unstable} and \ref{thm:coding_stable}]
The unstable manifold of the homoclinic orbits intersects with each cross-section $\Pi^{\mathrm{in}}_{i\sigma}$ along $G_{i\sigma}^u$ that lies in $Y$ (see Lemma \ref{Lem_909tfud34}). As discussed above, the evolution of symbols on $Y$ is described by the bijection $S_u$ given by Table \ref{table:bijection_Su}. Theorem \ref{thm:coding_unstable} is obtained from iterating any arbitrary symbol under the map $S_u$.

The stable manifold of the homoclinic orbits intersects with each cross-section $\Pi^{\mathrm{in}}_{i\sigma}$ along $G_{i\sigma}^s$ that lies in $R^s \cup \{(0,0)\}$ (see Lemma \ref{Lem_0ouihiu}). As discussed above, the evolution of symbols on $R^s$ under the forward iterations of the Poincar\'e map is described by the bijection $S_s$ given by Table \ref{table:bijection_Ss}. To obtain Theorem \ref{thm:coding_stable}, it is sufficient to iterate any arbitrary symbol by the map $S_s$. This completes the proof of both Theorems \ref{thm:coding_unstable} and \ref{thm:coding_stable}.
\end{proof}

\section{Existence of the invariant manifolds}\label{sec:existence_invariant_manifolds}

We establish the existence of the unstable invariant manifolds of the homoclinic bouquet and prove Lemma \ref{Lem_909tfud34} in this section. The case of the stable manifold (see Lemma~\ref{Lem_0ouihiu}) follows by reversing time and applying the same analysis as in the proof of Lemma \ref{Lem_989giwykd}.

Starting from a point in $Y$ on a cross-section $\Pi^{\mathrm{in}}_{1\sigma_{0}}$ for some $\sigma_{0} = \pm$, as we move forward in time, we intersect $\Pi^{\mathrm{in}}_{2\sigma_{1}}$, $\Pi^{\mathrm{in}}_{1\sigma_{2}}$, and $\Pi^{\mathrm{in}}_{2\sigma_{3}}$ for some $\sigma_{1},\sigma_{2},\sigma_{3} = \pm$, and then return to $\Pi^{\mathrm{in}}_{1\sigma_{0}}$ (see \eqref{eq_evolution_sigma}). Moreover, all the intersection points lie in $Y$. This dynamics is described by $T^4$, where $T$ is the Poincar\'e map (see Definition \ref{dfn:Poincare_map_T}). The proof of the existence of the unstable manifold of the homoclinic bouquet reduces to proving the existence of unstable invariant $C^1$ curves passing through $(0,0)$ for $T^{4}\big|_{Y \cap\Pi^{\mathrm{in}}_{1\sigma_{0}}}$. 

First, observe that
\begin{equation*}
T^{4}\big|_{Y \cap\Pi^{\mathrm{in}}_{1\sigma_{0}}} = T^{2}\big|_{Y \cap\Pi^{\mathrm{in}}_{1\sigma_{2}}} \circ T^{2}\big|_{Y \cap\Pi^{\mathrm{in}}_{1\sigma_{0}}}.
\end{equation*}
Since the system is invariant with respect to symmetries \eqref{eq:symmetry_1} and \eqref{eq:symmetry_2}, the maps 
$T^{2}|_{Y \cap\Pi^{\mathrm{in}}_{1\sigma_{0}}}$ and 
$T^{2}|_{Y \cap\Pi^{\mathrm{in}}_{1\sigma_{2}}}$ are equal. Note that we do not necessarily have $\sigma_{0} = \sigma_{2}$ (this holds only when $c_{2} > 0$, see Table \ref{table:bijection_Su} and Theorem \ref{thm:coding_unstable}), and therefore these two maps may have different domains and codomains. Here, by the equality we mean that, as mappings defined on a $(u_{2}, v_{2})$-coordinate plane, they coincide. Therefore, to prove the existence of the unstable curve for $T^{4}$, it is sufficient to prove the existence of the unstable curve passing through $(0,0)$ for $T^{2}|_{Y \cap\Pi^{\mathrm{in}}_{1\sigma_{0}}}$.

Take a point $(u_{20}, v_{20})$ in $Y \cap\Pi^{\mathrm{in}}_{1\sigma_{0}}$ and let $(\overline{u}_{20}, \overline{v}_{20}) \in Y \cap\Pi^{\mathrm{in}}_{1\sigma_{2}}$ be the image of this point under $T^{2}$. Then, there is a $\sigma_{1} = \pm$ such that
\begin{equation}\label{eq_composition_map}
T_{2\sigma_{1} 1\sigma_{2}} \circ T_{1\sigma_{0} 2\sigma_{1}}(u_{20}, v_{20}) = (\overline{u}_{20}, \overline{v}_{20}).
\end{equation}
The idea of the proof of the existence of the unstable manifold in a nutshell is to find appropriate coordinates in which the cross-form of map \eqref{eq_composition_map} has small derivatives and satisfies the assumptions of the Afraimovich--Shilnikov theorem. To this end, we first calculate the differential of map \eqref{eq_composition_map} in $(u_{2}, v_{2})$-coordinates in Section \ref{Sec_differential_composition_map}. We then show that this map admits a cross-map satisfying Afraimovich--Shilnikov theory in Section \ref{sec:Afraimovich_Shilnikov}. The proofs of two technical lemmas in Section \ref{Sec_differential_composition_map} are deferred to Appendices \ref{Sec_Proof_Lem_derivatives_local_map_12} and \ref{Sec_Proof_Lem_derivatives_local_map_21}.

Regarding the domain of the composition map \eqref{eq_composition_map}, note that by \eqref{eq096te6dft}, the size of the image of $(u_{20}, v_{20})$ under $T_{1\sigma_{0} 2\sigma_{1}}$, for a fixed $\delta$, is comparable to $\epsilon^{\gamma}$ which can be larger than $\epsilon$ as $\gamma <1$. Consequently, such point may not fall into the domain of $T_{2\sigma_{1} 1\sigma_{2}}$. Due to this, we choose $0 < \rho_{0} < \epsilon$ sufficiently small for which \eqref{eq_composition_map} is well-defined on
\begin{equation}\label{eq:restricted_domain_of_composition_map}
Y \cap B_{\rho_{0}},
\end{equation}
where $B_{\rho_{0}}$ is the $\rho_{0}$-ball on $\Pi^{\mathrm{in}}_{1\sigma_{0}}$ centered at $(0,0)$. We particularly consider this restricted domain in Section \ref{sec:Afraimovich_Shilnikov}.

\subsection{Estimation of the differential of the composition map}\label{Sec_differential_composition_map}

The aim of this section is to estimate the differential of $T_{2\sigma_{1} 1\sigma_{2}} \circ T_{1\sigma_{0} 2\sigma_{1}}$ for some $\sigma_{0},\sigma_{1}, \sigma_{2}\in \{-,+\}$. It makes the discussion here clearer if, alongside the $(u_{i}, v_{i})$ coordinates on the cross-sections that we have been using so far, we introduce a new notation as follows. Consider a point $(x_{0}, y_{0}) \in Y\cap \Pi^{\mathrm{in}}_{1\sigma_{0}}$ whose forward orbit intersects $\Pi^{\mathrm{in}}_{2\sigma_{1}}$ at some point, call it $(x_{1}, y_{1})$. Thus, $(x_{1}, y_{1}) = T_{1\sigma_{0} 2\sigma_{1}}(x_{0},y_{0})$, and $(x_{1}, y_{1})\in Y$. Let the forward orbit of $(x_{1}, y_{1})$ intersect $\Pi^{\mathrm{in}}_{1\sigma_{2}}$ at a point $(x_{2}, y_{2})$ for some $\sigma_{2} = \pm$. We then have $T_{2\sigma_{1} 1\sigma_{2}}(x_{1}, y_{1}) = (x_{2}, y_{2}) \in Y$, and so
\begin{equation*}
(x_{2}, y_{2}) = T_{2\sigma_{1} 1\sigma_{2}} \circ T_{1\sigma_{0} 2\sigma_{1}}(x_{0}, y_{0}).
\end{equation*}
Looking back at the $(u_{i}, v_{i})$ coordinates, we have $x_{0} = u_{20}$ and $y_{0} = v_{20}$, where $(u_{20}, v_{20}) \in \Pi^{\mathrm{in}}_{1\sigma_{0}}$ represents the initial point, $x_{1} = u_{10}$ and $y_{1} = v_{10}$, where $(u_{10}, v_{10}) \in \Pi^{\mathrm{in}}_{2\sigma_{1}}$ represents the first iteration, and $x_{2} = u_{20}$ and $y_{2} = v_{20}$, where $(u_{20}, v_{20}) \in \Pi^{\mathrm{in}}_{1\sigma_{2}}$ represents the second iteration.

Regarding the differential of $T_{2\sigma_{1} 1\sigma_{2}} \circ T_{1\sigma_{0} 2\sigma_{1}}$, using chain rule, we obtain
\begin{equation*}
\frac{\partial (x_{2}, y_{2})}{\partial (x_{0}, y_{0})} = DT_{2\sigma_{1} 1\sigma_{2}}(x_{1}, y_{1}) \cdot DT_{1\sigma_{0} 2\sigma_{1}}(x_{0}, y_{0}).
\end{equation*}
Thus, all we need to do is to estimate the differentials of $T_{1\sigma_{0} 2\sigma_{1}}$ and $T_{2\sigma_{1} 1\sigma_{2}}$.

\paragraph*{Estimation of $DT_{1\sigma_{0} 2\sigma_{1}}$.} Recall that $T_{1\sigma_{0} 2\sigma_{1}} = T^{\mathrm{glo}}_{2} \circ T^{\mathrm{loc}}_{1\sigma_{0} 2\sigma_{1}}$, and so $DT_{1\sigma_{0} 2\sigma_{1}} = DT^{\mathrm{glo}}_{2} \cdot DT^{\mathrm{loc}}_{1\sigma_{0} 2\sigma_{1}}$. Following \eqref{eq_67ihbu3y4g4byu}, when $(u_{20},v_{20})\in Y$, we have $v_{1\tau} = o(\lvert u_{1\tau}\rvert)$. Then, \eqref{eq_Taylor_expansion_global_map} gives
\begin{equation*}
T^{\mathrm{glo}}_{2}(u_{1\tau}, v_{1\tau}) = \left(a_{2} u_{1\tau} + b_{2} v_{1\tau} + O\left(\lvert u_{1\tau}\rvert^{3}\right),\, c_{2} u_{1\tau} + d_{2} v_{1\tau} + O\left(\lvert u_{1\tau}\rvert^{3}\right)\right)
\end{equation*}
and
\begin{equation*}
DT^{\mathrm{glo}}_{2}(u_{1\tau}, v_{1\tau}) = \left(\begin{matrix}
a_{2} + O(u_{1\tau}^{2}) & b_{2} + O(u_{1\tau}^{2})\\
c_{2} + O(u_{1\tau}^{2}) & d_{2} + O(u_{1\tau}^{2})
\end{matrix}\right).
\end{equation*}
Regarding the local map, we have the following lemma (the proof is postponed to Appendix \ref{Sec_Proof_Lem_derivatives_local_map_12}).

\begin{mylem}\label{Lem_derivatives_local_map_12}
Write $T^{\mathrm{loc}}_{1\sigma_{0} 2\sigma_{1}}(u_{20}, v_{20}) = \left(p_{12}(u_{20}, v_{20}), q_{12}(u_{20}, v_{20})\right)$. For $(u_{20}, v_{20})\in Y$, we have
\begin{equation}\label{eq_derivatives_local_map_12}
\begin{aligned}
\frac{\partial p_{12}}{\partial u_{20}} &= O(\lvert v_{20}\rvert^{\gamma}\delta^{2-\gamma}),\\
\frac{\partial p_{12}}{\partial v_{20}} &= \gamma\sigma_{0}\sigma_{1}\delta^{1-\gamma} \lvert v_{20}\rvert^{\gamma -1} [1 + O(\delta^{2})],\\
\frac{\partial q_{12}}{\partial u_{20}} &= \kappa\gamma^{-1} \sigma_{0}\sigma_{1} \lvert v_{20}\rvert^{1-\gamma} \delta^{\gamma -1} [1 + O(\delta^{2})],\\
\frac{\partial q_{12}}{\partial v_{20}} &= (\gamma^{-1} -1)\kappa \sigma_{0} u_{20} \lvert v_{20}\rvert^{-\gamma} \delta^{\gamma - 1}[1 + O(\delta^{2})].
\end{aligned}
\end{equation}
\end{mylem}

In the new notation, we write $(x_{1}, y_{1}) = T_{1\sigma_{0} 2\sigma_{1}}(x_{0}, y_{0})$, where $x_{0} = u_{20}$ and $y_{0} = v_{20}$. Define
\begin{equation}\label{eq_7db8i2g}
\zeta_{1} = \lvert y_{0}\rvert^{2\gamma - 1} \delta^{3-2\gamma},\qquad
\zeta_{2} = \lvert y_{0}\rvert^{2-2\gamma} \delta^{2\gamma -2},\qquad
\zeta_{3} = \lvert y_{0}\rvert^{4\gamma -2} \delta^{4-4\gamma}.
\end{equation}
By virtue of \eqref{eq_derivatives_local_map_12}, for $(x_{0}, y_{0}) \in Y$, we have $\frac{\partial p_{12}}{\partial x_{0}} = O(\zeta_{1} \cdot \frac{\partial q_{12}}{\partial x_{0}})$ and $\frac{\partial q_{12}}{\partial y_{0}} = O(\zeta_{2} \cdot \frac{\partial p_{12}}{\partial y_{0}})$. On the other hand, by \eqref{eq_45tu9yo}, we have $u_{1\tau}^{2} = O(\lvert y_{0}\rvert^{2\gamma} \delta^{2- 2\gamma})$. We then obtain
\begin{equation}\label{eq21rtyit76fqugjw}
\begin{aligned}
D T_{1\sigma_{0} 2\sigma_{1}}(x_{0}, y_{0}) &= \begin{pmatrix}
a_{2} + O(u_{1\tau}^{2}) & b_{2} + O(u_{1\tau}^{2})\\[1ex]
c_{2} + O(u_{1\tau}^{2}) & d_{2} + O(u_{1\tau}^{2})
\end{pmatrix} \begin{pmatrix}
\frac{\partial p_{12}}{\partial x_{0}} & \frac{\partial p_{12}}{\partial y_{0}}\\[1ex]
\frac{\partial q_{12}}{\partial x_{0}} & \frac{\partial q_{12}}{\partial y_{0}}
\end{pmatrix}\\[1ex]
&= \begin{pmatrix}
[b_{2} + O(\zeta_{1})] \frac{\partial q_{12}}{\partial x_{0}} &
[a_{2} + O(\zeta_{2})] \frac{\partial p_{12}}{\partial y_{0}}
\\[1ex]
[d_{2} + O(\zeta_{1})] \frac{\partial q_{12}}{\partial x_{0}} &
[c_{2} + O(\zeta_{2})] \frac{\partial p_{12}}{\partial y_{0}}
\end{pmatrix}.
\end{aligned}
\end{equation}
Although one could substitute estimates \eqref{eq_derivatives_local_map_12} for $\frac{\partial p_{12}}{\partial v_{20}}$ and $\frac{\partial q_{12}}{\partial u_{20}}$ to simplify \eqref{eq21rtyit76fqugjw}, we avoid doing so. A similar situation occurs for $D T_{2\sigma_{1} 1\sigma_{2}}$ in the next section. This is for technical reasons: in Section~\ref{sec:Afraimovich_Shilnikov}, particularly, in the estimation of the differential of $\mathbf{T}$, certain expressions exhibit cancellation of leading terms. Retaining the above forms of $D T_{1\sigma_{0} 2\sigma_{1}}$ and $D T_{2\sigma_{1} 1\sigma_{2}}$ provides better control of the remaining terms.

\paragraph*{Estimation of $DT_{2\sigma_{1} 1\sigma_{2}}$.}
Recall that $T_{2\sigma_{1} 1\sigma_{2}} = T^{\mathrm{glo}}_{1}\circ T^{\mathrm{loc}}_{2\sigma_{1} 1\sigma_{2}}$, and so $DT_{2\sigma_{1} 1\sigma_{2}} = DT^{\mathrm{glo}}_{1} \cdot DT^{\mathrm{loc}}_{2\sigma_{1} 1\sigma_{2}}$. Following \eqref{eq_98uriygvjbmdjer}, when $(u_{10},v_{10})\in Y$, we have $u_{2\tau} = o(\lvert v_{2\tau}\rvert)$. Then, \eqref{eq_Taylor_expansion_global_map} gives
\begin{equation}\label{eq56yhnto8y}
T^{\mathrm{glo}}_{1}(u_{2\tau}, v_{2\tau}) = \left(a_{1} u_{2\tau} + b_{1} v_{2\tau} + O\left(\lvert v_{2\tau}\rvert^{3}\right), c_{1} u_{2\tau} + d_{1} v_{2\tau} + O\left(\lvert v_{2\tau}\rvert^{3}\right)\right)
\end{equation}
and
\begin{equation*}
DT^{\mathrm{glo}}_{1}(u_{2\tau}, v_{2\tau}) = \left(\begin{matrix}
a_{1} + O(v_{2\tau}^{2}) & b_{1} + O(v_{2\tau}^{2})\\[1ex]
c_{1} + O(v_{2\tau}^{2}) & d_{1} + O(v_{2\tau}^{2})
\end{matrix}\right).
\end{equation*}
Regarding the local map, we have the following lemma which is proved in Appendix \ref{Sec_Proof_Lem_derivatives_local_map_21}.

\begin{mylem}\label{Lem_derivatives_local_map_21}
Write $T^{\mathrm{loc}}_{2\sigma_{1} 1\sigma_{2}}(u_{10}, v_{10}) = \left(p_{21}(u_{10}, v_{10}), q_{21}(u_{10}, v_{10})\right)$. For $(u_{10}, v_{10})\in Y$, we have
\begin{equation}\label{eq_derivatives_local_map_21}
\begin{aligned}
\frac{\partial p_{21}}{\partial u_{10}} &= O(\lvert v_{10}\rvert \delta),\\
\frac{\partial p_{21}}{\partial v_{10}} &= \sigma_{1} \sigma_{2} \frac{1}{\gamma}\lvert v_{10}\rvert^{\frac{1}{\gamma} -1} \delta^{1-\frac{1}{\gamma}} [1 + O(\delta^{2})],\\
\frac{\partial q_{21}}{\partial u_{10}} &= \kappa\gamma \sigma_{1} \sigma_{2} \lvert v_{10}\rvert^{1 -\frac{1}{\gamma}} \delta^{\frac{1}{\gamma} -1} [1 + O(\delta^{2})],\\
\frac{\partial q_{21}}{\partial v_{10}} &= (\gamma -1)\kappa \sigma_{1} u_{10} \lvert v_{10}\rvert^{-\frac{1}{\gamma}} \delta^{\frac{1}{\gamma} -1} [1 + O(\delta^{2})].
\end{aligned}
\end{equation}
\end{mylem}

In the new notation, we write $(x_{2}, y_{2}) = T_{2\sigma_{1} 1\sigma_{2}}(x_{1}, y_{1})$, where $x_{1} = u_{10}$ and $y_{1} = v_{10}$. Let $\eta = \lvert y_{1}\rvert \delta^{-1}$. Following this lemma, for $(x_{1}, y_{1}) \in Y$, we have $\frac{\partial p_{21}}{\partial u_{10}} = O({\eta}^{\frac{1}{\gamma}} \delta^{2} \cdot \frac{\partial q_{21}}{\partial u_{10}})$ and $\frac{\partial p_{21}}{\partial v_{10}} = O({\eta}^{\frac{2}{\gamma} - 2}\cdot \frac{\partial q_{21}}{\partial v_{10}})$. Moreover, by \eqref{eq_45tuvkb7taguer}, we have $v_{2\tau}^{2} = O(\eta^{4-\frac{2}{\gamma}} \delta^{2})$. We then have
\begin{equation}\label{eq989uob43}
\begin{aligned}
& D T_{2\sigma_{1} 1\sigma_{2}}(x_{1}, y_{1})\\[0.5ex]
&\quad = \left(\begin{matrix}
a_{1} + O(v_{2\tau}^{2}) & b_{1} + O(v_{2\tau}^{2})\\[1ex]
c_{1} + O(v_{2\tau}^{2}) & d_{1} + O(v_{2\tau}^{2})
\end{matrix}\right) \left(\begin{matrix}
\frac{\partial p_{21}}{\partial x_{1}} & \frac{\partial p_{21}}{\partial y_{1}}\\[1ex]
\frac{\partial q_{21}}{\partial x_{1}} & \frac{\partial q_{21}}{\partial y_{1}}
\end{matrix}\right)\\[1ex]
&\quad = \left(\begin{matrix}
[b_{1} + O(\eta^{\frac{1}{\gamma}} \delta^{2}) + O(\eta^{4-\frac{2}{\gamma}} \delta^{2})] \frac{\partial q_{21}}{\partial x_{1}}
& [b_{1} + O(\eta^{\frac{2}{\gamma} - 2}) + O(\eta^{4-\frac{2}{\gamma}} \delta^{2})] \frac{\partial q_{21}}{\partial y_{1}}
\\[1ex]
[d_{1} + O(\eta^{\frac{1}{\gamma}} \delta^{2}) + O(\eta^{4-\frac{2}{\gamma}} \delta^{2})] \frac{\partial q_{21}}{\partial x_{1}}
& [d_{1} + O(\eta^{\frac{2}{\gamma} - 2}) + O(\eta^{4-\frac{2}{\gamma}} \delta^{2})] \frac{\partial q_{21}}{\partial y_{1}}
\end{matrix}\right).
\end{aligned}
\end{equation}

\paragraph*{Estimation of $\pd{(x_{2}, y_{2})}{(x_{0}, y_{0})}$.} We calculated $\pd{(x_{2}, y_{2})}{(x_{1},y_{1})}$ at a point $(x_{1}, y_{1})$ in \eqref{eq989uob43}. However, to estimate $\pd{(x_{2}, y_{2})}{(x_{0}, y_{0})}$, we need to evaluate $x_{1}$ and $y_{1}$ in \eqref{eq989uob43} at $T_{1\sigma_{0} 2\sigma_{1}}(x_{0}, y_{0})$. To do so, note that \eqref{eq096te6dft} yields
\begin{equation*}
x_{1} = a_{2} \sigma_{0} \lvert y_{0}\rvert^{\gamma} \delta^{1-\gamma} [1 + O(\delta^{2})]
\qquad\mathrm{and}\qquad
y_{1} = c_{2} \sigma_{0} \lvert y_{0}\rvert^{\gamma} \delta^{1-\gamma} [1 + O(\delta^{2})].
\end{equation*}
Substituting this into \eqref{eq989uob43}, we obtain
\begin{equation*}
\pd{(x_{2}, y_{2})}{(x_{1},y_{1})} = \left(\begin{matrix}
[b_{1} + O(\zeta_{1} \zeta_{2}) + O(\zeta_{3})] \frac{\partial q_{21}}{\partial x_{1}}
& [b_{1} + O(\zeta_{2}) + O(\zeta_{3})] \frac{\partial q_{21}}{\partial y_{1}}
\\[1ex]
[d_{1} + O(\zeta_{1} \zeta_{2}) + O(\zeta_{3})] \frac{\partial q_{21}}{\partial x_{1}}
& [d_{1} + O(\zeta_{2}) + O(\zeta_{3})] \frac{\partial q_{21}}{\partial y_{1}}
\end{matrix}\right),
\end{equation*}
where $\zeta_{1}$, $\zeta_{2}$ and $\zeta_{3}$ are as in \eqref{eq_7db8i2g}. Finally, using this estimate and \eqref{eq21rtyit76fqugjw}, and applying the chain rule, we have
\begin{equation}\label{eq:derivative_composition_map}
\begin{aligned}
\frac{\partial x_{2}}{\partial x_{0}} &= \left[b_{1} b_{2} + O(\zeta_{1})\right] \frac{\partial q_{21}}{\partial x_{1}} \cdot \frac{\partial q_{12}}{\partial x_{0}} + \left[b_{1} d_{2} + O(\zeta_{1} + \zeta_{2})\right] \frac{\partial q_{21}}{\partial y_{1}} \cdot \frac{\partial q_{12}}{\partial x_{0}},\\
\frac{\partial x_{2}}{\partial y_{0}} &= \left[b_{1}a_{2} + O(\zeta_{2} + \zeta_{3})\right] \frac{\partial q_{21}}{\partial x_{1}} \cdot  \frac{\partial p_{12}}{\partial y_{0}} + \left[b_{1}c_{2} + O(\zeta_{2} + \zeta_{3})\right] \frac{\partial q_{21}}{\partial y_{1}}\cdot \frac{\partial p_{12}}{\partial y_{0}},\\
\frac{\partial y_{2}}{\partial x_{0}} &= \left[d_{1} b_{2} + O(\zeta_{1})\right] \frac{\partial q_{21}}{\partial x_{1}} \cdot \frac{\partial q_{12}}{\partial x_{0}} + \left[d_{1} d_{2} + O(\zeta_{1} + \zeta_{2})\right] \frac{\partial q_{21}}{\partial y_{1}} \cdot \frac{\partial q_{12}}{\partial x_{0}},\\
\frac{\partial y_{2}}{\partial y_{0}} &= \left[d_{1}a_{2} + O(\zeta_{2} + \zeta_{3})\right] \frac{\partial q_{21}}{\partial x_{1}} \cdot  \frac{\partial p_{12}}{\partial y_{0}} + \left[d_{1}c_{2} + O(\zeta_{2} + \zeta_{3})\right] \frac{\partial q_{21}}{\partial y_{1}}\cdot \frac{\partial p_{12}}{\partial y_{0}}.
\end{aligned}
\end{equation}

\subsection{Cross-map analysis: Afraimovich--Shilnikov theory and the proof of Lemma~\ref{Lem_909tfud34}}\label{sec:Afraimovich_Shilnikov}

In this section, we use the Afraimovich--Shilnikov theory to prove the existence of the unstable manifold passing through $(0,0)$ for the composition map \eqref{eq_composition_map}, thereby proving Lemma~\ref{Lem_909tfud34}. The main reason we use the Afraimovich--Shilnikov theory, rather than standard techniques from hyperbolic dynamics such as the graph transform method, is that map \eqref{eq_composition_map} is singular. In particular, the point $(0,0)$, which should play the role of a fixed point, is not even contained in its domain. Extending the domain to include this point does not resolve the issue by itself either, since the norm of the differential of the map tends to infinity as we converge to $(0,0)$. Moreover, the domain of \eqref{eq_composition_map} consists of multiple connected components rather than a single open neighborhood of $(0,0)$.

These difficulties are common in the study of dynamics near homoclinics to equilibria. A powerful technique to overcome these obstacles, is to find appropriate coordinates in which \eqref{eq_composition_map} admits a cross-form with small differential. Let us first introduce these coordinates here, and then we will elaborate the idea of cross-form further.

The composition map $T_{2\sigma_{1} 1\sigma_{2}} \circ T_{1\sigma_{0} 2\sigma_{1}}$ takes $(x_{0}, y_{0}) \in Y\cap\Pi_{1\sigma_{0}}^{\mathrm{in}}$ to $(x_{2}, y_{2}) \in Y\cap\Pi_{1\sigma_{2}}^{\mathrm{in}}$. We define new variables
\begin{equation}\label{eq:definition_w_z_variables}
w = \frac{y_{0}}{x_{0}},\qquad z = y_{0}^{\alpha},
\end{equation}
where
\begin{equation}\label{eq:formula_alpha}
\alpha = \frac{1}{2N+1} > 0
\end{equation}
for a sufficiently large positive integer $N$ that will be specified later. It is shown that the dynamics in $w$ direction is strongly contracting and it is strongly expanding in $z$ direction.

We can write the composition map in the new coordinates as $(w,z)\mapsto (\overline{w}, \overline{z})$, by
\begin{equation}\label{eq_map_wz_coordinates}
\overline{w} = f(w, z), \qquad \overline{z} = g(w, z),
\end{equation}
where $f$ and $g$ are defined by
\begin{equation}\label{eq:map_wz_xy}
\overline{w} = f(w, z) = \frac{y_{2}}{x_{2}}\\
\qquad\mathrm{and}\qquad
\overline{z} = g(w, z) = y_{2}^{\alpha}.
\end{equation}
Taking the restricted domain \eqref{eq:restricted_domain_of_composition_map} into account, if we choose
\begin{equation}\label{eq:radius_rectanglular_domain}
\rho < \min\{\rho_{0},\, \rho_{0}^{\alpha} (1+m^{2})^{\frac{-\alpha}{2}}\},
\end{equation}
then map \eqref{eq_map_wz_coordinates} is well-defined on
\begin{equation*}
\left([-m, - m^{-1}] \cup [m^{-1},m]\right) \times \left([-\rho,\rho]\setminus \{0\}\right).
\end{equation*}
Recall also that the images of the points under map \eqref{eq_composition_map} accumulate near the straight line with the slope $\frac{d_{1}}{b_{1}}$, i.e. $\overline{w} \approx \frac{d_{1}}{b_{1}}$. We can extend the domain of map \eqref{eq_map_wz_coordinates} to include $z=0$ by defining
\begin{equation}\label{eq:extended_map_in_wz}
\mathbf{T}(w,z) = \left\{\begin{array}{ll}
(f(w,z), g(w,z)) & \mathrm{if}\,\, z\neq 0,\\
(\frac{d_{1}}{b_{1}}, 0) & \mathrm{if}\,\, z = 0.
\end{array}
\right.
\end{equation}

This extension resolves the domain problem as $\mathbf{T}$ is now defined on a neighborhood of the segment $z=0$, however, this segment still remains a non-removable singularity (e.g. it is shown later in \eqref{eq_derivative_zbar_z} that $|\partial g / \partial z|$ tends to infinity as $z\rightarrow 0$). To deal with this problem, we prove in this section that $\mathbf{T}$ admits a cross-form.

We say $\mathbf{T}: (w,z)\mapsto (\overline{w}, \overline{z})$ admits a cross-form if there exist $C^{1}$ functions $F$ and $G$ such that
\begin{equation*}
\overline{w} = F(w, \overline{z}),\qquad
z = G(w, \overline{z}).
\end{equation*}
We call this map the cross-map of $\mathbf{T}$ and denote it by $\mathbf{T}^{\times}$. Although $\mathbf{T}$ extends continuously across $z=0$, it is not $C^{1}$ there. By contrast, the cross-map $\mathbf{T}^{\times}$ admits a $C^{1}$ extension across $\overline{z}=0$. If there exists a subset $\mathcal{R}$ of the domain of $\mathbf{T}$ for which $\mathbf{T}^{\times}(\mathcal{R}) \subset \mathcal{R}$, then, roughly speaking, the Afraimovich--Shilnikov theory implies that $\mathbf{T}$ possesses an invariant manifold passing through $(\frac{d_1}{b_1},0)$, provided that the differential of $\mathbf{T}^{\times}$ is sufficiently small. Moreover, this manifold contains the $\omega$-limit set of any point of $\mathcal{R}$ whose forward orbit lies entirely in $\mathcal{R}$.

Throughout, we first estimate the differential of $\mathbf{T}$. Next, we use these estimates to show that $\mathbf{T}^{\times}$ exists, and finally, we apply the Afraimovich--Shilnikov theory to $\mathbf{T}^{\times}$, thereby proving Lemma~\ref{Lem_909tfud34}.

\paragraph{Differential of the composition map in $(w,z)$-coordinates.}

Recall new variables \eqref{eq:definition_w_z_variables}, i.e. $w = \frac{y_{0}}{x_{0}}$ and $z = y_{0}^{\alpha}$. Expressing old variables in terms of new ones gives
\begin{equation*}
x_{0} = x_{0}(w,z) = z^{\frac{1}{\alpha}} w^{-1},
\qquad
y_{0} = y_{0}(w,z) = z^{\frac{1}{\alpha}}.
\end{equation*}
This yields
\begin{equation}\label{eq_68igtfygbi}
\frac{\partial x_{0}}{\partial w} = - z^{\frac{1}{\alpha}} w^{-2},\qquad
\frac{\partial x_{0}}{\partial z} = \frac{z^{\frac{1}{\alpha} - 1}}{\alpha w},\qquad
\frac{\partial y_{0}}{\partial w} = 0,\qquad
\frac{\partial y_{0}}{\partial z} = \frac{1}{\alpha} z^{\frac{1}{\alpha} - 1}.
\end{equation}
Recall map \eqref{eq:map_wz_xy}. To obtain $\frac{\partial (f,g)}{\partial (w,z)}$, it is enough to use $\frac{\partial (x_{2}, y_{2})}{\partial (x_{0}, y_{0})}$, estimated in \eqref{eq:derivative_composition_map}, together with \eqref{eq_68igtfygbi}. Taking into account that $\frac{\partial y_{0}}{\partial w} = 0$, we have
\begin{align}
\frac{\partial f}{\partial w} &= \frac{1}{x_{2}} \cdot \frac{\partial x_{0}}{\partial w}\cdot\left(\frac{\partial y_{2}}{\partial x_{0}} - \frac{\partial x_{2}}{\partial x_{0}}\cdot \overline{w}\right),\label{eq_657g7033ww}\\
\frac{\partial f}{\partial z} &= \frac{1}{x_{2}}\cdot \frac{\partial x_{0}}{\partial z} \cdot \left(\frac{\partial y_{2}}{\partial x_{0}} - \frac{\partial x_{2}}{\partial x_{0}} \cdot \overline{w}\right) + \frac{1}{x_{2}}\cdot \frac{\partial y_{0}}{\partial z} \cdot \left(\frac{\partial y_{2}}{\partial y_{0}} - \frac{\partial x_{2}}{\partial y_{0}} \cdot \overline{w}\right),\label{eq_657g7033wz}\\
\frac{\partial g}{\partial w} &= \alpha y_{2}^{\alpha - 1} \frac{\partial y_{2}}{\partial x_{0}}\cdot \frac{\partial x_{0}}{\partial w},\label{eq_657g7033zw}\\
\frac{\partial g}{\partial z} &= \alpha y_{2}^{\alpha - 1} \left(\frac{\partial y_{2}}{\partial x_{0}}\cdot \frac{\partial x_{0}}{\partial z} + \frac{\partial y_{2}}{\partial y_{0}}\cdot \frac{\partial y_{0}}{\partial z}\right).\label{eq_67uohjbnh}
\end{align}

Recall $\zeta_{1}$, $\zeta_{2}$ and $\zeta_{3}$ from \eqref{eq_7db8i2g}. Using $z = y_{0}^{\alpha}$, we also have
\begin{equation}\label{eq_zetas}
\zeta_{1} = \lvert z\rvert^{\frac{2\gamma - 1}{\alpha}} \delta^{3-2\gamma},\qquad
\zeta_{2} = \lvert z\rvert^{\frac{2-2\gamma}{\alpha}} \delta^{2\gamma -2},\qquad
\zeta_{3} = \lvert z\rvert^{\frac{4\gamma -2}{\alpha}} \delta^{4-4\gamma}.
\end{equation}

Note that $\zeta_{3} = \delta^{-2}\zeta_{1}^{2}$. Throughout the paper, we fixed a small $\delta >0$, and took $\epsilon >0$ sufficiently small. In this case, we have that $\zeta_{3} = O(\zeta_{1}^{2})$ which is also used throughout. Recall also the definition of the region $Y$, and note that since $(x_{0}, y_{0})$ and $(x_{2}, y_{2})$ lie in $Y$, we have that $|w|$ and $|\overline{w}|$ lie in the interval $[m^{-1}, m]$, which is also used throughout this section.

We start with calculating $\frac{\partial f}{\partial w}$. By \eqref{eq:derivative_composition_map}, \eqref{eq_68igtfygbi} and \eqref{eq_657g7033ww}, we have
\begin{equation}\label{eq:y6q8re7}
\begin{aligned}
\frac{\partial f}{\partial w} &= \frac{- z^{\frac{1}{\alpha}}}{x_{2} w^{2}}\Big(\left[d_{1} b_{2} -b_{1} b_{2} \overline{w} + O(\zeta_{1})\right] \frac{\partial q_{21}}{\partial x_{1}} \cdot \frac{\partial q_{12}}{\partial x_{0}}\\
&\qquad\qquad + \left[d_{1} d_{2} - b_{1} d_{2} \overline{w} + O(\zeta_{1} + \zeta_{2})\right] \frac{\partial q_{21}}{\partial y_{1}} \cdot \frac{\partial q_{12}}{\partial x_{0}}\Big).
\end{aligned}
\end{equation}
To simplify this expression, note that by \eqref{eq56yhnto8y}, we have
\begin{equation*}
x_{2} = b_{1} v_{2\tau} \left[1 + O(\lvert u_{2\tau} v_{2\tau}^{-1}\rvert) + O(v_{2\tau}^{2})\right]
\quad \mathrm{and}\quad
y_{2} = d_{1} v_{2\tau} \left[1 + O(\lvert u_{2\tau} v_{2\tau}^{-1}\rvert) + O(v_{2\tau}^{2})\right],
\end{equation*}
where $u_{2\tau}$ and $v_{2\tau}$ are given by \eqref{eq_178rtf6y3h23} and \eqref{eq_45tuvkb7taguer}, i.e.
\begin{equation}\label{eq_45tcgiytfyd}
u_{2\tau} = \sigma_{1} \lvert y_{1}\rvert^{\frac{1}{\gamma}} \delta^{1 - \frac{1}{\gamma}} [1 + O(\delta^{2})]
\qquad\mathrm{and}\qquad
v_{2\tau} = \kappa \gamma\sigma_{1} \sigma_{2} x_{1} \lvert y_{1}\rvert^{1-\frac{1}{\gamma}} \delta^{\frac{1}{\gamma} -1} [1 + O(\delta^{2})].
\end{equation}
However, for $(x_{0}, y_{0}) \in Y$, we have
\begin{equation}\label{eq_546yvgt7d3}
x_{1} = a_{2} \sigma_{0} \lvert y_{0}\rvert^{\gamma} \delta^{1-\gamma} [1 + O(\delta^{2})]
\qquad\mathrm{and}\qquad
y_{1} = c_{2} \sigma_{0} \lvert y_{0}\rvert^{\gamma} \delta^{1-\gamma} [1 + O(\delta^{2})].
\end{equation}
Substituting these expressions for $x_{1}$ and $y_{1}$ into \eqref{eq_45tcgiytfyd}, we have
\begin{equation*}
d_{1} - b_{1} \overline{w} = d_{1} - \frac{b_{1} y_{2}}{x_{2}} = O(\lvert u_{2\tau} v_{2\tau}^{-1}\rvert) + O(v_{2\tau}^{2}) = O(\zeta_{1}^{2}) + O(\zeta_{2}).
\end{equation*}
By virtue of this relation, we can simplify \eqref{eq:y6q8re7} further and obtain
\begin{equation}\label{eq:f_w}
\frac{\partial f}{\partial w} = \frac{- z^{\frac{1}{\alpha}}}{x_{2} w^{2}}\left(O(\zeta_{1} + \zeta_{2})\cdot \frac{\partial q_{21}}{\partial x_{1}} \cdot \frac{\partial q_{12}}{\partial x_{0}} + O(\zeta_{1} + \zeta_{2})\cdot \frac{\partial q_{21}}{\partial y_{1}} \cdot \frac{\partial q_{12}}{\partial x_{0}}\right).
\end{equation}
We can also express $x_{2}$ and $y_{2}$ in terms of $x_{0}$ and $y_{0}$ using \eqref{eq_546yvgt7d3} as follows
\begin{equation}\label{eq:8huy7trwuj45}
\begin{aligned}
x_{2} &= b_{1}\kappa\gamma\sigma_{1} \sigma_{2} x_{1} \lvert y_{1}\rvert^{1-\frac{1}{\gamma}} \delta^{\frac{1}{\gamma} -1} [1 + O(\delta^{2})]\\
&= b_{1}\kappa \gamma \sigma_{0}\sigma_{1}\sigma_{2} a_{2} \lvert c_{2}\rvert^{1-\frac{1}{\gamma}} \lvert y_{0}\rvert^{2\gamma -1} \delta^{2 -2\gamma} [1 + O(\delta^{2})],\\[1em]
y_{2} &= d_{1}\kappa\gamma\sigma_{1} \sigma_{2} x_{1} \lvert y_{1}\rvert^{1-\frac{1}{\gamma}} \delta^{\frac{1}{\gamma} -1} [1 + O(\delta^{2})]\\
&= d_{1}\kappa \gamma \sigma_{0}\sigma_{1}\sigma_{2} a_{2} \lvert c_{2}\rvert^{1-\frac{1}{\gamma}} \lvert y_{0}\rvert^{2\gamma -1} \delta^{2 -2\gamma} [1 + O(\delta^{2})].
\end{aligned}
\end{equation}
In particular,
\begin{equation*}
\frac{1}{x_2} = \kappa \sigma_{0}\sigma_{1}\sigma_{2} \gamma^{-1} b_{1}^{-1} a_{2}^{-1} \lvert c_{2}\rvert^{\frac{1}{\gamma}-1} \zeta_{1}^{-1} \delta [1 + O(\delta^{2})].
\end{equation*}
On the other hand, by virtue of \eqref{eq_546yvgt7d3} and \eqref{eq_derivatives_local_map_21}, we have
\begin{equation}\label{eq:o8iil4}
\begin{aligned}
\frac{\partial q_{12}}{\partial x_{0}} &= \kappa \sigma_{0}\sigma_{1} \gamma^{-1} \zeta_{2}^{\frac{1}{2}} [1 + O(\delta^{2})],\\
\frac{\partial q_{21}}{\partial x_{1}} &= \kappa \sigma_{1} \sigma_{2}\gamma |c_{2}|^{1 -\frac{1}{\gamma}} \zeta_{2}^{-\frac{1}{2}} [1 + O(\delta^{2})],\\
\frac{\partial q_{21}}{\partial y_{1}} &= (\gamma-1)\kappa \sigma_{0} \sigma_{1} a_{2} \lvert c_{2}\rvert^{-\frac{1}{\gamma}} \zeta_{2}^{-\frac{1}{2}} [1 + O(\delta^{2})].
\end{aligned}
\end{equation}
Thus,
\begin{equation*}
\begin{aligned}
\frac{\partial q_{21}}{\partial x_{1}} \cdot \frac{\partial q_{12}}{\partial x_{0}} &= \sigma_{0}\sigma_{2} \lvert c_{2}\rvert^{1 -\frac{1}{\gamma}} [1 + O(\delta^{2})],\\
\frac{\partial q_{21}}{\partial y_{1}} \cdot \frac{\partial q_{12}}{\partial x_{0}} &= (1 - \frac{1}{\gamma}) a_{2} \lvert c_{2}\rvert^{-\frac{1}{\gamma}} [1 + O(\delta^{2})].
\end{aligned}
\end{equation*}
Substituting these estimates into \eqref{eq:f_w}, we obtain
\begin{equation}\label{eq_87igy78tif}
\frac{\partial f}{\partial w} = O(\zeta_{1}\zeta_{2} + \zeta_{2}^{2}).
\end{equation}

We now estimate $\frac{\partial f}{\partial z}$. First, note that by \eqref{eq_657g7033ww}, \eqref{eq_657g7033wz} and \eqref{eq_68igtfygbi}, we have
\begin{equation*}
\frac{\partial f}{\partial z} = \frac{-w}{\alpha z}\cdot \frac{\partial f}{\partial w} + \frac{1}{x_{2}}\cdot \frac{\partial y_{0}}{\partial z} \cdot \left(\frac{\partial y_{2}}{\partial y_{0}} - \frac{\partial x_{2}}{\partial y_{0}} \cdot \overline{w}\right).
\end{equation*}
Substituting the obtained estimates into this relation yields
\begin{equation*}
\frac{\partial f}{\partial z} = \frac{-w}{\alpha z}\cdot \frac{\partial f}{\partial w} + \frac{z^{\frac{1}{\alpha}-1}}{\alpha x_{2}}\left(O(\zeta_{1} + \zeta_{2})\cdot \frac{\partial q_{21}}{\partial x_{1}} \cdot \frac{\partial p_{12}}{\partial y_{0}} + O(\zeta_{1} + \zeta_{2})\cdot \frac{\partial q_{21}}{\partial y_{1}}\cdot \frac{\partial p_{12}}{\partial y_{0}}\right).
\end{equation*}
However, \eqref{eq_derivatives_local_map_12} yields $\frac{\partial p_{12}}{\partial y_{0}} = \gamma\sigma_{0}\sigma_{1} \zeta_{2}^{-\frac{1}{2}} [1 + O(\delta^{2})]$. Noting this relation and \eqref{eq:o8iil4}, we have
\begin{equation*}
\begin{aligned}
\frac{\partial q_{21}}{\partial x_{1}} \cdot \frac{\partial p_{12}}{\partial y_{0}} &= \kappa \gamma^{2}\sigma_{0} \sigma_{2} \lvert c_{2}\rvert^{1 -\frac{1}{\gamma}} \zeta_{2}^{-1} [1 + O(\delta^{2})],\\
\frac{\partial q_{21}}{\partial y_{1}}\cdot \frac{\partial p_{12}}{\partial y_{0}} &= \gamma (\gamma -1)\kappa a_{2} \lvert c_{2}\rvert^{-\frac{1}{\gamma}} \zeta_{2}^{-1} [1 + O(\delta^{2})].
\end{aligned}
\end{equation*}
Thus, taking \eqref{eq_87igy78tif} into account, we obtain
\begin{equation}\label{eq_48uuhg8ih8}
\frac{\partial f}{\partial z} = \frac{1}{\alpha z}\left[O(\zeta_{1}) + O(\zeta_{2})\right].
\end{equation}

To find $\frac{\partial g}{\partial w}$, notice that using the obtained estimates above as well as \eqref{eq:derivative_composition_map}, we have
\begin{equation}\label{eq_1t4e5ryguho}
\frac{\partial y_{2}}{\partial x_{0}} = \sigma_{0}\sigma_{2} \gamma^{-1} C_{1} [C_{2} + O(\delta^{2})],
\end{equation}
where
\begin{equation*}
C_{1} := d_{1} a_{2} \gamma\lvert c_{2}\rvert^{1-\frac{1}{\gamma}}
\qquad \mathrm{and}\qquad
C_{2} := \frac{b_{2}}{a_{2}} + \sigma_{0}\sigma_{2}(1 - \gamma^{-1}) \frac{d_{2}}{\lvert c_{2}\rvert}
\end{equation*}
are constants. Following \eqref{eq:8huy7trwuj45} and the relation $\alpha = \frac{1}{2N+1}$ (see \eqref{eq:formula_alpha}), we also have
\begin{equation}\label{eq:8w6tyg24}
y_{2}^{\alpha - 1} = C_{1}^{\alpha - 1} \zeta_{1}^{\alpha - 1} \delta^{1-\alpha} [1 + O(\delta^{2})].
\end{equation}
Substituting \eqref{eq_68igtfygbi}, \eqref{eq_1t4e5ryguho} and \eqref{eq:8w6tyg24} into \eqref{eq_657g7033zw} and noting that $z = \sigma_{1} |z|$ (see \eqref{eq_evolution_sigma}), we have
\begin{equation}\label{eq_iuyubr4gr3e}
\frac{\partial g}{\partial w} = -\gamma^{-1} \sigma_{0}\sigma_{1}\sigma_{2} \alpha C_{1}^{\alpha} w^{-2} \zeta_{2} \zeta_{3}^{\frac{\alpha}{2}} \left[C_{2} + O(\delta^{2})\right].
\end{equation}

Finally, we find $\frac{\partial g}{\partial z}$. By \eqref{eq_68igtfygbi} and \eqref{eq_67uohjbnh}, we have $\frac{\partial g}{\partial z} = z^{\frac{1}{\alpha} - 1} y_{2}^{\alpha - 1} (\frac{\partial y_{2}}{\partial x_{0}}\cdot w^{-1} + \frac{\partial y_{2}}{\partial y_{0}})$. Using the obtained estimates above as well as \eqref{eq:derivative_composition_map}, we obtain
\begin{equation}\label{eq:9yigyj3}
\frac{\partial y_{2}}{\partial y_{0}} = \sigma_{0} \sigma_{2} \kappa C_{1} \left[\gamma + \sigma_{0} \sigma_{2}\mathrm{sgn}(c_{2}) (\gamma -1)\right] \zeta_{2}^{-1} [1 + O(\delta^{2})].
\end{equation}
However, \eqref{eq_evolution_sigma} yields $\sigma_{2} = \sigma_{0}\mathrm{sgn}(c_{2})$. This simplifies \eqref{eq:9yigyj3} further, and gives
\begin{equation}\label{eq_576uyfvbi}
\frac{\partial y_{2}}{\partial y_{0}} = \sigma_{0} \sigma_{2} \kappa (2\gamma -1) C_{1} \zeta_{2}^{-1} [1 + O(\delta^{2})].
\end{equation}
Note that $2\gamma -1 \neq 0$ for $0.5 < \gamma <1$. By \eqref{eq_1t4e5ryguho} and \eqref{eq_576uyfvbi}, if $\epsilon > 0$ is sufficiently small (i.e. $\lvert z\rvert$ is sufficiently small), then $\frac{\partial y_{2}}{\partial x_{0}}\cdot w^{-1} = O(\zeta_{2}^{-1} \delta^{2})$. This gives
\begin{equation}\label{eq_derivative_zbar_z}
\frac{\partial g}{\partial z} = \sigma_{0} \sigma_{2} \kappa (2\gamma -1) C_{1}^{\alpha} \zeta_{3}^{\frac{\alpha}{2}} {|z|}^{-1} [1 + O(\delta^{2})].
\end{equation}

So far, we have estimated the differential of the map $\mathbf{T}$. We are now ready to prove Lemma \ref{Lem_909tfud34}.

\begin{proof}[Proof of Lemma \ref{Lem_909tfud34}]
We show that $\mathbf{T}$ admits a cross-form, and the cross-map $\mathbf{T}^{\times}$ satisfies the assumptions of Theorem 4.3\footnote{The original formulation of Afraimovich--Shilnikov theory \cite{AnnulusPrinciple} was developed on annuli for proving the existence of invariant tori (see also Chapter 4 of \cite{Dimabook}). The formulation that we use here is valid for Banach spaces as well.} in \cite{Dimabook}.

Note that $z=0$ corresponds to the point $(0,0)$ on the cross-section, and so $z=0$ if and only if $\overline{z} = 0$. Fixing $\delta$, for $z\neq 0$, relation \eqref{eq:8huy7trwuj45} yields
\begin{equation*}
|g(w,z)| \asymp |z|^{2\gamma-1},
\end{equation*}
i.e. there are $K_{1}, K_{2} >0$ such that $K_{1} |z|^{2\gamma-1} \leq |g(w,z)| \leq K_{2} |z|^{2\gamma-1}$. Consequently, $g$ extends continuously across zero by setting $g(w,0) = 0$, Moreover, it follows from \eqref{eq_derivative_zbar_z} that $\frac{\partial g}{\partial z}$ does not vanish on the domain of \eqref{eq_map_wz_coordinates}, its sign for both $z<0$ and $z>0$ is the same, and $|\frac{\partial g}{\partial z}| \to \infty$ as $z\to 0$. Thus, $z\mapsto g(w,z)$ is strictly monotone, and by virtue of the implicit function theorem, we can express $z$ in terms of $(w, \overline{z})$, i.e. $z = G(w, \overline{z})$ for some $C^{1}$ function $G$, and we have
\begin{equation*}
\pd{G}{\overline{z}}(w,\overline{z}) = \Bigl(\frac{\partial g}{\partial z}(w, z)\Bigr)^{-1} = \sigma_{0} \sigma_{2} \kappa (2\gamma -1)^{-1} C_{1}^{-\alpha} \zeta_{3}^{-\frac{\alpha}{2}} |z| [1 + O(\delta^{2})].
\end{equation*}
Having $G$ and defining $F(w, \overline{z}) := f(w, G(w, \overline{z}))$, we have that $\mathbf{T}$ admits a cross-form. Note also that substituting $z = G(w, \overline{z})$ into \eqref{eq_map_wz_coordinates} gives
\begin{align}
\overline{w} &= f(w, G(w, \overline{z})),\nonumber\\
\overline{z} &= g(w, G(w, \overline{z})),\label{eq_8i284ho}
\end{align}
By differentiating \eqref{eq_8i284ho} with respect to $w$ and taking \eqref{eq_iuyubr4gr3e} and \eqref{eq_derivative_zbar_z} into account, we have
\begin{equation*}
\begin{aligned}
\pd{G}{w}(w,\overline{z}) &= -\left(\pd{g}{z}(w, z)\right)^{-1} \cdot \pd{g}{w}(w, z)\\
&= \kappa \alpha \gamma^{-1} (2\gamma -1)^{-1} w^{-2} \zeta_{2} z [C_{2} + O(\delta^{2})].
\end{aligned}
\end{equation*}
Differentiating $F(w, \overline{z}) = f(w, G(w, \overline{z}))$ with respect to $w$ and $\overline{z}$, and using \eqref{eq_87igy78tif} and \eqref{eq_48uuhg8ih8}, we obtain
\begin{equation*}
\begin{aligned}
\pd{F}{w}(w, \overline{z}) &= \pd{f}{w}(w,z) + \pd{f}{z}(w,z)\cdot \pd{G}{w}(w, \overline{z})\\
&= O(\zeta_{1}\zeta_{2} + \zeta_{2}^{2})
\end{aligned}
\end{equation*}
and
\begin{equation*}
\begin{aligned}
\pd{F}{\overline{z}}(w, \overline{z}) &= \pd{f}{z}(w,z)\cdot \pd{G}{\overline{z}}(w, \overline{z})\\
&= \alpha^{-1} C_{1}^{-\alpha} O(\zeta_{1}^{1-\alpha} + \zeta_{1}^{-\alpha} \zeta_{2}).
\end{aligned}
\end{equation*}

\begin{figure}[!tbh]
\centering
\includegraphics[width=\linewidth]{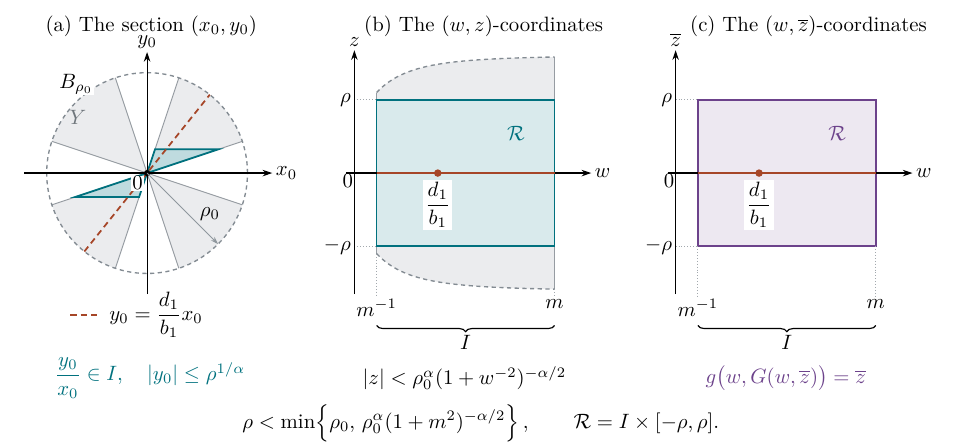}
\caption{Domains for the cross-map construction, shown for $d_1/b_1>0$ and $I=[m^{-1},m]$.
\textup{(a)} The gray sectors are $Y\cap B_{\rho_0}$; the teal region is specified by the displayed inequalities, with the origin included. The dashed line has slope $d_1/b_1$.
\textup{(b)} The rectangle $\mathcal R$ in $(w,z)$-coordinates. The curved boundaries correspond to $\partial B_{\rho_0}$, and the highlighted segment $I\times\{0\}$ represents the section origin.
\textup{(c)} The same rectangle in the cross-map variables $(w,\overline{z})$.}	
\label{fig:cross_map_domains}
\end{figure}

We need a subset $\mathcal{R}$ of the domain of $\mathbf{T}$ in the form of a rectangle for which $\mathbf{T}^{\times}(\mathcal{R}) \subset \mathcal{R}$. Recall that the image of the points in the domain under $\mathbf{T}$ accumulates near $\overline{w} = \frac{d_1}{b_1}$. Let $I$ be one of the two intervals $[-m, -m^{-1}]$ and $[m^{-1}, m]$ that contains $\frac{d_1}{b_1}$. We then can take $\mathcal{R}$ to be $I\times [-\rho,\rho]$, where $\rho$ is as in \eqref{eq:radius_rectanglular_domain}; we have contraction in $w$ and strong expansion in $z$ which implies $\mathbf{T}^{\times}(\mathcal{R})\subset \mathcal{R}$. Let
\begin{equation*}
\begin{alignedat}{2}
K_{1} &:= \sup_{(w,\overline{z}) \in \mathcal{R}}\left\|\pd{F}{w}\right\| \cdot \left\|\pd{G}{\overline{z}}\right\|,\qquad
& K_{2} := \sup_{(w,\overline{z}) \in \mathcal{R}}\left\|\pd{F}{\overline{z}}\right\|,\\
K_{3} &:= \sup_{(w,\overline{z}) \in \mathcal{R}}\left\|\pd{G}{w}\right\|,\qquad
& K_{4} := \sup_{(w,\overline{z}) \in \mathcal{R}}\left\|\pd{F}{w}\right\|.
\end{alignedat}
\end{equation*}
For $(w,\overline{z}) \in \mathcal{R}$, we have $|z| \leq |\overline{z}| \leq \rho$. Noting this and \eqref{eq_zetas}, we have that for any $\alpha > 0$ the values of $K_1$, $K_3$, and $K_4$ can become arbitrarily small by choosing $\rho > 0$ small enough. We further have that if
\begin{equation}\label{eq:appropriate_alpha}
0 < \alpha <
\min\left\{\frac{2-2\gamma}{2\gamma-1},1\right\},
\end{equation}
then $K_2$ becomes arbitrarily small by choosing $\rho > 0$ small enough. Therefore, choosing any $\alpha$ satisfying \eqref{eq:appropriate_alpha}, we can choose $\epsilon > 0$ sufficiently small (recall that $\rho < \epsilon$) such that $\sqrt{K_1}+\sqrt{K_{2} K_{3}}<1$ and $K_{4}+\sqrt{K_{2}K_{3}}<1$.

For $w\in I$, define
\begin{equation*}
F(w,0)=\frac{d_1}{b_1}, \qquad G(w,0)=0.
\end{equation*}
The preceding asymptotics show that these extensions are continuous. Under condition~\eqref{eq:appropriate_alpha}, the derivative estimates above also give
\begin{equation*}
\frac{\partial F}{\partial w},\quad
\frac{\partial F}{\partial \overline{z}},\quad
\frac{\partial G}{\partial w},\quad
\frac{\partial G}{\partial \overline{z}}
\longrightarrow 0
\qquad\text{as }\overline{z}\to 0,
\end{equation*}
uniformly for $w\in I$. Consequently, $F$ and $G$ extend as $C^1$ functions across $\overline{z}=0$, with all four partial derivatives vanishing there.

Thus, the assumptions of Theorem~4.3 of \cite{Dimabook} are fulfilled, yielding an invariant curve of the form $\{w=W(z)\}$ for some $C^1$ Lipschitz function $W$. Taking this curve to the $(x,y)$-coordinates gives the desired unstable invariant curve. However, as the change of coordinates
\eqref{eq:definition_w_z_variables} is singular at the section origin $(x,y)=(0,0)$, to finalize the proofs of parts~\eqref{item:unstable_curve_1} and~\eqref{item:unstable_curve_2} of Lemma~\ref{Lem_909tfud34}, we need to show that the resulting curve in the $(x,y)$-coordinates is $C^1$ at $(0,0)$.

The curve $\{w=W(z)\}$ in $(x,y)$-coordinates is given by $x(y) = \frac{y}{W(y^{\alpha})}$. Note that $W(0) = \frac{d_{1}}{b_{1}} \neq 0$ by \eqref{eq:extended_map_in_wz}. Taking derivative, we have
\begin{equation*}
x^{\prime}(y) = \frac{W(y^{\alpha}) - \alpha y^{\alpha} W^{\prime}(y^{\alpha})}{[W(y^{\alpha})]^{2}},
\end{equation*}
which gives $x^{\prime}(y) \rightarrow \frac{b_{1}}{d_{1}}$ as $y\rightarrow 0$. This proves the $C^1$ regularity of the invariant curve including at the section origin in $(x,y)$-coordinates, as desired.

Note that we showed the $C^1$-smoothness of the curves on type-1 cross-sections. The $C^1$-smoothness of the invariant curves on the type-2 sections are obtained by applying the flow of the system on the invariant curves on type-1 sections.

The branches on the type-2 sections are obtained by applying the corresponding transition maps to the type-1 branches. By the derivative estimates in Lemma 5.1 and the derivative of \(T_2^{\mathrm{glo}}\), their tangent slopes converge to \(c_2/a_2\) at the section origin; hence the two branches join there to form a \(C^1\) curve.

Part \eqref{item:unstable_curve_3} of Lemma~\ref{Lem_909tfud34} follows from Theorem 4.3 of \cite{Dimabook} as well. More precisely, a careful examination of the proof of Theorem 4.2 (the original Afraimovich--Shilnikov formulation) of \cite{Dimabook} (see also Proposition D.2 of \cite{Bakrani2022JDE}) gives the following statement.

\begin{mylem}\label{lem:distance_from_invariant_curve}
Let the rectangle $\mathcal{R}$ be bounded. Then for every $\epsilon^{\prime}>0$, there exists an integer $N>0$ such that for every $n\geq N$ and every point $(w, z) \in \mathcal{R}$ for which $\{\mathbf{T}^{i}(w, z)\}_{i=0}^{n} \subset \mathcal{R}$, we have that the distance of $\mathbf{T}^{n}(w, z)$ from the invariant curve is at most $\epsilon^{\prime}$.
\end{mylem}
Following this statement, and since $\epsilon^{\prime}>0$ is arbitrary, if all the backward iterations of a point in $\mathcal{R}$ under $\mathbf{T}$ remain in $\mathcal{R}$, it has to belong to the invariant curve. This completes the proof of Lemma~\ref{Lem_909tfud34}.
\end{proof}

\section{Proof of Theorem \ref{thm:superhomoclinic}}\label{sec:proof_superhomoclinic}

We prove Theorem \ref{thm:superhomoclinic} in this section. The proof is based on the cross-map analysis in Section \ref{sec:Afraimovich_Shilnikov} and the convergence supplied by Lemma \ref{lem:distance_from_invariant_curve}.

\paragraph{Preliminary setting.}
For $\sigma\in\{+,-\}$, let $G^u_{1\sigma}$ and $G^s_{1\sigma}$ be the curves given by Lemmas~\ref{Lem_909tfud34} and \ref{Lem_0ouihiu}. These are the traces on $\Pi^{\mathrm{in}}_{1\sigma}$ of the local unstable and stable manifolds of the bouquet. They lie in $Y$ and $R^s \cup \{(0,0)\}$, respectively.

Let $r_u$ and $r_s$ be as in \eqref{eq:bumps_per_visit}. By \eqref{eq_evolution_sigma} and \eqref{eq:0oh8ygi34w}, the first return to $\Pi^{\mathrm{in}}_{1\sigma}$ in forward time on $Y$ is $T^{r_u}$, and the first return in backward time on $R^s$ is $T^{-r_s}$, where $T$ is the Poincar\'e map (see Definition \ref{dfn:Poincare_map_T}).

Let $\rho$ be as in \eqref{eq:radius_rectanglular_domain}. Choose $r >0$ such that $r^{\alpha} < \rho$. Thus, $r < \rho < \epsilon < \delta$. This choice ensures compatibility with the cross-map rectangle: $|v_2|\leq r$ implies $|z|=|v_2|^\alpha\leq r^\alpha<\rho$. Write $\overline B_{r}$ for the closed ball of radius $r$ centered at $(0,0)$ on $\Pi^{\mathrm{in}}_{1\sigma}$. Decreasing $r$ further if necessary, each branch of $G^u_{1\sigma} \setminus\{(0,0)\}$ and $G^s_{1\sigma}\setminus\{(0,0)\}$ meets $\partial B_{r}$ transversely in exactly one point. We choose $r$ sufficiently small for the analogous time-reversed construction as well.

\paragraph{Convergence of spanning arcs.}
Consider a simple arc $l_0\subset\overline B_{r}$ passing through $(0,0)$, with its two endpoints on $\partial B_{r}$ and all remaining points in $B_{r}$. Suppose that the two branches of $l_0\setminus\{(0,0)\}$ lie in opposite cone components of $Y$. Thus, $l_0$ spans the ball $\overline B_{r}$. Construct $l_{n+1}$ from $l_n$ by applying $T^{r_u}$ on the two punctured branches of $l_n\setminus\{(0,0)\}$, taking their closures at $(0,0)$, and retaining on each image branch the portion from $(0,0)$ to its first intersection with $\partial B_{r}$. The expansion in Lemma~\ref{Lem_66fiytfi5}, together with the local domains established in Section~\ref{sec:existence_invariant_manifolds}, ensures that these image branches reach the boundary $\partial B_{r}$. Consequently, every $l_n$ is again a spanning arc, and $l_n\subset Y$.

We claim that these arcs approach $G^{u}_{1\sigma}$ uniformly, in the sense that
\begin{equation}\label{eq:C_spanning_convergence}
\sup_{p\in l_n}\operatorname{dist}(p,G^u_{1\sigma}) \longrightarrow0\qquad\text{as }n\to\infty.
\end{equation}
Recall the map $\mathbf{T}$ given by \eqref{eq:extended_map_in_wz}. This map is the representation of $T^{2}$ restricted to $\Pi^{\mathrm{in}}_{1\sigma}\cap Y$ in $(w,z)$-coordinates. In the proof of Lemma \ref{Lem_909tfud34}, we proved that $\mathbf{T}$ admits a cross-form, and introduced a bounded rectangle $\mathcal{R}$ for which, $\mathbf{T}$ together with its cross-map and the rectangle $\mathcal{R}$ fulfill the assumptions of Theorem 4.3 in \cite{Dimabook}. Consequently, Lemma \ref{lem:distance_from_invariant_curve} asserts that \eqref{eq:C_spanning_convergence} holds in $(w,z)$-coordinates which readily implies the uniform convergence in $(u_{2}, v_{2})$-coordinates as well. This proves the claim.

By reversing time, the analogous construction with $T^{-r_s}$ applied to a spanning arc in $R^s$ gives arcs approaching $G^s_{1\sigma}$ uniformly in distance. Here the two initial branches lie in the two components of
$R^s\cap\Pi^{\mathrm{in}}_{1\sigma}$, and closure at the section origin is understood as above.

\paragraph{Approximating arcs from the equilibrium manifolds.}
The normal form~\eqref{eq:normal_form} straightens the local unstable manifold of $O$ to $\{u_1=u_2=0\}$. Its intersection with $\Pi^{\mathrm{out}}_{1\sigma}$ is the $v_2$-axis. The image of this axis under $T_1^{\mathrm{glo}}$ is a curve in $W^u_{\mathrm{glo}}(O)$ passing through $(0,0)$ on $\Pi^{\mathrm{in}}_{1\sigma}$, with tangent slope $d_1/b_1$. Let $W_0^u$ be its spanning arc in $\overline B_{r}$. Its punctured branches lie in opposite cone components of $Y$. Applying $T^{r_u}$ with the truncation just described gives spanning arcs
$W_n^u\subset W^u_{\mathrm{glo}}(O)$ approaching $G^u_{1\sigma}$.

Similarly, the local stable manifold of $O$ is $\{v_1=v_2=0\}$ and intersects $\Pi^{\mathrm{in}}_{1\sigma}$ in the $u_2$-axis. Let $W_0^s$ be the spanning segment of this axis in $\overline B_{r}$. Its punctured branches lie in $R^s$. Applying $T^{-r_s}$ with the same truncation convention gives spanning arcs $W_m^s\subset W^s_{\mathrm{glo}}(O)$ approaching $G^s_{1\sigma}$. In particular,
\begin{equation*}
\begin{aligned}
\sup_{p\in W_n^u}\operatorname{dist}(p,G^u_{1\sigma})
&\longrightarrow0 &&\text{as }n\to\infty,\\
\sup_{p\in W_m^s}\operatorname{dist}(p,G^s_{1\sigma})
&\longrightarrow0 &&\text{as }m\to\infty.
\end{aligned}
\end{equation*}
Each iteration in these constructions adds exactly one visit to the chosen cross-section.

\paragraph{Crossing argument.}
Let $\Gamma$ be a super-homoclinic orbit supplied by Assumption~\ref{assum:exist_superhomoclinic}. Its backward and forward tails intersect some sections $\Pi^{\mathrm{in}}_{1\sigma_1}$ and $\Pi^{\mathrm{in}}_{1\sigma_2}$ infinitely often along $G^u_{1\sigma_1}$ and $G^s_{1\sigma_2}$ for some $\sigma_{1},\sigma_{2} = \pm$, respectively, with the intersection points tending to the corresponding section origins. Choose tail points
\begin{equation*}
 p_u\in\Gamma\cap G^u_{1\sigma_1},
 \qquad p_s\in\Gamma\cap G^s_{1\sigma_2},
\end{equation*}
strictly inside the corresponding $r$-balls, away from the section origins, and with $p_u$ preceding $p_s$ in forward time. The trajectories along $\Gamma$ define a Poincar\'e map
\begin{equation*}
 P:U_u\longrightarrow U_{s}, \qquad P(p_u)=p_s,
\end{equation*}
where $U_u\subset\Pi^{\mathrm{in}}_{1\sigma_1}$ and $U_s\subset\Pi^{\mathrm{in}}_{1\sigma_2}$ are neighborhoods of $p_u$ and $p_s$. From now on, $W_n^u$ denotes the approximating arcs on $\Pi^{\mathrm{in}}_{1\sigma_1}$, and $W_m^s$ those on $\Pi^{\mathrm{in}}_{1\sigma_2}$.

\begin{figure}[!tbh]
\centering
\includegraphics[width=\linewidth]{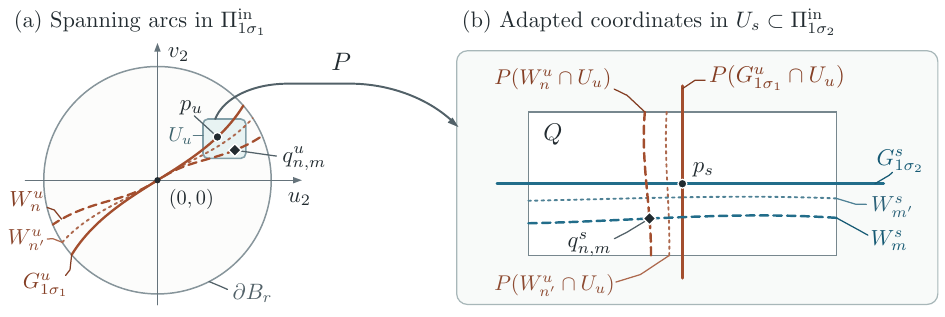}
\caption{Crossing construction for Theorem~\ref{thm:superhomoclinic}.
\textup{(a)} In the original section coordinates, the arcs $W_n^u$ span $\overline B_r$ and approach $G^u_{1\sigma_1}$. The passage $P:U_u\to U_s$ along $\Gamma$ sends $p_u$ to $p_s$.
\textup{(b)} In adapted coordinates, the solid bouquet traces are the axes through $p_s$, the center of $Q$. For sufficiently large $n,m$, a subarc of $P(W_n^u\cap U_u)$ joins the top and bottom of $Q$, and a subarc of $W_m^s$ joins its left and right sides. Short dashes show selected later arcs ($n'>n$, $m'>m$) approaching the corresponding solid traces. The diamonds mark an intersection $q^s_{n,m}$ of the long-dashed arcs and its preimage under $P$, the point $q^u_{n,m}$, on the same homoclinic orbit to $O$.}
\label{fig:theorem_C_crossing}
\end{figure}

The transverse intersection of the bouquet manifolds implies that $P(G^u_{1\sigma_1}\cap U_u)$ and $G^s_{1\sigma_2}$ intersect transversely at $p_s$. Choose adapted $C^1$ local coordinates on $U_s$ in which $p_s=(0,0)$ and these two curves are, respectively, the vertical and horizontal coordinate axes (see Figure \ref{fig:theorem_C_crossing}). The origin in these adapted local coordinates is $p_s$, not the cross-section origin.

In the adapted coordinates, choose a small rectangle $Q$ centered at $p_s$ that $P(G^u_{1\sigma_1}\cap U_u)$ and $G^s_{1\sigma_2}$ extend beyond its corresponding sides. For sufficiently large $n$ and $m$, the approximating arcs $P(W^u_n \cap U_u)$ and $W^s_m$ contain subarcs joining the top and bottom sides and the left and right sides of $Q$, respectively. Thus, there exist $N_{u}$ and $N_{s}$ such that for any $n \geq N_{u}$ and $m \geq N_{s}$ the curves $P(W^{u}_{n} \cap U_u)$ and $W^{s}_{m}$ intersect at a point
\begin{equation*}
q_{n,m}^{s}\in W^u_{\mathrm{glo}}(O)\cap W^s_{\mathrm{glo}}(O)
\end{equation*}
on a homoclinic orbit to $O$. Let $q_{n,m}^{u} := P^{-1}(q_{n,m}^{s})$. We show below that the homoclinic orbit passing through these two points is the desired $\Gamma_{n,m}$ in Theorem \ref{thm:superhomoclinic}.

\paragraph{Visit and bump counts.} Taking a point on $W_0^u$ (excluding the section origin on $\Pi^{\mathrm{in}}_{1\sigma_{1}}$) and moving backward in time, we take a single bump along $\Gamma_{1\sigma_{1}}$ and once we reach $\Pi^{\mathrm{out}}_{1\sigma_{1}}$, we are already on the local unstable manifold of $O$ and hence we converge to $O$. For $i\geq 1$, tracing a point on $W_i^u$ backward to $W_{i-1}^u$ completes $r_u$ bumps along the bouquet. During this passage, the successive incoming-section labels follow the code provided by Theorem~\ref{thm:coding_unstable} in reverse order. This implies that, moving backward in time, a point on $W_i^u$ visits $\Pi^{\mathrm{in}}_{1\sigma_{1}}$ exactly $i$ times (excluding the starting point) and contributes $i r_{u} +1$ bumps. 

Regarding the forward dynamics, taking a point on $W_0^s$ (excluding the section origin on $\Pi^{\mathrm{in}}_{1\sigma_{2}}$), we are already on the local stable manifold of $O$ and hence we converge to $O$ directly. For $i\geq 1$, if we take a point on $W_i^s$ and move forward in time, we take $r_s$ bumps along the bouquet according to the coding provided by Theorem \ref{thm:coding_stable}, and reach $W_{i-1}^s$. This implies that, moving forward in time, a point on $W_i^s$ visits $\Pi^{\mathrm{in}}_{1\sigma_{2}}$ exactly $i$ times (excluding the starting point) and contributes $i r_{s}$ bumps. This completes the proof of Theorem~\ref{thm:superhomoclinic}.

\begin{myrem}
The same crossing argument shows that, for all sufficiently large $n$, $P(W_n^u\cap U_u)$ intersects $G^s_{1\sigma_2}$. Such an intersection lies on an orbit converging to $O$ in backward time and asymptotic to the bouquet $B$ in forward time. Similarly, intersections of $P(G^u_{1\sigma_1}\cap U_u)$ with $W_m^s$ give orbits asymptotic to $B$ in backward time and converging to $O$ in forward time.
\end{myrem}

\appendix

\section{Proof of Lemma \ref{Lem_derivatives_local_map_12}}\label{Sec_Proof_Lem_derivatives_local_map_12}

Let us drop the subscript $12$ from $p$ and $q$, so
\begin{equation*}
T^{\mathrm{loc}}_{1\sigma_{0} 2\sigma_{1}}(u_{20}, v_{20}) = \left(p(u_{20}, v_{20}), q(u_{20}, v_{20})\right).
\end{equation*}
Let $(u_{1}(r), u_{2}(r), v_{1}(r), v_{2}(r))$, where $r=(t, \tau, u_{10}, u_{20}, v_{1\tau}, v_{2\tau})$, be the solution of system \eqref{eq:normal_form} that satisfies boundary conditions \eqref{eq_boundary_conditions_1}. Here, we have $u_{10} = \sigma_{0}\delta$ and $v_{2\tau} = \sigma_{1} \delta$. Recall from Lemma \ref{Lem_coordinates_on_cross_sections} that $v_{10}$ is a function of $(u_{20}, v_{20})$, i.e.
\begin{equation*}
v_{10} = h^{\mathrm{in}}_{1\sigma_{0}}(u_{20}, v_{20}) = \kappa \sigma_{0} \gamma^{-1} \delta^{-1} u_{20} v_{20} [1 + O(\delta^{2})] + O\left(\delta v_{20}^{2}\right).
\end{equation*}
Moreover, the traveling time $\tau$ is a function of $(u_{20}, v_{20})$, too, i.e. $\tau = \tau(u_{20}, v_{20})$. We then have
\begin{equation*}
\begin{aligned}
p(u_{20}, v_{20}) &= u_{1}\left(\tau(u_{20}, v_{20}), \tau(u_{20}, v_{20}), \sigma_{0}\delta, u_{20}, q(u_{20}, v_{20}), \sigma_{1}\delta\right),\\
h^{\mathrm{in}}_{1\sigma_{0}}(u_{20}, v_{20}) &= v_{1}\left(0, \tau(u_{20}, v_{20}), \sigma_{0}\delta, u_{20}, q(u_{20}, v_{20}), \sigma_{1}\delta\right),\\
v_{20} &= v_{2}\left(0, \tau(u_{20}, v_{20}), \sigma_{0}\delta, u_{20}, q(u_{20}, v_{20}), \sigma_{1}\delta\right).
\end{aligned}
\end{equation*}
Differentiating these equations with respect to $u_{20}$ and $v_{20}$ yields
\begin{align}
\frac{\partial p}{\partial u_{20}} &= \left(\frac{\partial u_{1}}{\partial t}\Big\vert_{t=\tau} + \frac{\partial u_{1}}{\partial \tau}\Big\vert_{t=\tau}\right)\cdot \frac{\partial \tau}{\partial u_{20}} + \frac{\partial u_{1}}{\partial u_{20}}\Big\vert_{t=\tau} + \frac{\partial u_{1}}{\partial v_{1\tau}}\Big\vert_{t=\tau}\cdot \frac{\partial q}{\partial u_{20}},\label{eq6iggggg6}\\
\frac{\partial p}{\partial v_{20}} &= \left(\frac{\partial u_{1}}{\partial t}\Big\vert_{t=\tau} + \frac{\partial u_{1}}{\partial \tau}\Big\vert_{t=\tau}\right)\cdot \frac{\partial \tau}{\partial v_{20}} + \frac{\partial u_{1}}{\partial v_{1\tau}}\Big\vert_{t=\tau}\cdot \frac{\partial q}{\partial v_{20}},\label{eq6iggggg5}\\
\frac{\partial h^{\mathrm{in}}_{1\sigma_{0}}}{\partial u_{20}} &= \frac{\partial v_{1}}{\partial \tau}\Big\vert_{t=0} \cdot \frac{\partial \tau}{\partial u_{20}} + \frac{\partial v_{1}}{\partial u_{20}}\Big\vert_{t=0} + \frac{\partial v_{1}}{\partial v_{1\tau}}\Big\vert_{t=0} \cdot \frac{\partial q}{\partial u_{20}},\label{eq6iggggg4}\\
\frac{\partial h^{\mathrm{in}}_{1\sigma_{0}}}{\partial v_{20}} &= \frac{\partial v_{1}}{\partial \tau}\Big\vert_{t=0} \cdot \frac{\partial \tau}{\partial v_{20}} + \frac{\partial v_{1}}{\partial v_{1\tau}}\Big\vert_{t=0} \cdot \frac{\partial q}{\partial v_{20}},\label{eq6iggggg2}\\
0 &= \frac{\partial v_{2}}{\partial \tau}\Big\vert_{t=0} \cdot \frac{\partial \tau}{\partial u_{20}} + \frac{\partial v_{2}}{\partial u_{20}}\Big\vert_{t=0} + \frac{\partial v_{2}}{\partial v_{1\tau}}\Big\vert_{t=0} \cdot \frac{\partial q}{\partial u_{20}},\label{eq6iggggg3}\\
1 &= \frac{\partial v_{2}}{\partial \tau}\Big\vert_{t=0} \cdot \frac{\partial \tau}{\partial v_{20}} + \frac{\partial v_{2}}{\partial v_{1\tau}}\Big\vert_{t=0} \cdot \frac{\partial q}{\partial v_{20}}.\label{eq6iggggg1}
\end{align}
We find the estimates in Lemma \ref{Lem_derivatives_local_map_12} by finding the necessary estimates in the above formulas. To simplify the calculations, we also define a new variable
\begin{equation*}
\eta := \lvert v_{20} \rvert \delta^{-1}.
\end{equation*}
Note also that, since $\sgn(v_{20}) = \sigma_{1}$, then $v_{20} = \sigma_{1}\eta\delta$.

We start with $\frac{\partial u_{1}}{\partial t}\big\vert_{t=\tau}$. By virtue of \eqref{eq:normal_form}, we have
\begin{equation*}
\begin{aligned}
\frac{\partial u_{1}}{\partial t}\Big\vert_{t=\tau} &= -\lambda_{1} u_{1\tau} + O\left(\lvert u_{1\tau}^{2} v_{1\tau}\rvert + \lvert u_{1\tau} u_{2\tau} \delta\rvert + \lvert u_{2\tau} v_{1\tau} \delta\rvert + \lvert u_{1\tau}^{3} \delta^{2}\rvert + \lvert u_{2\tau}^{2} v_{1\tau}^{3}\rvert + \lvert u_{1\tau} u_{2\tau}^{2} v_{1\tau}^{2}\rvert\right)\\
&= -\lambda_{1} u_{1\tau}[1 + O(\delta^{2})] + O\left(\lvert u_{2\tau} v_{1\tau} \delta\rvert\right)\\
&= -\sigma_{0}\lambda_{1}\delta \eta^{\gamma} [1 + O(\delta^{2})],
\end{aligned}
\end{equation*}
where the second line is derived using \eqref{eq_45tu9yo} and \eqref{eq_67iby20iubhv}, and noting that $u_{2\tau} = O(\delta)$.

Regarding $\frac{\partial h^{\mathrm{in}}_{1\sigma_{0}}}{\partial u_{20}}$ and $\frac{\partial h^{\mathrm{in}}_{1\sigma_{0}}}{\partial v_{20}}$, by virtue of \eqref{eq56yigkre4r} and for $(u_{20}, v_{20}) \in Y$, we have
\begin{equation*}
\frac{\partial h^{\mathrm{in}}_{1\sigma_{0}}}{\partial u_{20}} = \kappa \sigma_{0} \gamma^{-1} \delta^{-1} v_{20} [1 + O(\delta^{2})]
\qquad \mathrm{and}\qquad
\frac{\partial h^{\mathrm{in}}_{1\sigma_{0}}}{\partial v_{20}} = \kappa \sigma_{0} \gamma^{-1} \delta^{-1} u_{20} [1 + O(\delta^{2})].
\end{equation*}

To find the other required estimates, write system \eqref{eq:normal_form} in the form
\begin{equation}\label{eq405rfctf7tu}
\dot{u}_{i} = -\lambda_{i} u_{i} + F_{i}(u_{1}, u_{2}, v_{1}, v_{2}),\qquad
\dot{v}_{i} = \lambda_{i} v_{i} + G_{i}(u_{1}, u_{2}, v_{1}, v_{2}),
\end{equation}
for $i=1,2$, where $F_{i}$ and $G_{i}$, together with their first derivatives, vanish at the origin $O$. Let $x = (u_{1}(r), u_{2}(r), v_{1}(r), v_{2}(r))$, where $r=(t, \tau, u_{10}, u_{20}, v_{1\tau}, v_{2\tau})$, be the solution of system \eqref{eq:normal_form} satisfying \eqref{eq_boundary_conditions_1}. Let $\theta \in \{\tau, u_{10}, u_{20}, v_{1\tau}, v_{2\tau}\}$, and define
\begin{equation*}
U_{i} := \frac{\partial u_{i}(r)}{\partial \theta},\qquad
V_{i} := \frac{\partial v_{i}(r)}{\partial \theta}.
\end{equation*}
Differentiating \eqref{eq405rfctf7tu} with respect to $\theta$ gives the variational equation
\begin{equation*}
\begin{aligned}
\dot{U}_{i} &= -\lambda_{i} U_{i} + \sum_{j=1}^{2} \left(\frac{\partial F_{i}}{\partial u_{j}}(x) \cdot U_{j} + \frac{\partial F_{i}}{\partial v_{j}}(x) \cdot V_{j}\right),\\
\dot{V}_{i} &= \lambda_{i} V_{i} + \sum_{j=1}^{2} \left(\frac{\partial G_{i}}{\partial u_{j}}(x) \cdot U_{j} + \frac{\partial G_{i}}{\partial v_{j}}(x) \cdot V_{j}\right)
\end{aligned}
\end{equation*}
for $i=1,2$. Similar to the case discussed in Section~ \ref{sec:local_analysis_near_equilibrium}, this system has a unique solution $(U_{1}(t), U_{2}(t), V_{1}(t), V_{2}(t))$ satisfying boundary conditions
\begin{equation}\label{eq:boundary_conditions_derivative}
U_{1}(0) = U_{10},\qquad  U_{2}(0) = U_{20},\qquad  V_{1}(\tau) = V_{1\tau},\qquad  V_{2}(\tau) = V_{2\tau},
\end{equation}
for any real $U_{10}$, $U_{20}$, $V_{1\tau}$ and $V_{2\tau}$ (see, e.g. Theorem 2.10 of \cite{Dimabook} and its proof). Moreover, this solution is the fixed point of the integral operator
\begin{equation*}
\mathfrak{T}: (U_{1}(t), U_{2}(t), V_{1}(t), V_{2}(t)) \mapsto (\overline{U}_{1}(t), \overline{U}_{2}(t), \overline{V}_{1}(t), \overline{V}_{2}(t)),
\end{equation*}
defined by
\begin{equation}\label{eq43bt4re}
\begin{aligned}
\overline{U}_{i}(t) &= e^{-\lambda_{i} t} U_{i0} + \int_{0}^{t} e^{\lambda_{i}(s-t)}\sum_{j=1}^{2} \left(\frac{\partial F_{i}}{\partial u_{j}}(x(s)) \cdot U_{j}(s) + \frac{\partial F_{i}}{\partial v_{j}}(x(s)) \cdot V_{j}(s)\right) ds,\\
\overline{V}_{i}(t) &= e^{-\lambda_{i} (\tau - t)} V_{i\tau} - \int_{t}^{\tau} e^{\lambda_{i}(t-s)} \sum_{j=1}^{2} \left(\frac{\partial G_{i}}{\partial u_{j}}(x(s)) \cdot U_{j}(s) + \frac{\partial G_{i}}{\partial v_{j}}(x(s)) \cdot V_{j}(s)\right) ds.
\end{aligned}
\end{equation}
The operator $\mathfrak{T}$ is a contraction and its fixed point can be estimated by iterations of $\mathfrak{T}$ \cite{PoincareBirkhoff, Dimabook}. In Section \ref{sec:proofs_technical_lemmas}, we use this method to prove the following technical lemmas.
\begin{mylem}\label{Lem_620yu3tyohk}
The following hold.
\begin{enumerate}[(i)]
\item For the case of $\frac{\partial u_{1}}{\partial u_{20}}$, $\frac{\partial v_{1}}{\partial u_{20}}$ and $\frac{\partial v_{2}}{\partial u_{20}}$, the boundary conditions are $U_{10} = V_{1\tau} = V_{2\tau} = 0$ and $U_{20} = 1$. Solving this boundary value problem, for $i=1,2$, we obtain
\begin{equation*}
\begin{gathered}
\frac{\partial u_{1}}{\partial u_{20}}\Big\vert_{t=\tau} = e^{-\lambda_{2}\tau} O(\delta^{2}),\qquad
\frac{\partial u_{2}}{\partial u_{20}}\Big\vert_{t=\tau} = e^{-\lambda_{2}\tau} [1 + O(\delta^{2})],\qquad
\frac{\partial v_{i}}{\partial u_{20}}\Big\vert_{t=0} = e^{-\lambda_{2}\tau} O(\delta^{2}).
\end{gathered}
\end{equation*}

\item For the case of $\frac{\partial u_{1}}{\partial v_{1\tau}}$, $\frac{\partial v_{1}}{\partial v_{1\tau}}$ and $\frac{\partial v_{2}}{\partial v_{1\tau}}$, the boundary conditions are $U_{10} = U_{20} = V_{2\tau} = 0$ and $V_{1\tau} = 1$. Solving this boundary value problem, for $i=1,2$, we obtain
\begin{equation*}
\frac{\partial u_{i}}{\partial v_{1\tau}}\Big\vert_{t=\tau} = e^{-\lambda_{2}\tau} O(\delta^{2}),\qquad
\frac{\partial v_{1}}{\partial v_{1\tau}}\Big\vert_{t=0} = e^{-\lambda_{1}\tau} + e^{- \lambda_{2}\tau} O(\delta^{2}),\qquad
\frac{\partial v_{2}}{\partial v_{1\tau}}\Big\vert_{t=0} = e^{-\lambda_{2}\tau} O(\delta^{2}).
\end{equation*}
\end{enumerate}
\end{mylem}

To deal with the derivatives with respect to $\tau$, note that it follows from the first equation of \eqref{eq5676fuyg8} that $U_{10} = U_{20} = 0$. Differentiating the second equation of \eqref{eq5676fuyg8} with respect to $\tau$ also yields $V_{i\tau} = -\lambda_{i} v_{i\tau} - G_{i}(u_{1\tau}, u_{2\tau}, v_{1\tau}, v_{2\tau})$ for $i=1,2$. By virtue of the form of system \eqref{eq:normal_form} as well as the conditions $u_{10} = \sigma_{0}\delta$, $v_{2\tau} = \sigma_{1} \delta$ and $(u_{20}, v_{20})\in Y$, we have
\begin{equation*}
\begin{aligned}
G_{1}(u_{1\tau}, u_{2\tau}, v_{1\tau}, v_{2\tau}) = O(\lvert v_{1\tau}\rvert \delta^{2})\qquad \mathrm{and}\qquad
G_{2}(u_{1\tau}, u_{2\tau}, v_{1\tau}, v_{2\tau}) = O(\delta^{3}).
\end{aligned}
\end{equation*}
Thus, the boundary conditions are
\begin{equation}\label{eq_098iuh}
U_{10} = 0,\qquad
U_{20} = 0,\qquad
V_{1\tau} = -\lambda_{1} v_{1\tau}[1 + O(\delta^{2})],\qquad
V_{2\tau} = -\sigma_{1}\lambda_{2} \delta[1 + O(\delta^{2})].
\end{equation}
\begin{mylem}\label{Lem_derivatives_resp_tau_12}
Assume $(u_{20}, v_{20})\in Y$. Solving the boundary value problem \eqref{eq_098iuh}, we obtain
\begin{equation*}
\frac{\partial u_{i}}{\partial \tau}\Big\vert_{t=\tau} = e^{-\lambda_{i}\tau} O(\delta^{3}),\qquad
\frac{\partial v_{i}}{\partial \tau}\Big\vert_{t=0} = -\lambda_{i} e^{-\lambda_{i}\tau} v_{i\tau} [1 + O(\delta^{2})],\qquad (i=1,2).
\end{equation*}
\end{mylem}
We prove this lemma in Section \ref{sec:proofs_technical_lemmas}. Taking \eqref{eq567ug79osw} into account, Lemmas \ref{Lem_620yu3tyohk} and \ref{Lem_derivatives_resp_tau_12} yield
\begin{equation*}
\begin{gathered}
\frac{\partial u_{1}}{\partial u_{20}}\Big\vert_{t=\tau} = O(\eta \delta^{2}),\qquad
\frac{\partial u_{2}}{\partial u_{20}}\Big\vert_{t=\tau} = \eta [1 + O(\delta^{2})],\qquad
\frac{\partial v_{i}}{\partial u_{20}}\Big\vert_{t=0} = O(\eta \delta^{2}),\\
\frac{\partial u_{i}}{\partial v_{1\tau}}\Big\vert_{t=\tau} = O(\eta \delta^{2}),\qquad
\frac{\partial v_{1}}{\partial v_{1\tau}}\Big\vert_{t=0} = \eta^{\gamma} [1 + O(\delta^{2})],\qquad
\frac{\partial v_{2}}{\partial v_{1\tau}}\Big\vert_{t=0} = O(\eta \delta^{2}),\\
\frac{\partial u_{1}}{\partial \tau}\Big\vert_{t=\tau} = O\left(\eta^{\gamma} \delta^{3}\right),\qquad
\frac{\partial u_{2}}{\partial \tau}\Big\vert_{t=\tau} = O\left(\eta \delta^{3}\right),\\
\frac{\partial v_{1}}{\partial \tau}\Big\vert_{t=0} = -\lambda_{2} \kappa \sigma_{0}\sigma_{1} u_{20} \eta [1 + O(\delta^{2})],\qquad
\frac{\partial v_{2}}{\partial \tau}\Big\vert_{t=0} = -\lambda_{2} \sigma_{1}\eta\delta [1 + O(\delta^{2})], \qquad (i=1,2).
\end{gathered}
\end{equation*}

So far, we have found the estimates that we needed to obtain the estimates in Lemma \ref{Lem_derivatives_local_map_12}. By substituting the obtained estimates into \eqref{eq6iggggg1}, we get
\begin{equation*}
1 = -\lambda_{2} \sigma_{1}\eta\delta [1 + O(\delta^{2})] \cdot \frac{\partial \tau}{\partial v_{20}} + O(\eta \delta^{2}) \cdot \frac{\partial q}{\partial v_{20}}.
\end{equation*}
This implies
\begin{equation}\label{eq_68yivhbsbu}
\frac{\partial \tau}{\partial v_{20}} = \frac{-\sigma_{1}}{\lambda_{2} \eta\delta}[1 + O(\delta^{2})] + O(\delta) \cdot \frac{\partial q}{\partial v_{20}}.
\end{equation}
Substituting this estimate and the other estimates obtained above into \eqref{eq6iggggg2} gives
\begin{equation*}
\begin{aligned}
\kappa \gamma^{-1} \sigma_{0} \delta^{-1} u_{20} [1 + O(\delta^{2})] &= -\lambda_{2} \kappa \sigma_{0}\sigma_{1} u_{20} \eta [1 + O(\delta^{2})] \cdot \left(\frac{-\sigma_{1}}{\lambda_{2} \eta\delta}[1 + O(\delta^{2})] + O(\delta) \cdot \frac{\partial q}{\partial v_{20}}\right)\\
&\quad+ \eta^{\gamma} [1 + O(\delta^{2})]\cdot \frac{\partial q}{\partial v_{20}}.
\end{aligned}
\end{equation*}
Noting that $(u_{20}, v_{20}) \in Y$, estimate \eqref{eq_derivatives_local_map_12} for $\frac{\partial q}{\partial v_{20}}$ is obtained from the above equation. Then, by \eqref{eq_68yivhbsbu}, we get
\begin{equation*}
\frac{\partial \tau}{\partial v_{20}} = \frac{-\sigma_{1}}{\lambda_{2} \eta \delta}[1 + O(\delta^{2})].
\end{equation*}

Substituting the obtained estimates above into \eqref{eq6iggggg3} gives
\begin{equation*}
0 = -\lambda_{2} \sigma_{1}\eta\delta [1 + O(\delta^{2})] \cdot \frac{\partial \tau}{\partial u_{20}} + O(\eta \delta^{2}) + O(\eta \delta^{2}) \cdot \frac{\partial q}{\partial u_{20}}
\end{equation*}
and therefore $\frac{\partial \tau}{\partial u_{20}} = O(\delta) + O(\delta) \cdot \frac{\partial q}{\partial u_{20}}$. Using these estimates, relation \eqref{eq6iggggg4} becomes
\begin{equation*}
\begin{aligned}
\kappa\gamma^{-1} \sigma_{0}\sigma_{1} \eta [1 + O(\delta^{2})] &= -\lambda_{2} \kappa \sigma_{0}\sigma_{1} u_{20} \eta [1 + O(\delta^{2})] \cdot \left(O(\delta) + O(\delta) \cdot \frac{\partial q}{\partial u_{20}}\right)\\
&\quad + O(\eta \delta^{2}) + \eta^{\gamma} [1 + O(\delta^{2})] \cdot \frac{\partial q}{\partial u_{20}}.
\end{aligned}
\end{equation*}
Having $(u_{20}, v_{20}) \in Y$, simplifying this equation gives estimate \eqref{eq_derivatives_local_map_12} for $\frac{\partial q}{\partial u_{20}}$. We also get $\frac{\partial \tau}{\partial u_{20}} = O(\delta)$.

Substituting the obtained estimates above into \eqref{eq6iggggg5} gives
\begin{equation*}
\begin{aligned}
\frac{\partial p}{\partial v_{20}} &= \left[-\sigma_{0}\lambda_{1} \eta^{\gamma}\delta [1 + O(\delta^{2})] + O\left(\eta^{\gamma} \delta^{3}\right)\right]\cdot \frac{-\sigma_{1}}{\lambda_{2} \eta\delta}[1 + O(\delta^{2})]\\
&\quad+ O(\eta \delta^{2})\cdot (1-\gamma)\kappa \gamma^{-1} \sigma_{0} u_{20} \eta^{-\gamma} \delta^{-1}[1 + O(\delta^{2})].
\end{aligned}
\end{equation*}
Simplifying this expression gives estimate \eqref{eq_derivatives_local_map_12} for $\frac{\partial p}{\partial v_{20}}$. Similarly, using the obtained estimates, \eqref{eq6iggggg6} gives
\begin{equation*}
\begin{aligned}
\frac{\partial p}{\partial u_{20}} &= \left[-\sigma_{0}\lambda_{1}\delta \eta^{\gamma} [1 + O(\delta^{2})] + O\left(\eta^{\gamma} \delta^{3}\right)\right]\cdot O(\delta) + O(\eta \delta^{2})\\
&\quad+ O(\eta \delta^{2})\cdot \kappa \gamma^{-1} \sigma_{0}\sigma_{1} \eta^{1-\gamma} [1 + O(\delta^{2})].
\end{aligned}
\end{equation*}
Simplifying this expression gives estimate \eqref{eq_derivatives_local_map_12} for $\frac{\partial p}{\partial u_{20}}$. This finishes the proof of Lemma \ref{Lem_derivatives_local_map_12}.

\section{Proof of Lemma \ref{Lem_derivatives_local_map_21}}\label{Sec_Proof_Lem_derivatives_local_map_21}

The proof of Lemma \ref{Lem_derivatives_local_map_21} is similar to the proof of Lemma \ref{Lem_derivatives_local_map_12}. For simplicity, we drop the subscript $21$ from $p$ and $q$, i.e.
\begin{equation*}
T^{\mathrm{loc}}_{2\sigma_{1} 1\sigma_{2}}(u_{10}, v_{10}) = \left(p(u_{10}, v_{10}), q(u_{10}, v_{10})\right).
\end{equation*}
Let $(u_{1}(r), u_{2}(r), v_{1}(r), v_{2}(r))$, where $r=(t, \tau, u_{10}, u_{20}, v_{1\tau}, v_{2\tau})$, be the solution of system \eqref{eq:normal_form} that satisfies boundary conditions \eqref{eq_boundary_conditions_1}. Here, we have $u_{20} = \sigma_{1}\delta$ and $v_{1\tau} = \sigma_{2} \delta$. Recall from Lemma \ref{Lem_coordinates_on_cross_sections} that $v_{20}$ is a function of $(u_{10}, v_{10})$, i.e.
\begin{equation*}
v_{20} = h^{\mathrm{in}}_{2\sigma_{1}}(u_{10}, v_{10}) = \kappa \sigma_{1} \gamma \delta^{-1} u_{10} v_{10} [1 + O(\delta^{2})] + O\left(\delta v_{10}^{2}\right).
\end{equation*}
Moreover, the traveling time $\tau$ is a function of $(u_{10}, v_{10})$, too, i.e. $\tau = \tau(u_{10}, v_{10})$. We then have
\begin{equation*}
\begin{aligned}
p(u_{10}, v_{10}) &= u_{2}\left(\tau(u_{10}, v_{10}), \tau(u_{10}, v_{10}), u_{10}, \sigma_{1}\delta, \sigma_{2}\delta, q(u_{10}, v_{10})\right),\\
h^{\mathrm{in}}_{2\sigma_{1}}(u_{10}, v_{10}) &= v_{2}\left(0, \tau(u_{10}, v_{10}), u_{10}, \sigma_{1}\delta, \sigma_{2}\delta, q(u_{10}, v_{10})\right),\\
v_{10} &= v_{1}\left(0, \tau(u_{10}, v_{10}), u_{10}, \sigma_{1}\delta, \sigma_{2}\delta, q(u_{10}, v_{10})\right).
\end{aligned}
\end{equation*}
Differentiating these equations with respect to $u_{10}$ and $v_{10}$ yields
\begin{align}
\frac{\partial p}{\partial u_{10}} &= \left(\frac{\partial u_{2}}{\partial t}\Big\vert_{t=\tau} + \frac{\partial u_{2}}{\partial \tau}\Big\vert_{t=\tau}\right)\cdot \frac{\partial \tau}{\partial u_{10}} + \frac{\partial u_{2}}{\partial u_{10}}\Big\vert_{t=\tau} + \frac{\partial u_{2}}{\partial v_{2\tau}}\Big\vert_{t=\tau}\cdot \frac{\partial q}{\partial u_{10}},\label{eq45tv55}\\
\frac{\partial p}{\partial v_{10}} &= \left(\frac{\partial u_{2}}{\partial t}\Big\vert_{t=\tau} + \frac{\partial u_{2}}{\partial \tau}\Big\vert_{t=\tau}\right)\cdot \frac{\partial \tau}{\partial v_{10}} + \frac{\partial u_{2}}{\partial v_{2\tau}}\Big\vert_{t=\tau}\cdot \frac{\partial q}{\partial v_{10}},\label{eq45tv66}\\
\frac{\partial h^{\mathrm{in}}_{2\sigma_{1}}}{\partial u_{10}} &= \frac{\partial v_{2}}{\partial \tau}\Big\vert_{t=0}\cdot \frac{\partial \tau}{\partial u_{10}} + \frac{\partial v_{2}}{\partial u_{10}}\Big\vert_{t=0} + \frac{\partial v_{2}}{\partial v_{2\tau}}\Big\vert_{t=0} \cdot\frac{\partial q}{\partial u_{10}},\label{eq45tv11}\\
\frac{\partial h^{\mathrm{in}}_{2\sigma_{1}}}{\partial v_{10}} &= \frac{\partial v_{2}}{\partial \tau}\Big\vert_{t=0}\cdot \frac{\partial \tau}{\partial v_{10}} + \frac{\partial v_{2}}{\partial v_{2\tau}}\Big\vert_{t=0} \cdot\frac{\partial q}{\partial v_{10}},\label{eq45tv22}\\
0 &= \frac{\partial v_{1}}{\partial \tau}\Big\vert_{t=0}\cdot \frac{\partial \tau}{\partial u_{10}} + \frac{\partial v_{1}}{\partial u_{10}}\Big\vert_{t=0} + \frac{\partial v_{1}}{\partial v_{2\tau}}\Big\vert_{t=0} \cdot\frac{\partial q}{\partial u_{10}},\label{eq45tv33}\\
1 &= \frac{\partial v_{1}}{\partial \tau}\Big\vert_{t=0}\cdot \frac{\partial \tau}{\partial v_{10}} + \frac{\partial v_{1}}{\partial v_{2\tau}}\Big\vert_{t=0} \cdot\frac{\partial q}{\partial v_{10}}.\label{eq45tv44}
\end{align}
Similar to the proof of Lemma \ref{Lem_derivatives_local_map_12}, we find the estimates in Lemma \ref{Lem_derivatives_local_map_21} by finding the necessary estimates in the above formulas. To simplify the calculations, we also define a new variable
\begin{equation*}
\eta := \lvert v_{10} \rvert \delta^{-1}.
\end{equation*}
Note also that since $\sgn(v_{10}) = \sigma_{2}$, then $v_{10} = \sigma_{2}\eta\delta$.

We start with estimating $\frac{\partial u_{2}}{\partial t}\big\vert_{t=\tau}$. By virtue of \eqref{eq:normal_form}, we have
\begin{equation*}
\begin{aligned}
\frac{\partial u_{2}}{\partial t}\Big\vert_{t=\tau} &= -\lambda_{2} u_{2\tau} + O\Big(\lvert u_{1\tau}^{2} u_{2\tau}\rvert + \lvert u_{1\tau} u_{2\tau} \delta\rvert + \lvert u_{2\tau}^{3}\rvert + \lvert u_{2\tau}^{2} v_{2\tau}\rvert\\
&\quad+ \lvert u_{1\tau}^{2} \delta^{2} v_{2\tau}\rvert + \lvert u_{1\tau}^{4} v_{2\tau}\rvert + \lvert u_{1\tau}^{3} \delta v_{2\tau}\rvert\Big).
\end{aligned}
\end{equation*}
On the other hand, \eqref{eq6yfigk8ihou} yields $u_{1\tau} = O(e^{-\lambda_{1}\tau} \delta)$. Thus, by \eqref{eq_178rtf6y3h23} and \eqref{eq_45tuvkb7taguer}, for $(u_{10}, v_{10}) \in Y$, we obtain
\begin{equation*}
\frac{\partial u_{2}}{\partial t}\Big\vert_{t=\tau} = -\lambda_{2} u_{2\tau} [1 + O(\delta^{2})] + O(\lvert u_{1\tau}^{2} \delta^{2} v_{2\tau}\rvert) = -\sigma_{1} \lambda_{2} \eta^{\frac{1}{\gamma}} \delta [1 + O(\delta^{2})].
\end{equation*}
Regarding $\frac{\partial h^{\mathrm{in}}_{2\sigma_{1}}}{\partial u_{10}}$ and $\frac{\partial h^{\mathrm{in}}_{2\sigma_{1}}}{\partial v_{10}}$, by virtue of \eqref{eq56yigkre4r} and for $(u_{10}, v_{10}) \in Y$, we have
\begin{equation*}
\frac{\partial h^{\mathrm{in}}_{2\sigma_{1}}}{\partial u_{10}} = \kappa \sigma_{1} \gamma \delta^{-1} v_{10} [1 + O(\delta^{2})]
\qquad \mathrm{and}\qquad
\frac{\partial h^{\mathrm{in}}_{2\sigma_{1}}}{\partial v_{10}} = \kappa \sigma_{1} \gamma \delta^{-1} u_{10} [1 + O(\delta^{2})].
\end{equation*}

As in the proof of Lemma \ref{Lem_derivatives_local_map_12}, the other necessary estimates are found by solving some boundary value problems. We start with the following lemma which is proved in Section \ref{sec:proofs_technical_lemmas}.

\begin{mylem}\label{Lem_78iuhr4}
Recall boundary conditions \eqref{eq:boundary_conditions_derivative} and the integral operator $\mathfrak{T}$ given by \eqref{eq43bt4re}.
\begin{enumerate}[(i)]
\item For the case of $\frac{\partial u_{2}}{\partial u_{10}}$, $\frac{\partial v_{1}}{\partial u_{10}}$ and $\frac{\partial v_{2}}{\partial u_{10}}$, the boundary conditions are $U_{10} = 1$ and $U_{20} = V_{1\tau} = V_{2\tau} = 0$. Solving this boundary value problem, for $i=1,2$, we obtain
\begin{equation*}
\frac{\partial u_{1}}{\partial u_{10}}\Big\vert_{t=\tau} = e^{-\lambda_{1}\tau} [1 + O(\delta^{2})],\qquad
\frac{\partial u_{2}}{\partial u_{10}}\Big\vert_{t=\tau} = e^{-\lambda_{2}\tau} O(\delta^{2}),\qquad
\frac{\partial v_{i}}{\partial u_{10}}\Big\vert_{t=0} = e^{-\lambda_{2}\tau} O(\delta^{2}).
\end{equation*}
\item For the case of $\frac{\partial u_{2}}{\partial v_{2\tau}}$, $\frac{\partial v_{1}}{\partial v_{2\tau}}$ and $\frac{\partial v_{2}}{\partial v_{2\tau}}$, the boundary conditions are $U_{10} = U_{20} = V_{1\tau} = 0$ and $V_{2\tau} = 1$. Solving this boundary value problem, for $i=1,2$, we obtain
\begin{equation*}
\frac{\partial u_{i}}{\partial v_{2\tau}}\Big\vert_{t=\tau} = e^{-\lambda_{2}\tau} O(\delta^{2}),\qquad
\frac{\partial v_{1}}{\partial v_{2\tau}}\Big\vert_{t=0} = e^{-\lambda_{2}\tau} O(\delta^{2}),\qquad
\frac{\partial v_{2}}{\partial v_{2\tau}}\Big\vert_{t=0} = e^{-\lambda_{2}\tau} [1 + O(\delta^{2})].
\end{equation*}
\end{enumerate}
\end{mylem}

We now estimate $\frac{\partial u_{2}}{\partial \tau}$, $\frac{\partial v_{1}}{\partial \tau}$ and $\frac{\partial v_{2}}{\partial \tau}$. The first equation of \eqref{eq5676fuyg8} implies $U_{10} = U_{20} = 0$. Differentiating the second equation of \eqref{eq5676fuyg8} with respect to $\tau$ yields $V_{i\tau} = -\lambda_{i} v_{i\tau} - G_{i}(u_{1\tau}, u_{2\tau}, v_{1\tau}, v_{2\tau})$ for $i=1,2$. Considering the form of system \eqref{eq:normal_form} as well as the conditions $u_{20} = \sigma_{1}\delta$, $v_{1\tau} = \sigma_{2} \delta$ and $(u_{10}, v_{10})\in Y$, we obtain
\begin{equation*}
\begin{aligned}
G_{1}(u_{1\tau}, u_{2\tau}, v_{1\tau}, v_{2\tau}) = O(\delta^{3})\qquad \mathrm{and}\qquad
G_{2}(u_{1\tau}, u_{2\tau}, v_{1\tau}, v_{2\tau}) = O(\lvert v_{2\tau}\rvert \delta^{2}).
\end{aligned}
\end{equation*}
Thus, the boundary conditions are
\begin{equation}\label{eq_1359y7hirg}
U_{10} = 0,\qquad
U_{20} = 0,\qquad
V_{1\tau} = -\sigma_{2}\lambda_{1} \delta [1 + O(\delta^2)],\qquad
V_{2\tau} = -\lambda_{2} v_{2\tau} [1 + O(\delta^2)].
\end{equation}

\begin{mylem}\label{Lem_derivatives_resp_tau_21}
Assume $(u_{10}, v_{10})\in Y$. Solving the boundary value problem \eqref{eq_1359y7hirg}, we obtain
\begin{equation*}
\frac{\partial u_{i}}{\partial \tau}\Big\vert_{t=\tau} = e^{-\lambda_{i}\tau} O\left(\delta^{3}\right),\qquad
\frac{\partial v_{i}}{\partial \tau}\Big\vert_{t=0} = -\lambda_{i} e^{-\lambda_{i}\tau} v_{i\tau} [1 + O(\delta^{2})],\qquad (i=1,2).
\end{equation*}
\end{mylem}
We prove this lemma in Section \ref{sec:proofs_technical_lemmas}. Taking \eqref{eq9er578iuytry} into account, Lemmas \ref{Lem_78iuhr4} and \ref{Lem_derivatives_resp_tau_21} yield
\begin{equation*}
\begin{gathered}
\frac{\partial u_{1}}{\partial u_{10}}\Big\vert_{t=\tau} = \eta [1 + O(\delta^{2})],\qquad
\frac{\partial u_{2}}{\partial u_{10}}\Big\vert_{t=\tau} = O(\eta^{\frac{1}{\gamma}} \delta^{2}),\qquad
\frac{\partial v_{i}}{\partial u_{10}}\Big\vert_{t=0} = O(\eta^{\frac{1}{\gamma}} \delta^{2}),\\
\frac{\partial u_{i}}{\partial v_{2\tau}}\Big\vert_{t=\tau} = O(\eta^{\frac{1}{\gamma}} \delta^{2}),\qquad
\frac{\partial v_{1}}{\partial v_{2\tau}}\Big\vert_{t=0} = O(\eta^{\frac{1}{\gamma}} \delta^{2}),\qquad
\frac{\partial v_{2}}{\partial v_{2\tau}}\Big\vert_{t=0} = \eta^{\frac{1}{\gamma}} [1 + O(\delta^{2})],\\
\frac{\partial u_{1}}{\partial \tau}\Big\vert_{t=\tau} = O(\eta \delta^{3}),
\qquad\qquad
\frac{\partial u_{2}}{\partial \tau}\Big\vert_{t=\tau} = O(\eta^{\frac{1}{\gamma}} \delta^{3}),\\
\frac{\partial v_{1}}{\partial \tau}\Big\vert_{t=0} = -\lambda_{1} \sigma_{2}\eta\delta [1 + O(\delta^{2})],\qquad
\frac{\partial v_{2}}{\partial \tau}\Big\vert_{t=0} = -\lambda_{1} \kappa\sigma_{1} \sigma_{2} u_{10} \eta [1 + O(\delta^{2})],\qquad (i=1,2).
\end{gathered}
\end{equation*}

So far, we have found the estimates that we needed to obtain the estimates in Lemma \ref{Lem_derivatives_local_map_21}. By substituting the obtained estimates into \eqref{eq45tv33}, we get
\begin{equation*}
0 = -\lambda_{1} \sigma_{2}\eta\delta [1 + O(\delta^{2})]\cdot \frac{\partial \tau}{\partial u_{10}} + O(\eta^{\frac{1}{\gamma}} \delta^{2}) + O(\eta^{\frac{1}{\gamma}} \delta^{2}) \cdot\frac{\partial q}{\partial u_{10}}.
\end{equation*}
This implies $\frac{\partial \tau}{\partial u_{10}} = \eta^{\frac{1}{\gamma} -1} \delta [O(1) + O(1) \cdot\frac{\partial q}{\partial u_{10}}]$. Substituting these estimates into \eqref{eq45tv11} yields
\begin{equation*}
\begin{aligned}
\kappa \sigma_{1} \sigma_{2} \gamma\eta [1 + O(\delta^{2})] &= -\lambda_{1} \kappa\sigma_{1} \sigma_{2} u_{10} \eta [1 + O(\delta^{2})] \cdot \eta^{\frac{1}{\gamma} -1} \delta\left(O(1) + O(1) \cdot\frac{\partial q}{\partial u_{10}}\right)\\
&\quad + O(\eta^{\frac{1}{\gamma}} \delta^{2}) + \eta^{\frac{1}{\gamma}} [1 + O(\delta^{2})] \cdot\frac{\partial q}{\partial u_{10}}.
\end{aligned}
\end{equation*}
Simplifying this expression gives estimate \eqref{eq_derivatives_local_map_21} for $\frac{\partial q}{\partial u_{10}}$. We also obtain $\frac{\partial \tau}{\partial u_{10}} = O(\delta)$.

Substituting the obtained estimates into \eqref{eq45tv44} yields
\begin{equation*}
1 = -\lambda_{1} \sigma_{2}\eta\delta [1 + O(\delta^{2})]\cdot \frac{\partial \tau}{\partial v_{10}} + O(\eta^{\frac{1}{\gamma}} \delta^2)\cdot\frac{\partial q}{\partial v_{10}}
\end{equation*}
and therefore
\begin{equation*}
\frac{\partial \tau}{\partial v_{10}} = \frac{-\sigma_{2}}{\lambda_{1}\eta\delta}[1 + O(\delta^{2})] + O(\eta^{\frac{1}{\gamma} -1} \delta)\cdot\frac{\partial q}{\partial v_{10}}.
\end{equation*}
Substituting the obtained estimates into \eqref{eq45tv22}, we have
\begin{equation*}
\begin{aligned}
\frac{\kappa \gamma \sigma_{1} u_{10}}{\delta} [1 + O(\delta^{2})] &= -\lambda_{1} \kappa\sigma_{1} \sigma_{2} u_{10} \eta [1 + O(\delta^{2})] \left(\frac{-\sigma_{2}}{\lambda_{1}\eta\delta}[1 + O(\delta^{2})] + O(\eta^{\frac{1}{\gamma} -1} \delta)\cdot\frac{\partial q}{\partial v_{10}}\right)\\
&\quad + \eta^{\frac{1}{\gamma}} [1 + O(\delta^{2})] \cdot\frac{\partial q}{\partial v_{10}}.
\end{aligned}
\end{equation*}
Simplifying this expression gives estimate \eqref{eq_derivatives_local_map_21} for $\frac{\partial q}{\partial v_{10}}$. We also have $\frac{\partial \tau}{\partial v_{10}} = \frac{-1}{\lambda_{1} v_{10}}[1 + O(\delta^{2})]$.

Substituting the obtained estimates into \eqref{eq45tv55}, we have
\begin{equation*}
\begin{aligned}
\frac{\partial p}{\partial u_{10}} &= \left(-\sigma_{1} \lambda_{2} \eta^{\frac{1}{\gamma}} \delta [1 + O(\delta^{2})] + O(\eta^{\frac{1}{\gamma}} \delta^{3})\right)\cdot O(\delta) + O(\eta^{\frac{1}{\gamma}} \delta^{2})\\
&\quad+ O(\eta^{\frac{1}{\gamma}} \delta^{2})\cdot \kappa\gamma \sigma_{1} \sigma_{2} \eta^{1 -\frac{1}{\gamma}} [1 + O(\delta^{2})].
\end{aligned}
\end{equation*}
Simplifying this expression gives estimate \eqref{eq_derivatives_local_map_21} for $\frac{\partial p}{\partial u_{10}}$. Similarly, by the obtained estimates, \eqref{eq45tv66} gives
\begin{equation*}
\begin{aligned}
\frac{\partial p}{\partial v_{10}} &= \left(-\sigma_{1} \lambda_{2} \eta^{\frac{1}{\gamma}} \delta [1 + O(\delta^{2})] + O(\eta^{\frac{1}{\gamma}} \delta^{3})\right)\cdot \frac{-\sigma_{2}}{\lambda_{1} \eta\delta}[1 + O(\delta^{2})]\\
&\quad + O(\eta^{\frac{1}{\gamma}} \delta^{2})\cdot (\gamma -1)\kappa \sigma_{1} u_{10} \delta^{-1} \eta^{-\frac{1}{\gamma}} [1 + O(\delta^{2})].
\end{aligned}
\end{equation*}
Simplifying this expression gives estimate \eqref{eq_derivatives_local_map_21} for $\frac{\partial p}{\partial v_{10}}$. This ends the proof of Lemma \ref{Lem_derivatives_local_map_21}.

\section{Proofs of technical lemmas}\label{sec:proofs_technical_lemmas}

This section is dedicated to the proofs of Lemmas \ref{lem:formula_trajectories}, \ref{Lem_620yu3tyohk}, \ref{Lem_derivatives_resp_tau_12}, \ref{Lem_78iuhr4} and \ref{Lem_derivatives_resp_tau_21}.

For the Shilnikov boundary value problems below, we give the integral operators, the trial bounds and the estimates obtained from them. We use Maple to expand the polynomial majorants, integrate the resulting exponential terms and simplify the expressions for the solutions and their first derivatives. All these calculations have been implemented symbolically. Worked hand calculations for similar normal forms are given in \cite[Section~3.3.2 and the proof of Lemma~3.34]{BakraniPhDthesis}. We begin with the trajectory estimates.

\begin{proof}[Proof of Lemma \ref{lem:formula_trajectories}]
We adopt the notation of Section \ref{sec:NF_and_Shilnikov_boundary_value}. Consider the set
\begin{equation}\label{eq67uyg2hoeu}
\mathcal{A} = \left\{x(t): \quad \lvert u_{1}(t)\rvert \leq \eta_{1}, \quad \lvert u_{2}(t)\rvert \leq \eta_{2}, \quad \lvert v_{1}(t)\rvert \leq \nu_{1}, \quad \lvert v_{2}(t)\rvert \leq \nu_{2}\right\},
\end{equation}
where $x(t) := (u_{1}(t), u_{2}(t), v_{1}(t), v_{2}(t)) \in \mathcal{I}$, i.e. a continuous function from $[0,\tau]$ to $\Omega$, and
\begin{equation}\label{eq_87iu45kfe}
\begin{aligned}
\eta_{1} &:= 2 e^{-\lambda_{1}t}\lvert u_{10}\rvert + e^{-2\lambda_{2}t - \lambda_{1}(\tau - t)} \left\vert u_{20} v_{1\tau}\right\vert,\\
\eta_{2} &:= 2 e^{-\lambda_{2}t} \lvert u_{20}\rvert + e^{-2\lambda_{1}t - \lambda_{2}(\tau - t)} \lvert u_{10} v_{2\tau}\rvert,\\
\nu_{1} &:= 2 e^{-\lambda_{1}(\tau - t)}\lvert v_{1\tau}\rvert + e^{-2\lambda_{2}(\tau - t) - \lambda_{1}t} \lvert v_{2\tau} u_{10}\rvert,\\
\nu_{2} &:= 2 e^{-\lambda_{2}(\tau - t)} \lvert v_{2\tau}\rvert + e^{-2\lambda_{1}(\tau - t) - \lambda_{2}t} \lvert v_{1\tau} u_{20}\rvert.
\end{aligned}
\end{equation}
We claim that the set $\mathcal{A}$ is invariant under the integral operator $\mathfrak{T}$ given by \eqref{eq5676fuyg8}, i.e. $\mathfrak{T}(\mathcal{A}) \subset \mathcal{A}$. Recall system \eqref{eq:normal_form}, and let $M_{1} > 0$ be an upper bound for all the $\lvert f_{ij} \rvert$ and $\lvert g_{ij} \rvert$ on $\Omega$. Let $x(t) \in \mathcal{A}$. Taking $F_{i}$ and $G_{i}$ from \eqref{eq56guyhw} into account, we then have
\begin{equation*}
\left\vert F_{i}\left(x(t)\right)\right\vert \leq M_{1} P_{i},\qquad
\left\vert G_{i}\left(x(t)\right)\right\vert \leq M_{1} Q_{i},\qquad
(i=1,2),
\end{equation*}
where
\begin{align*}
P_{1} &:= \eta_{1}^{2} \nu_{1} +  \eta_{1} \eta_{2} \nu_{2} + \eta_{2}^{3} \nu_{1} \nu_{2} +  \eta_{2}^{2} \nu_{1} \nu_{2}^{2} +  \eta_{1}^{3} \nu_{2}^{2} +  \eta_{2}^{2} \nu_{1}^{3} +  \eta_{1} \eta_{2}^{2} \nu_{1}^{2},\\
P_{2} &:= \eta_{1}^{2}\eta_{2} + \eta_{1}\eta_{2}\nu_{1} + \eta_{2}^{3} + \eta_{2}^{2}\nu_{2} + \eta_{1}^{2} \nu_{1}^{2} \nu_{2} + \eta_{1}^{4}\nu_{2} + \eta_{1}^{3} \nu_{1} \nu_{2}\\
Q_{1} &:= \nu_{1}^{2}\eta_{1} + \nu_{1} \nu_{2} \eta_{2} + \nu_{2}^{3} \eta_{1} \eta_{2} + \nu_{2}^{2} \eta_{1} \eta_{2}^{2} + \nu_{1}^{3}\eta_{2}^{2} + \nu_{2}^{2} \eta_{1}^{3} + \nu_{1} \nu_{2}^{2} \eta_{1}^{2},\\
Q_{2} &:= \nu_{1}^{2} \nu_{2} + \nu_{1} \nu_{2} \eta_{1} + \nu_{2}^{3} + \nu_{2}^{2} \eta_{2} + \nu_{1}^{2} \eta_{1}^{2} \eta_{2} + \nu_{1}^{4} \eta_{2} + \nu_{1}^{3} \eta_{1} \eta_{2}.
\end{align*}
Taking the integral operator $\mathfrak{T}$ given by\eqref{eq5676fuyg8} into account, we obtain
\begin{equation*}
\begin{aligned}
\lvert\overline{u}_{i}(t) - e^{-\lambda_{i}t} u_{i0}\rvert
&\leq M_{1} e^{-\lambda_{i}t} \int_{0}^{t} e^{\lambda_{i}s} P_{i} ds,\\
\lvert \overline{v}_{i}(t) - e^{-\lambda_{i}(\tau - t)} v_{i\tau}\rvert
&\leq M_{1} e^{\lambda_{i}t} \int_{t}^{\tau} e^{-\lambda_{i}s} Q_{i} ds
\qquad (i=1,2).
\end{aligned}
\end{equation*}
We substitute \eqref{eq_87iu45kfe} into these inequalities, integrate, and simplify. We do not present the details of this calculation here; instead, we state the final result obtained using symbolic computations in Maple. A similar detailed hand calculation can be found in \cite{BakraniPhDthesis} (see, for example, Section 3.3.2). For some $M_{2} > 0$, we have
\begin{subequations}\label{eq_098ukihe_uuu}
\begin{align}
&\begin{aligned}
&\lvert\overline{u}_{1}(t) - e^{-\lambda_{1}t} u_{10}\rvert\\
&\quad\leq M_{2}\Bigl(e^{\left(-\lambda_{1}-2 \lambda_{2}\right) \tau +\lambda_{1} t} \lvert u_{20}^{2} v_{1\tau} v_{2\tau}^{2}\rvert
   + e^{-\lambda_{1} t +\left(-\lambda_{1}-\lambda_{2}\right) \tau} \lvert u_{20}^{2} v_{1\tau} v_{2\tau}\rvert\\
&\qquad {}+ e^{-3 \lambda_{1} \tau +\left(3 \lambda_{1}-2 \lambda_{2}\right) t} \lvert u_{20}^{2} v_{1\tau}^{3}\rvert
   + t e^{-\lambda_{1} t -\lambda_{2} \tau} \lvert u_{10} u_{20} v_{2\tau}\rvert\\
&\qquad {}+ t e^{-\lambda_{1} t -\lambda_{1} \tau} \lvert u_{10}^{2} v_{1\tau}\rvert
   + e^{-\lambda_{1} t -3 \lambda_{1} \tau} \lvert u_{20}^{3} v_{1\tau}^{2}\rvert\\
&\qquad {}+ e^{-\lambda_{1} t -2 \lambda_{1} \tau} \lvert u_{10} u_{20} v_{1\tau}^{2}\rvert
   + e^{-2 \lambda_{2} \tau +\left(-3 \lambda_{1}+2 \lambda_{2}\right) t} \lvert u_{10}^{3} v_{2\tau}\rvert\\
&\qquad {}+ e^{-\lambda_{1} t -2 \lambda_{1} \tau} \lvert u_{10} u_{20}^{2} v_{1\tau}\rvert
   + e^{-2 \lambda_{2} \tau +\left(-3 \lambda_{1}+2 \lambda_{2}\right) t} \lvert u_{10}^{2} v_{2\tau}^{2}\rvert\Bigr),
\end{aligned}\label{eq_098ukihe_uuu:u1}
\\[1ex]
&\begin{aligned}
&\lvert\overline{u}_{2}(t) - e^{-\lambda_{2}t} u_{20}\rvert\\
&\quad\leq M_{2}\Bigl(e^{\left(-\lambda_{1}-\lambda_{2}\right) \tau +\left(-2 \lambda_{1}+\lambda_{2}\right) t} \lvert u_{10}^{2} v_{1\tau} v_{2\tau} \rvert
   + t e^{-\lambda_{2} (t + \tau)} \lvert u_{20}^{2} v_{2\tau} \rvert\\
&\qquad {}+ e^{-2 \lambda_{2} \tau +\left(-2 \lambda_{1}+\lambda_{2}\right) t} \lvert u_{10} u_{20} v_{2\tau}^{2} \rvert
   + e^{\left(-2 \lambda_{1}-\lambda_{2}\right) \tau +\lambda_{2} t} \lvert u_{10}^{2} v_{1\tau}^{2} v_{2\tau} \rvert\\
&\qquad {}+ e^{-\lambda_{2} t} \lvert u_{20}^{3} \rvert
   + e^{-\lambda_{2} t} \lvert u_{10}^{2} u_{20} \rvert\\
&\qquad {}+ e^{-\lambda_{2} t -\lambda_{2} \tau} \lvert u_{10}^{3} v_{2\tau} \rvert
   + e^{-\lambda_{2} t -2 \lambda_{1} \tau} \lvert u_{20}^{2} v_{1\tau}^{2} \rvert\\
&\qquad {}+ e^{-3 \lambda_{2} \tau +\left(-4 \lambda_{1}+3 \lambda_{2}\right) t} \lvert u_{10}^{2} v_{2\tau}^{3} \rvert
   + t e^{-\lambda_{2} t -\lambda_{1} \tau} \lvert u_{10} u_{20} v_{1\tau} \rvert\Bigr),
\end{aligned}\label{eq_098ukihe_uuu:u2}
\\[1ex]
&\begin{aligned}
&\lvert \overline{v}_{1}(t) - e^{-\lambda_{1}(\tau - t)} v_{1\tau}\rvert\\
&\quad\leq M_{2}\Bigl(e^{-2 \lambda_{2} \tau +\left(-3 \lambda_{1}+2 \lambda_{2}\right) t} \lvert u_{10}^{3} v_{2\tau}^{2} \rvert
   + e^{-3 \lambda_{1} \tau +\left(3 \lambda_{1}-2 \lambda_{2}\right) t} \lvert u_{20}^{2} v_{1\tau}^{2} \rvert\\
&\qquad {}+ e^{-3 \lambda_{1} \tau +\left(3 \lambda_{1}-2 \lambda_{2}\right) t} \lvert u_{20} v_{1\tau}^{3} \rvert
   + e^{-\lambda_{1} t -2 \lambda_{2} \tau} \lvert u_{10} u_{20}^{2} v_{2\tau}^{2} \rvert\\
&\qquad {}+ e^{\lambda_{1} t -4 \lambda_{1} \tau} \lvert u_{10}^{2} v_{2\tau}^{3} \rvert
   + (\tau - t) e^{\left(-\lambda_{1}-\lambda_{2}\right) \tau +\lambda_{1} t} \lvert u_{20} v_{1\tau} v_{2\tau} \rvert\\
&\qquad {}+ (\tau - t) e^{\lambda_{1} t -2 \lambda_{1} \tau} \lvert u_{10} v_{1\tau}^{2} \rvert
   + e^{\left(-2 \lambda_{1}-\lambda_{2}\right) \tau +\lambda_{1} t} \lvert u_{10} u_{20} v_{2\tau}^{2} \rvert\\
&\qquad {}+ e^{\lambda_{1} t -3 \lambda_{1} \tau} \lvert u_{10} v_{1\tau} v_{2\tau}^{2} \rvert
   + e^{\lambda_{1} t -3 \lambda_{1} \tau} \lvert u_{10}^{2} v_{1\tau} v_{2\tau} \rvert\Bigr),
\end{aligned}\label{eq_098ukihe_uuu:v1}
\\[1ex]
&\begin{aligned}
&\lvert \overline{v}_{2}(t) - e^{-\lambda_{2}(\tau - t)} v_{2\tau}\rvert\\
&\quad\leq M_{2}\Bigl(e^{-3 \lambda_{1} \tau +\left(2 \lambda_{1}-\lambda_{2}\right) t} \lvert u_{10} u_{20} v_{1\tau}^{2} \rvert
   + e^{\lambda_{2} t -2 \lambda_{2} \tau} \lvert u_{20} v_{1\tau}^{3} \rvert\\
&\qquad {}+ e^{\left(-2 \lambda_{1}-\lambda_{2}\right) \tau +\left(2 \lambda_{1}-\lambda_{2}\right) t} \lvert u_{20}^{2} v_{1\tau} v_{2\tau} \rvert
   + (\tau - t) e^{\lambda_{2} t -2 \lambda_{2} \tau} \lvert u_{20} v_{2\tau}^{2} \rvert\\
&\qquad {}+ (\tau - t) e^{\left(-\lambda_{1}-\lambda_{2}\right) \tau +\lambda_{2} t} \lvert u_{10} v_{1\tau} v_{2\tau} \rvert
   + e^{\left(-2 \lambda_{1}-\lambda_{2}\right) \tau +\lambda_{2} t} \lvert u_{10}^{2} v_{2\tau}^{2} \rvert\\
&\qquad {}+ e^{-\lambda_{2} t -2 \lambda_{1} \tau} \lvert u_{10}^{2} u_{20} v_{1\tau}^{2} \rvert
   + e^{\lambda_{2} t -\lambda_{2} \tau} \lvert v_{1\tau}^{2} v_{2\tau} \rvert\\
&\qquad {}+ e^{-4 \lambda_{1} \tau +\left(4 \lambda_{1}-3 \lambda_{2}\right) t} \lvert u_{20}^{3} v_{1\tau}^{2} \rvert
   + e^{\lambda_{2} t -\lambda_{2} \tau} \lvert v_{2\tau}^{3} \rvert\Bigr).
\end{aligned}\label{eq_098ukihe_uuu:v2}
\end{align}
\end{subequations}
These estimates hold for any sufficiently small $\delta >0$, and so, we can make $\mathcal{A}$ invariant by choosing $\delta >0$ small enough from the beginning.

In the course of the proof of the invariance of $\mathcal{A}$, we have also shown that the image of any element of $\mathcal{A}$ under $\mathfrak{T}$ satisfies \eqref{eq_098ukihe_uuu:u1}--\eqref{eq_098ukihe_uuu:v2}. However, $\mathfrak{T}$ is a contraction on the complete metric space $\mathcal{I}$ (see, e.g. \cite{Dimabook}). Therefore, by the Banach fixed point theorem, the solution of system \eqref{eq:normal_form} that satisfies boundary conditions \eqref{eq_boundary_conditions_1} lies in $\mathcal{A}$, too, and since it is the fixed point of $\mathfrak{T}$, it satisfies \eqref{eq_098ukihe_uuu:u1}--\eqref{eq_098ukihe_uuu:v2}. On the other hand, all the terms in the right side of the first (resp. second) inequality that have a coefficient $u_{10}$ (resp. $u_{20}$) are bounded above by $e^{-\lambda_{1} t} \lvert u_{10}\rvert \delta^{2}$ (resp. $e^{-\lambda_{2} t} \lvert u_{20}\rvert \delta^{2}$). Moreover, the terms in the right side of the third (resp. fourth) inequality that have a coefficient $v_{1\tau}$ (resp. $v_{2\tau}$) are bounded above by $e^{-\lambda_{1} (\tau - t)} \lvert v_{1\tau}\rvert \delta^{2}$ (resp. $e^{-\lambda_{2} (\tau - t)} \lvert v_{2\tau}\rvert \delta^{2}$). Taking these two into account, we obtain estimates \eqref{eq:estimate_trajectory_1} and \eqref{eq:estimate_trajectory_2} from \eqref{eq_098ukihe_uuu:u1}--\eqref{eq_098ukihe_uuu:v2} by choosing an appropriate $M > 0$. This ends the proof of Lemma \ref{lem:formula_trajectories}.
\end{proof}

In the rest, we solve some boundary value problems to prove Lemmas \ref{Lem_620yu3tyohk}, \ref{Lem_derivatives_resp_tau_12}, \ref{Lem_78iuhr4} and \ref{Lem_derivatives_resp_tau_21}. Some of the calculations here are performed using symbolic computation in Maple. We do not include the details of these calculations; however, a similar detailed hand calculation can be found in \cite{BakraniPhDthesis} (see the proof of Lemma 3.34 therein).

Consider the integral operator $\mathfrak{T}$ given by \eqref{eq43bt4re}. The normal form \eqref{eq:normal_form} gives the following derivative estimates.
\begin{align*}
\shortintertext{\textit{Derivatives of $F_1$.}}
\partial F_{1}/\partial u_{1} &= \mathcal{O}\left(u_{1}v_{1}, u_{2}v_{2}, u_{2}^{2} v_{1}^{2}, u_{1}^{2}v_{2}^{2}\right),\\
\partial F_{1}/\partial u_{2} &= \mathcal{O}\left(u_{1} v_{2}, u_{2}^{2} v_{1} v_{2}, u_{2} v_{1} v_{2}^{2}, u_{1}^{2} u_{2} v_{1}, u_{1} u_{2} v_{1}^{2}, u_{2} v_{1}^{3}\right),\\
\partial F_{1}/\partial v_{1} &= \mathcal{O}\left(u_{1}^{2}, u_{1}u_{2}v_{1}v_{2}, u_{2}^{3}v_{2}, u_{2}^{2}v_{2}^{2}, u_{2}^{2} v_{1}^{2}, u_{1} u_{2}^{2} v_{1}\right),\\
\partial F_{1}/\partial v_{2} &= \mathcal{O}\left(u_{1}^{2}v_{1}, u_{1}u_{2}, u_{2}^{3}v_{1}, u_{2}^{2}v_{1}v_{2}, u_{1}^{3}v_{2}\right).\\
\shortintertext{\textit{Derivatives of $F_2$.}}
\partial F_{2}/\partial u_{1} &= \mathcal{O}\left(u_{1}u_{2}, u_{2}v_{1}, u_{1}^{3} v_{2}, u_{1}^{2} v_{1} v_{2}, u_{1} v_{1}^{2} v_{2}\right),\\
\partial F_{2}/\partial u_{2} &= \mathcal{O}\left(u_{1}^{2}, u_{1}v_{1}, u_{2}^{2}, u_{2}v_{2}\right),\\
\partial F_{2}/\partial v_{1} &= \mathcal{O}\left(u_{1}u_{2}, u_{1}^{3}v_{2}, u_{1}^{2}v_{1}v_{2}, u_{2}^{3}v_{1}, u_{2}^{2}v_{1}v_{2}\right),\\
\partial F_{2}/\partial v_{2} &= \mathcal{O}\left(u_{2}^{2}, u_{1}^{4}, u_{1}^{3} v_{1}, u_{1}^{2} u_{2} v_{2}, u_{1}^{2} v_{1}^{2}, u_{1}u_{2}v_{1}v_{2}\right).\\
\shortintertext{\textit{Derivatives of $G_1$.}}
\partial G_{1}/\partial u_{1} &= \mathcal{O}\left(v_{1}^{2}, v_{1}v_{2}u_{1}u_{2}, v_{2}^{3} u_{2}, v_{2}^{2} u_{2}^{2}, v_{2}^{2} u_{1}^{2}, v_{1} v_{2}^{2} u_{1}\right),\\
\partial G_{1}/\partial u_{2} &= \mathcal{O}\left(v_{1}^{2} u_{1}, v_{1} v_{2}, v_{2}^{3} u_{1}, v_{2}^{2} u_{1} u_{2}, v_{1}^{3}u_{2}\right),\\
\partial G_{1}/\partial v_{1} &= \mathcal{O}\left(v_{1} u_{1}, v_{2} u_{2}, v_{2}^{2} u_{1}^{2}, v_{1}^{2} u_{2}^{2}\right),\\
\partial G_{1}/\partial v_{2} &= \mathcal{O}\left(v_{1} u_{2}, v_{2}^{2} u_{1} u_{2}, v_{2} u_{1} u_{2}^{2}, v_{1}^{2} v_{2} u_{1}, v_{1} v_{2} u_{1}^{2}, v_{2} u_{1}^{3}\right).\\
\shortintertext{\textit{Derivatives of $G_2$.}}
\partial G_{2}/\partial u_{1} &= \mathcal{O}\left(v_{1}v_{2}, v_{1}^{3} u_{2}, v_{1}^{2} u_{1} u_{2}, v_{2}^{3} u_{1}, v_{2}^{2} u_{1} u_{2}\right),\\
\partial G_{2}/\partial u_{2} &= \mathcal{O}\left(v_{2}^{2}, v_{1}^{4}, v_{1}^{3} u_{1}, v_{1}^{2} v_{2} u_{2}, v_{1}^{2} u_{1}^{2}, v_{1} v_{2} u_{1} u_{2}\right),\\
\partial G_{2}/\partial v_{1} &= \mathcal{O}\left(v_{1} v_{2}, v_{2}u_{1}, v_{1}^{3} u_{2}, v_{1}^{2} u_{1} u_{2}, v_{1} u_{1}^{2} u_{2}\right),\\
\partial G_{2}/\partial v_{2} &= \mathcal{O}\left(v_{1}^{2}, v_{1} u_{1}, v_{2}^{2}, v_{2} u_{2}\right).
\end{align*}
Here $\mathcal{O}(m_{1}, \ldots, m_{k})$ stands for a linear combination of monomials $m_{1}$,..., $m_{k}$ for which the coefficients are ${C}^{\infty}$ functions of $(u_{1}, u_{2}, v_{1}, v_{2})$. Following the proof of Lemma \ref{lem:formula_trajectories}, the solution $x = (u_{1}(t), u_{2}(t), v_{1}(t), v_{2}(t))$ of system \eqref{eq:normal_form} to the boundary conditions \eqref{eq_boundary_conditions_1} satisfies the estimates in \eqref{eq67uyg2hoeu}. Recall \eqref{eq_87iu45kfe}, and define the corresponding majorants by
\begin{align*}
\shortintertext{\textit{Majorants for $F_1$.}}
A_{1u_{1}} &:= {\eta}_{1}{\nu}_{1} + {\eta}_{2}{\nu}_{2} + {\eta}_{2}^{2} {\nu}_{1}^{2} + {\eta}_{1}^{2}{\nu}_{2}^{2},\\
A_{1u_{2}} &:= {\eta}_{1} {\nu}_{2} + \eta_{2}^{2} {\nu}_{1} {\nu}_{2} + \eta_{2} {\nu}_{1} {\nu}_{2}^{2} + {\eta}_{1}^{2} {\eta}_{2} {\nu}_{1} + {\eta}_{1} {\eta}_{2} {\nu}_{1}^{2} + {\eta}_{2} {\nu}_{1}^{3},\\
A_{1v_{1}} &:= {\eta}_{1}^{2} + {\eta}_{2}^{3} {\nu}_{2} + {\eta}_{2}^{2} {\nu}_{2}^{2} + {\eta}_{2}^{2} {\nu}_{1}^{2} + {\eta}_{1} {\eta}_{2}^{2} {\nu}_{1},\\
A_{1v_{2}} &:= {\eta}_{1}^{2}{\nu}_{1} + {\eta}_{1}{\eta}_{2} + {\eta}_{2}^{3}{\nu}_{1} + {\eta}_{2}^{2}{\nu}_{1} {\nu}_{2} + {\eta}_{1}^{3}{\nu}_{2}.\\
\shortintertext{\textit{Majorants for $F_2$.}}
A_{2u_{1}} &:= {\eta}_{1}{\eta}_{2} + {\eta}_{2}{\nu}_{1} + {\eta}_{1}^{3} {\nu}_{2} + {\eta}_{1}^{2} {\nu}_{1} {\nu}_{2} + {\eta}_{1} {\nu}_{1}^{2} {\nu}_{2},\\
A_{2u_{2}} &:= {\eta}_{1}^{2} + {\eta}_{1}{\nu}_{1} + {\eta}_{2}^{2} + {\eta}_{2}{\nu}_{2},\\
A_{2v_{1}} &:= {\eta}_{1}{\eta}_{2} + {\eta}_{1}^{3}{\nu}_{2} + {\eta}_{1}^{2}{\nu}_{1}{\nu}_{2} + {\eta}_{2}^{3}{\nu}_{1} + {\eta}_{2}^{2}{\nu}_{1}{\nu}_{2},\\
A_{2v_{2}} &:= {\eta}_{2}^{2} + {\eta}_{1}^{4} + {\eta}_{1}^{3} {\nu}_{1} + {\eta}_{1}^{2} {\eta}_{2} {\nu}_{2} + {\eta}_{1}^{2} {\nu}_{1}^{2} + {\eta}_{1}{\eta}_{2}{\nu}_{1}{\nu}_{2}.\\
\shortintertext{\textit{Majorants for $G_1$.}}
B_{1u_1} &:= {\nu}_{1}^{2} + {\nu}_{2}^{3} {\eta}_{2} + {\nu}_{2}^{2} {\eta}_{2}^{2} + {\nu}_{2}^{2} {\eta}_{1}^{2} + {\nu}_{1} {\nu}_{2}^{2} {\eta}_{1},\\
B_{1u_2} &:= {\nu}_{1}^{2} {\eta}_{1} + {\nu}_{1} {\nu}_{2} + {\nu}_{2}^{3} {\eta}_{1} + {\nu}_{2}^{2} {\eta}_{1} {\eta}_{2} + {\nu}_{1}^{3}{\eta}_{2},\\
B_{1v_1} &:= {\nu}_{1} {\eta}_{1} + {\nu}_{2} {\eta}_{2} + {\nu}_{2}^{2} {\eta}_{1}^{2} + {\nu}_{1}^{2} {\eta}_{2}^{2},\\
B_{1v_2} &:= {\nu}_{1} {\eta}_{2} + {\nu}_{2}^{2} {\eta}_{1} {\eta}_{2} + {\nu}_{2} {\eta}_{1} {\eta}_{2}^{2} + {\nu}_{1}^{2} {\nu}_{2} {\eta}_{1} + {\nu}_{1} {\nu}_{2} {\eta}_{1}^{2} + {\nu}_{2} {\eta}_{1}^{3}.\\
\shortintertext{\textit{Majorants for $G_2$.}}
B_{2u_1} &:= {\nu}_{1}{\nu}_{2} + {\nu}_{1}^{3} {\eta}_{2} + {\nu}_{1}^{2} {\eta}_{1} {\eta}_{2} + {\nu}_{2}^{3} {\eta}_{1} + {\nu}_{2}^{2} {\eta}_{1} {\eta}_{2},\\
B_{2u_2} &:= {\nu}_{2}^{2} + {\nu}_{1}^{4} + {\nu}_{1}^{3} {\eta}_{1} + {\nu}_{1}^{2} {\nu}_{2} {\eta}_{2} + {\nu}_{1}^{2} {\eta}_{1}^{2} + {\nu}_{1} {\nu}_{2} {\eta}_{1} {\eta}_{2},\\
B_{2v_1} &:= {\nu}_{1} {\nu}_{2} + {\nu}_{2}{\eta}_{1} + {\nu}_{1}^{3} {\eta}_{2} + {\nu}_{1}^{2} {\eta}_{1} {\eta}_{2} + {\nu}_{1} {\eta}_{1}^{2} {\eta}_{2},\\
B_{2v_2} &:= {\nu}_{1}^{2} + {\nu}_{1} {\eta}_{1} + {\nu}_{2}^{2} + {\nu}_{2} {\eta}_{2}.
\end{align*}
There exists a constant $M_{1} > 0$ such that for $i,j = 1,2$, we have
\begin{equation*}
\begin{alignedat}{2}
\left|\frac{\partial F_i}{\partial u_j}(x)\right|
&\leq M_1 A_{iu_j},
\qquad
&\left|\frac{\partial F_i}{\partial v_j}(x)\right|
&\leq M_1 A_{iv_j},
\\[0.8ex]
\left|\frac{\partial G_i}{\partial u_j}(x)\right|
&\leq M_1 B_{iu_j},
\qquad
&\left|\frac{\partial G_i}{\partial v_j}(x)\right|
&\leq M_1 B_{iv_j}.
\end{alignedat}
\end{equation*}
Here we have used the fact that, in a sufficiently small neighborhood of the origin, $2\lvert u_{1}u_{2}v_{1}v_{2}\rvert \leq u_{1}^{2}+u_{2}^{2}v_{2}^{2}$. Thus, the contribution of this monomial to $\partial F_{1}/\partial v_{1}$ is already covered by the displayed majorant. The analogous contribution to $\partial G_{1}/\partial u_{1}$ is treated in the same way. Similar inequalities allow further simplification of the majorants, but the reduction in the number of terms is modest, so we retain the expressions above.

Consequently, by \eqref{eq43bt4re}, if $\lvert U_{i}(t)\rvert \leq \Theta_{i}(t)$ and $\lvert V_{i}(t)\rvert \leq \Upsilon_{i}(t)$ for some functions $\Theta_{i}$ and $\Upsilon_{i}$ and for $i=1,2$, then
\begin{equation}\label{eq_iyg4kwblk}
\begin{aligned}
\lvert\overline{U}_{i}(t) - e^{-\lambda_{i} t} U_{i0}\rvert &\leq M_{1}\int_{0}^{t} e^{\lambda_{i}(s-t)} \sum_{j=1}^{2} \Big(A_{iu_{j}} \Theta_{j}(s) + A_{iv_{j}} \Upsilon_{j}(s)\Big) ds,\\
\lvert\overline{V}_{i}(t) - e^{-\lambda_{i} (\tau -t)} V_{i\tau}\rvert &\leq M_{1}\int_{t}^{\tau} e^{\lambda_{i}(t-s)} \sum_{j=1}^{2} \Big(B_{iu_{j}} \Theta_{j}(s) + B_{iv_{j}} \Upsilon_{j}(s)\Big) ds.
\end{aligned}
\end{equation}

Our strategy for proving Lemmas \ref{Lem_620yu3tyohk}, \ref{Lem_derivatives_resp_tau_12}, \ref{Lem_78iuhr4} and \ref{Lem_derivatives_resp_tau_21} is the same as in the proof of Lemma \ref{lem:formula_trajectories}. For each lemma, we introduce a set
\begin{equation}\label{eq:set_A}
\mathcal{A} := \left\{\left(U_{1}(t), U_{2}(t), V_{1}(t), V_{2}(t)\right):\quad
\lvert U_{i}(t)\rvert \leq \Theta_{i}(t),\quad
\lvert V_{i}(t)\rvert \leq \Upsilon_{i}(t) \quad\mathit{for}\,\, i=1,2\right\},
\end{equation}
and using Maple, we estimate symbolically the image of any element of $\mathcal{A}$ under $\mathfrak{T}$, and then by virtue of \eqref{eq_iyg4kwblk}, we conclude that $\mathcal{A}$ is invariant under $\mathfrak{T}$. However, since $\mathfrak{T}$ is a contraction (see, e.g. \cite{Dimabook}), the fixed point of $\mathfrak{T}$, which is indeed the solution to the boundary value problem \eqref{eq:boundary_conditions_derivative}, belongs to $\mathcal{A}$, and so, it satisfies the estimates that we have found for the image of elements in $\mathcal{A}$ under $\mathfrak{T}$. This is how we prove the lemmas in this section.

\begin{proof}[Proof of Lemmas \ref{Lem_620yu3tyohk} and \ref{Lem_78iuhr4}]
Consider the set $\mathcal{A}$ given by \eqref{eq:set_A}, where $\Theta_{i}(t) = 2 e^{-\lambda_{i} t}$ and $\Upsilon_{i}(t) = 2 e^{-\lambda_{i}(\tau - t)}$. Using Maple, we apply \eqref{eq_iyg4kwblk} and obtain
\begin{equation*}
\begin{aligned}
\lvert \overline{U}_{1}(t) - e^{-\lambda_{1} t} U_{10}\rvert &\leq M\Big(
t e^{-\lambda_{1} t -\lambda_{1} \tau} \delta^{2} + e^{-\lambda_{1} t -2 \lambda_{1} \tau} \delta^{3} + e^{-3 \lambda_{1} \tau +\left(3 \lambda_{1}-2 \lambda_{2}\right) t} \delta^{4}\Big),\\
\lvert \overline{U}_{2}(t) - e^{-\lambda_{2} t} U_{20}\rvert &\leq
M e^{-\lambda_{2} t} \delta^{2},\\
\lvert \overline{V}_{1}(t) - e^{-\lambda_{1}(\tau - t)} V_{1\tau}\rvert &\leq
M \Big((\tau - t) e^{\lambda_{1} t -2 \lambda_{1} \tau} \delta^{2} + e^{\lambda_{1} t -3 \lambda_{1} \tau} \delta^{3} + e^{-2 \lambda_{2} \tau +\left(-3 \lambda_{1}+2 \lambda_{2}\right) t} \delta^{4}\Big),\\
\lvert \overline{V}_{2}(t) - e^{-\lambda_{2}(\tau - t)} V_{2\tau}\rvert &\leq M e^{-\lambda_{2}(\tau - t)} \delta^{2}
\end{aligned}
\end{equation*}
for some $M>0$, which holds for any sufficiently small $\delta >0$. Choosing $\delta >0$ sufficiently small from the beginning, the set $\mathcal{A}$ becomes invariant under $\mathfrak{T}$, and so the solution to the boundary value problem \eqref{eq:boundary_conditions_derivative} satisfies
\begin{equation*}
\lvert U_{i}(t) - e^{-\lambda_{i} t} U_{i0}\rvert \leq M e^{-\lambda_{2} t} \delta^{2}
\qquad \mathrm{and} \qquad
\lvert V_{i}(t) - e^{-\lambda_{i}(\tau - t)} V_{i\tau}\rvert \leq M e^{-\lambda_{2} (\tau - t)} \delta^{2}.
\end{equation*}
Lemmas \ref{Lem_620yu3tyohk} and \ref{Lem_78iuhr4} follow from this estimate by substituting the corresponding values of $U_{10}$, $U_{20}$, $V_{1\tau}$, and $V_{2\tau}$. This ends the proof.
\end{proof}

To prove Lemmas \ref{Lem_derivatives_resp_tau_12} and \ref{Lem_derivatives_resp_tau_21}, we first state the following lemma.
\begin{mylem}
Consider the set
\begin{equation}\label{eq78t7yqgdeft}
\mathcal{A} := \left\{\left(U_{1}(t), U_{2}(t), V_{1}(t), V_{2}(t)\right):\quad
\lvert U_{i}(t)\rvert \leq e^{-\lambda_{i}t} \delta,\quad
\lvert V_{i}(t)\rvert \leq \left(\lambda_{i} + 1\right) e^{-\lambda_{i}(\tau - t)} \lvert v_{i\tau}\rvert\right\},
\end{equation}
and suppose $U_{10} = U_{20} = 0$. There exists $M >0$ such that for any sufficiently small $\delta >0$, we have
\begin{subequations}\label{eq687iughwbyt4}
\begin{align}
&\begin{aligned}
\lvert \overline{U}_{1}(t)\rvert
&\leq M\Bigl(t e^{-\lambda_{1} t -\lambda_{1} \tau} \delta^{3} + e^{-\lambda_{1} t -2 \lambda_{1} \tau} \delta^{4}\\
&\qquad {}+ e^{\left(-\lambda_{1}-\lambda_{2}\right) \tau +\left(-\lambda_{1}+\lambda_{2}\right) t} \delta^{4}
   + e^{-\lambda_{1} t -\lambda_{1} \tau} \delta^{5}\Bigr),
\end{aligned}\label{eq687iughwbyt4:U1}
\\[1ex]
&\lvert \overline{U}_{2}(t)\rvert\leq M e^{-\lambda_{2} t} \delta^{3},\label{eq687iughwbyt4:U2}
\\[1ex]
&\begin{aligned}
&\lvert \overline{V}_{1}(t) - e^{-\lambda_{1}(\tau - t)} V_{1\tau}\rvert\\
&\quad\leq M\Bigl(e^{-\lambda_{1} (\tau -t)} \delta^{2} \lvert v_{1\tau}\rvert
   + e^{\left(-2 \lambda_{1}-\lambda_{2}\right) \tau +\lambda_{1} t} \delta \lvert u_{20} v_{2\tau}^{3}\rvert\\
&\qquad {}+ e^{-2 \lambda_{2} \tau +\left(-3 \lambda_{1}+2 \lambda_{2}\right) t} \delta \lvert u_{10}^{2} v_{2\tau}^{2}\rvert
   + e^{-\lambda_{1} t -2 \lambda_{2} \tau} \delta \lvert u_{10} u_{20} v_{2\tau}^{2}\rvert\\
&\qquad {}+ e^{\left(-2 \lambda_{1}-\lambda_{2}\right) \tau +\lambda_{1} t} \delta \lvert u_{10} v_{2\tau}^{2}\rvert
   + e^{-\lambda_{1} t -2 \lambda_{2} \tau} \delta \lvert u_{20}^{2} v_{2\tau}^{2}\rvert\Bigr),
\end{aligned}\label{eq687iughwbyt4:V1}
\\[1ex]
&\begin{aligned}
&\lvert \overline{V}_{2}(t) - e^{-\lambda_{2}(\tau - t)} V_{2\tau}\rvert\\
&\quad\leq M\Bigl(e^{-\lambda_{2}(\tau -t)} \delta^{2} \lvert v_{2\tau}\rvert
   + e^{-\lambda_{2} t -2 \lambda_{1} \tau} \delta \lvert u_{10} u_{20} v_{1\tau}^{2}\rvert\\
&\qquad {}+ e^{\lambda_{2} (t -2\tau)} \delta \lvert v_{1\tau}^{4}\rvert
   + e^{\lambda_{2}(t -2\tau)} \lvert u_{20} v_{1\tau}^{3}\rvert\\
&\qquad {}+ e^{-3 \lambda_{1} \tau +\left(2 \lambda_{1}-\lambda_{2}\right) t} \delta \lvert u_{20} v_{1\tau}^{2}\rvert
   + e^{-\lambda_{2} t -2 \lambda_{1} \tau} \delta \lvert u_{10}^{2} v_{1\tau}^{2}\rvert\\
&\qquad {}+ e^{-3 \lambda_{1} \tau +\left(2 \lambda_{1}-\lambda_{2}\right) t} \delta \lvert u_{10} v_{1\tau}^{3}\rvert\Bigr).
\end{aligned}\label{eq687iughwbyt4:V2}
\end{align}
\end{subequations}
\end{mylem}

Setting $\Theta_{i}(t) = e^{-\lambda_{i}t} \delta$ and $\Upsilon_{i}(t) = \left(\lambda_{i} + 1\right) e^{-\lambda_{i}(\tau - t)} \lvert v_{i\tau}\rvert$ and noting that $U_{10} = U_{20} = 0$, this lemma is directly obtained from \eqref{eq_iyg4kwblk} using symbolic calculation with Maple.

\begin{proof}[Proof of Lemma \ref{Lem_derivatives_resp_tau_12}]
We claim that under the conditions $u_{10} = \sigma_{0}\delta$, $v_{2\tau} = \sigma_{1} \delta$ and $(u_{20}, v_{20})\in Y$, the set $\mathcal{A}$ is invariant under $\mathfrak{T}$. Recall the boundary conditions \eqref{eq_098iuh}. It readily follows from \eqref{eq687iughwbyt4} that if $\delta >0$ is small enough, then $\overline{U}_1$, $\overline{U}_2$, and $\overline{V}_2$ satisfy the estimates in \eqref{eq78t7yqgdeft}. We now show that $\overline{V}_1$ satisfies the estimate in \eqref{eq78t7yqgdeft}, too. By \eqref{eq687iughwbyt4}, we have
\begin{equation*}
\begin{aligned}
&\lvert \overline{V}_{1}(t) - e^{-\lambda_{1}(\tau - t)} V_{1\tau}\rvert\\
&\qquad\qquad \leq M \Big(e^{-\lambda_{1}(\tau - t)} \lvert v_{1\tau}\delta^{2}\rvert + 2 e^{\left(-2 \lambda_{1}-\lambda_{2}\right) \tau +\lambda_{1} t} \delta^{4} + 3 e^{-2 \lambda_{2} \tau +\left(-3 \lambda_{1}+2 \lambda_{2}\right) t} \delta^{5}\Big).
\end{aligned}
\end{equation*}
Recall \eqref{eq_67iby20iubhv}, i.e. $v_{1\tau} = \kappa \sigma_{0}\sigma_{1} \gamma^{-1} u_{20} \lvert v_{20}\rvert^{1-\gamma} \delta^{\gamma - 1}[1 + O(\delta^{2})]$. Taking \eqref{eq567ug79osw} into account, there exist positive $K$ and $K^{\prime}$ such that
\begin{equation*}
\lvert v_{1\tau}\rvert \geq K^{\prime} \lvert v_{20}\rvert^{2-\gamma} \delta^{\gamma -1} > K e^{-2\lambda_{2}\tau + \lambda_{1}\tau} \delta
\end{equation*}
holds for any sufficiently small $\delta >0$. This gives
\begin{equation*}
2 e^{\left(-2 \lambda_{1}-\lambda_{2}\right) \tau +\lambda_{1} t} \delta^{4} + 3 e^{-2 \lambda_{2} \tau +\left(-3 \lambda_{1}+2 \lambda_{2}\right) t} \delta^{5}
\leq 5K^{-1} e^{-\lambda_{1}(\tau - t)} \lvert v_{1\tau}\rvert \delta^{3}.
\end{equation*}
This implies $\lvert \overline{V}_{1}(t) - e^{-\lambda_{1}(\tau - t)} V_{1\tau}\rvert = e^{-\lambda_{1}(\tau - t)} O(\lvert v_{1\tau}\delta^{2}\rvert)$. Considering \eqref{eq_098iuh}, and choosing $\delta >0$ small enough, we obtain the estimate in \eqref{eq78t7yqgdeft} for $\overline{V}_1$, as desired. This shows the invariance of $\mathcal{A}$. Thus, the solution to the boundary conditions \eqref{eq_098iuh} belongs to $\mathcal{A}$, and satisfies \eqref{eq687iughwbyt4}. The estimates in Lemma \ref{Lem_derivatives_resp_tau_12} follow from \eqref{eq687iughwbyt4}.
\end{proof}

\begin{proof}[Proof of Lemma \ref{Lem_derivatives_resp_tau_21}]
The proof is quite similar to the proof of Lemma \ref{Lem_derivatives_resp_tau_12}. We show that under the conditions $u_{20} = \sigma_{1}\delta$, $v_{1\tau} = \sigma_{2} \delta$ and $(u_{10}, v_{10})\in Y$, the set $\mathcal{A}$ is invariant under $\mathfrak{T}$. Considering boundary conditions \eqref{eq_1359y7hirg}, it readily follows from \eqref{eq687iughwbyt4} that if $\delta >0$ is small enough, then $\overline{U}_1$, $\overline{U}_2$, and $\overline{V}_1$ satisfy the estimates in \eqref{eq78t7yqgdeft}. We now show that $\overline{V}_2$ satisfies the estimate in \eqref{eq78t7yqgdeft}, too. By \eqref{eq687iughwbyt4}, we have
\begin{equation}\label{eq:67uhufivr4hi35j}
\begin{aligned}
\lvert \overline{V}_{2}(t) - e^{-\lambda_{2}(\tau - t)} V_{2\tau}\rvert &\leq M \Big(e^{-\lambda_{2}(\tau - t)} \lvert v_{2\tau}\rvert\delta^{2} + 2e^{-3 \lambda_{1} \tau +\left(2 \lambda_{1}-\lambda_{2}\right) t} \delta^{4}\\
&\quad + 2 e^{\lambda_{2} t -2 \lambda_{2} \tau} \delta^{4} + 2e^{-\lambda_{2} t -2 \lambda_{1} \tau} \delta^{4}\Big).
\end{aligned}
\end{equation}
Since $(u_{10}, v_{10})\in Y$, by \eqref{eq_45tuvkb7taguer}, we have $v_{2\tau} = \kappa\sigma_{1} \sigma_{2} \gamma u_{10} |v_{10}|^{1-\frac{1}{\gamma}} \delta^{\frac{1}{\gamma} -1} [1 + O(\delta^{2})]$. 
Note also that since $(u_{10}, v_{10})\in Y$, we have $|u_{10}| \geq m^{-1} |v_{10}|$ by \eqref{eq:definition_regions}. Thus, there exists $K^{\prime} >0$ such that
\begin{equation*}
\lvert v_{2\tau}\rvert
\geq K^{\prime} \lvert v_{10}\rvert^{2-\frac{1}{\gamma}} \delta^{\frac{1}{\gamma} -1}
\end{equation*}
holds for any sufficiently small $\delta > 0$. By \eqref{eq024rfghvgyi} and \eqref{eq9er578iuytry}, there exists $K > 0$ such that
\begin{equation*}
\lvert v_{2\tau}\rvert
\geq K^{\prime} \lvert v_{10}\rvert^{2-\frac{1}{\gamma}} \delta^{\frac{1}{\gamma} -1}
\geq K e^{(\lambda_{2} - 2\lambda_{1})\tau} \delta
\end{equation*}
holds for any sufficiently small $\delta > 0$. Taking the right side of \eqref{eq:67uhufivr4hi35j} into account, we have
\begin{equation*}
\begin{aligned}
& M \Bigl(e^{-\lambda_{2}(\tau - t)} \lvert v_{2\tau}\rvert\delta^{2} + 2e^{-3 \lambda_{1} \tau +\left(2 \lambda_{1}-\lambda_{2}\right) t} \delta^{4} + 2 e^{\lambda_{2} t -2 \lambda_{2} \tau} \delta^{4} + 2e^{-\lambda_{2} t -2 \lambda_{1} \tau} \delta^{4}\Bigr)\\
&\qquad\qquad\qquad\qquad \leq M\left(1 + 6K^{-1}\delta\right) e^{-\lambda_{2}(\tau - t)} \lvert v_{2\tau}\rvert\delta^{2}.
\end{aligned}
\end{equation*}
This implies $\lvert \overline{V}_{2}(t) - e^{-\lambda_{2}(\tau - t)} V_{2\tau}\rvert = e^{-\lambda_{2}(\tau - t)} O(\lvert v_{2\tau}\delta^{2}\rvert)$. Considering \eqref{eq_1359y7hirg}, and choosing $\delta >0$ small enough, we obtain the estimate in \eqref{eq78t7yqgdeft} for $\overline{V}_2$, as desired. This shows the invariance of $\mathcal{A}$. Thus, the solution to the boundary conditions \eqref{eq_1359y7hirg} belongs to $\mathcal{A}$, and satisfies \eqref{eq687iughwbyt4}. The estimates in Lemma \ref{Lem_derivatives_resp_tau_21} follow from \eqref{eq687iughwbyt4}.
\end{proof}

{
\small
\bibliography{Myreferences}
\bibliographystyle{abbrvurl}
}

\end{document}